\documentclass[ejs]{imsart}

\RequirePackage[OT1]{fontenc}
\RequirePackage{amsthm,amsmath}
\RequirePackage[numbers]{natbib}
\RequirePackage[colorlinks,citecolor=blue,urlcolor=blue]{hyperref}
\usepackage{graphicx}

\usepackage{amsmath,amsfonts,amsthm,epsf,stackrel}
\usepackage{amssymb}
\usepackage{graphicx,subfigure}
\usepackage{float}
\usepackage{enumerate}
\usepackage{hyperref}
\usepackage{leftidx}
\usepackage{multirow}

\usepackage[mathscr]{euscript}
\usepackage{mathrsfs}

\usepackage{caption}

\usepackage{amsfonts}
\usepackage{upgreek}
\usepackage{contour}

\newtheorem{theorem}{Theorem}[section]
\newtheorem{lemma}[theorem]{Lemma}
\newtheorem{proposition}[theorem]{Proposition}

\newtheorem*{thma}{Lemma 2.3}

\newcommand{\tendd}{\stackrel{{\cal D}}{\longrightarrow}}
\newcommand{\bs}{\boldsymbol}
\newcommand{\ve}{\varepsilon}

\newcommand{\ignore}[1]{}
\DeclareMathAlphabet{\mathpzc}{OT1}{pzc}{m}{it}

\usepackage{xcolor}
\pubyear{2020}
\volume{0}
\issue{0}
\firstpage{1}
\lastpage{8}
\arxiv{2012.05949}

\startlocaldefs
\numberwithin{equation}{section}
\theoremstyle{plain}

\endlocaldefs

\begin{document}

\begin{frontmatter}
\title{Optimal selection of sample-size dependent common subsets of covariates for multi-task regression prediction}
\runtitle{Optimal selection}

\begin{aug}
\author{\fnms{David} \snm{Azriel}\ead[label=e1]{davidazr@technion.ac.il}}

\address{Faculty of Industrial Engineering and Management,\\
Technion – Israel Institute of Technology, Haifa, Israel\\
\printead{e1}}

\author{\fnms{Yosef} \snm{Rinott}\ead[label=e2]{yosef.rinott@mail.huji.ac.il}}

\address{Department of Statistics and
	Center for the Study of Rationality \\
The Hebrew University, Jerusalem, Israel\\
\printead{e2}}

\runauthor{Azriel and Rinott}

\affiliation{Technion and The Hebrew University}

\end{aug}

\begin{abstract}
	
An analyst 	is given a training set consisting of regression datasets $D_j$ of different sizes, which are distributed according to some $G_j$, $j=1,\ldots,\cal J$,  where the distributions $G_j$ are assumed to form a random sample generated by some common source. In particular, the $D_j$'s have a common set of covariates and they are all labeled.  
The training  set is used by the analyst for selection of subsets of covariates denoted by ${\mathpzc P}^*(n)$, whose role is described next.

 The multi-task problem we consider is as follows: given a number of random labeled datasets (which may be in the training set or not) $D_{J_k}$ of size $n_k$, $k=1,\ldots,K$, estimate separately for each dataset the regression coefficients on the subset of covariates ${\mathpzc P}^*(n_k)$ and then predict future dependent variables given their covariates.

 Naturally, a large sample size $n_k$ of $D_{J_k}$ allows a larger subset of covariates, and the dependence of the size of the selected covariate subsets on $n_k$ is needed in order to achieve good prediction and avoid overfitting.
Subset selection  is notoriously difficult and computationally demanding, and requires large samples; using  all the regression  datasets  in the training set together amounts to borrowing strength toward  better selection under suitable assumptions. Furthermore, using common subsets for all regressions having a given sample size standardizes and simplifies the data collection  and avoids having to select and use a different subset for each prediction task.
Our approach is efficient when  the relevant covariates for prediction are  common to the different regressions,  while the models' coefficients may vary between different regressions. 

Last but not least, we propose a simple and meaningful measure, GENO, that allows comparisons of the predictive value of different subsets of covariates by comparing the sample size they require in order to achieve the same prediction error.
\end{abstract}


\tableofcontents
\addcontentsline{toc}{chapter}{Appendix B: A table of notation}
\end{frontmatter}

\section{Introduction} \label{sec:Int}
\subsection{A general description of the problem} \label{sec:gendesc}

This paper concerns data consisting of a class of regression datasets, and a multi-task of predictions in different regressions.  The emphasis is on selection of common subsets of covariates for prediction in the different regressions, which depend on the regression datasets' sample sizes. As documented in the classical model selection literature, the size or dimension of models for prediction should depend on the sample size, for example through a penalty function that depends on the dimension of the model and the sample size.   References and further details will be provided following  a description of our motivating problem.

In order to improve service, a hospital wants to develop a tool for predicting the actual duration of planned  visits of any particular patient to any doctor in the hospital.
Given a sample of size $n$ of different patients'  visits  to any particular doctor ($n$ can vary between doctors) with covariates such as the past durations of the patients' visits, the nature of the visits, the time scheduled etc., and a response variable, which is the  actual  duration, the goal is to predict the duration of the next visit of a given patient to the particular doctor. Our objective is to select an optimal subset  of covariates, denoted by ${\mathpzc P}^*(n)$,  to be used for the prediction of a future visit's duration. We shall provide a procedure that selects the optimal set with high probability.  The number of covariates in the set ${\mathpzc P}^*(n)$ depends naturally on $n$, with a large $n$ allowing more variables in the regression, taking account of the need to find the right balance between efficiency of models, and the pitfall of overfitting.  For any given $n$, we want to select a standard  set of covariates to be used for any doctor in the hospital for whom we have a regression dataset of size $n$. 
However,  we allow different regression coefficients for different doctors since different doctors may be influenced differently by the patient's background. An intercept  for each doctor represents her or his general tendency for longer or shorter visits.

Standardization is desirable for more than one reason, and will be discussed in more detail later.  First, it obviously simplifies the data collection and maintenance. Second,  performing separate model selection for each doctor may be computationally demanding. Third, model selection is notorious for being difficult and to require much data. The idea we present is to perform model selection on the basis of a sample of doctors as described below, and thus borrow strength from different datasets  and obtain better subset selection. 

In order to perform the subset selection we assume we have a training sample of $\cal J$ doctors each  consisting of a dataset containing the covariates and the response (actual duration) of $N_j$ visits, $j=1,\ldots,\cal J$. Under certain assumptions, we use these data to select  subsets of covariates for different value of $n$. 
We then use the selected subset for prediction for any doctor (who in general may not be in the training sample) on the basis of a sample of visits (of some size $n$) as described above. When predicting for a doctor in the sample, say doctor $j$, it is natural to take $n=N_j$.  The  subset selection procedure and its properties are the focus of this paper.

One may suggest to concatenate the whole training sample and perform a single regression with the same coefficients for all doctors, but allowing a different intercept for each doctor. In certain cases this may result in good subset selection. However, suppose, for example, that for about half of the  doctors the covariate  ``duration of  previous visit"  has a positive coefficient in the regression and for the other half it is negative. It is easy to conceive of a justification for  each possibility. In this case, this covariate may not enter the model if the regression  is computed by concatenating the data into a single model. However, allowing different regression coefficients for different doctors, the variable may enter the model and contribute to the prediction with a different coefficient for different doctors. Thus, allowing  the regression coefficients to vary between individual regressions adds flexibility to the model, and in particular it improves the prediction in our dataset (verified by cross validation, see Section \ref{sec:comparison}). Of course an informal screening of variables is often done, either at the stage of collecting the data, or before conducting formal variable selection and analysis. In particular, researchers may decide to avoid certain variables or interaction terms in order to keep the selection process feasible.

The classical theory of model selection in regression deals with the selection of a subset of covariates (or features) that are useful for prediction based on a single regression dataset.  Numerous model selection methods have been suggested; AIC (\citet{Akaike}), Mallows $C_p$ (\citet{Mallows}), and FIC (\citet{FIC}) are prominent examples. These methods apply to a single regression dataset of a given size, for which a model is to be selected and then used for prediction.  For a well-known Bayesian approach to model selection, see  \citet{Schwarz}. A large body of literature emerged following these articles. In the setup of a single dataset, serious issues of optimality arise; see, e.g., \citet{Yang}.

Breiman's celebrated paper \citep{Breiman} starts with a similar training set of regression datasets with similar assumptions, however, both the subset of covariates and the regression coefficients used for prediction are common to all regressions.
A very closely related setup appears in \citet{Jordan}
which ``addresses the problem of recovering a common set of covariates that are relevant simultaneously to several classification problems." The paper focuses on  classification or discrimination problems, but regression is also mentioned. References cited in this paper, which deal with the same problem, are referred  to  as ``transfer learning" or ``multi-task learning" in the machine learning literature. They
demonstrate that learning multiple related tasks from data simultaneously can be advantageous in terms of predictive performance relative to learning these tasks independently.  In \citet{Jordan} the goal is to decide which variables are  ``relevant
to the overall class of prediction problems without making a commitment to a specific value of a parameter," that is, allowing different parameters for the different prediction tasks, and to
``borrow strength across multiple estimation problems in order to support a decision that a covariate is to be selected." A large number of papers and review articles on multi-task learning have appeared, mostly in the past decade. For a recent survey containing numerous applications and references, see, for example, \citet{Zhang}.   In the latter paper  transfer learning refers (in our setup) to predicting for a single target file that may not be in the training sample, while multi-task learning refers to predicting for every dataset in the training sample. Here we consider both possibilities.

Our paper  differs from  \citet{Jordan} and more generally from the multi-task literature in several ways: our emphasis is on regression, our focus is not on algorithms but rather on asymptotic consistency and optimality type results; however, the main difference is that in the spirit of model selection,   the  selected common covariate subsets,  depend on the sample sizes of the different regression prediction tasks, thus avoiding under and overfitting. 
The issue of different sample sizes $n$ appears in \citet{Zhang} in a context where face databases have different image sizes; however the proposed solution is to project these databases to a common subspace, which results in loss of information, rather than taking the task size into account as we propose.

\subsection{A formal setup}\label{sec:forset}
Our setup is formalized as follows. We assume we have a training sample ${\cal T}$ of regression datasets all having the same set of covariates. Thus ${\cal T} =\{D_j: j=1,\ldots,\cal J\}$  with $D_j=\{({\bf X}_{ij},Y_{ij})\},\,i=1,\ldots,N_j,\,\,j=1,\ldots,\cal J$, where ${\bf X}_{ij} \in \mathbb{R}^d$ is a column vector of $d$ random covariate values of the $i$th subject in the $j$th dataset, and $Y_{ij} \in \mathbb{R}$ is a response variable. For each $j$, the $N_j$ vectors  $({\bf X}_{ij},Y_{ij})$ are  iid from some distribution $G_j \in \cal G$, where $\cal G$ is a set (population) of distributions of size $|{\cal G}|=\cal K$. We assume that ${\cal J \le \cal K} \le \infty$, and that $\{G_j\}_{j=1}^{\cal J}$ is a random sample from $\cal G$. Now consider a new regression dataset of some size $n$,  $D_J=\{({\bf X}_{iJ},Y_{iJ})\}_{i=1}^n$ distributed according to  $G_J$, a random element of $\cal G$, which may but need not be in the training set $\cal T$.
If $D_J$ is in $\cal T$ then it is natural to assume that $n=N_J$. 

We consider the following task:
for $({\bf X}, Y) \sim G_J$ independent of the above datasets,  we want to predict $Y$ from a given ${\bf X}$ using the sample $D_J$. It is natural to be interested in the multi-task of prediction for many random $D_J$'s; however, it suffices to study the prediction error for one such $D_J$.
Since $G_J$ is random, we clearly need to consider different possible values of $n$, and  random covariates.  Our treatment of random covariates is based on a generalization of Mallows $C_p$ to random covariates that was inspired by notes generously given to us by Larry Brown (see \cite{Brown}).

As usual, the prediction model  involves two components, the subset of variables to be used, and their regression coefficients. The regression coefficients will be estimated by standard least squares based on the sample $D_J$, and thus will vary between $D_J$'s. For the subset selection, a task that is known to  require large samples,  we shall pool the whole training set. Such pooling can be efficient if the set of distributions $\cal G$, which may be finite or generated by a probability model (superpopulation model), is sufficiently homogeneous (to be discussed in Section \ref{sec:Jind})  in a way that justifies a common model selection. Besides some technical conditions for such homogeneity, a user would have to apply common sense to decide if one can borrow strength and learn the subset selection from the pooled sample $\cal T$ rather than from the individual dataset $D_J$. As mentioned before, numerous examples appear in \citet{Zhang} and the references therein. 
Our goal is to select for each possible value of $n$, a subset of covariates based on the pooled training sample $\cal T$  and use it for prediction, using least squares estimates,  computed for each regression dataset  $D_J$ separately. Thus, we 
select subsets for prediction that are common to all regressions having the same sample size, but we allow different parameters for the regressions.

Given a distribution $G_j \in \cal G$, 
let $m_j({\bf X}_{ij}):= E_{G_j}(Y_{ij}\mid {\bf X}_{ij})$ be the conditional expectation under $G_j$. We do not assume a linear model or any particular model for $m_j$ when we analyze our procedures, but for the sake of prediction we shall approximate $m_j({\bf X}_{ij})$ by a linear function ${\bf X}'_{ij}{\bs \beta}_j$, where ${\bs \beta}_j$ is the vector of projection coefficients under $G_j$.   We shall require minimal assumptions on $G_j$ such as moment conditions, to be specified later.

When $m_j({\bf X}_{ij})$ is not linear then ${\bf X}_{ij}$ is not ancillary, and its marginal distribution matters; see, e.g, \citet{Buja}. In this case, conditioning on ${\bf X}$ or considering it as nonrandom leads to loss of information. For a recent discussion on fixed versus random ${\bf X}$ in the context of model selection see \citet{Rosset}. 
When $m_j({\bf X}_{ij})={\bf X}'_{ij}{\bs \beta}_j$, allowing linear models with different coefficients in different regressions is called \emph{heterogeneous regression ANCOVA}; see, e.g., \citet{Rutherford}, Chapter 8, and the references therein.  Related models appear under titles such as repeated measure regression (see, e.g., \citet{VC}), often with mixed effects.

Given datasets $\{{\bf X}_{ij}, Y_{ij}\}$ from $G_j$, consider the subset of covariates  $\mathpzc P$ of size  $p \le d$. We may sometimes refer to $\mathpzc P$ as a \textit{model}.   Let ${\bf X}_{ij}^{(\mathpzc P)}$ denote the subvector of ${\bf X}_{ij}$ consisting of the covariates in $\mathpzc P$. Let ${\bs \beta}^{(\mathpzc P)}_j$ denote the  linear projection coefficient vector  and let  $\widehat{{\bs \beta}}^{(\mathpzc P)}_{j,n}$ be its least squares  estimator based on $n$ observations, where we assume that $n>p$. 
In Section \ref{sec:BG} we discuss the case of discrete covariates in which exact (or perfect) multicollinearity may occur with a positive probability, and the least squares  estimators are not unique.

For now we focus on the case that $\cal J=K$, that is we observe all datasets in $\cal G$. (In Section \ref{sec:supersuper} we  consider prediction of an out of (the training) sample dataset, in the spirit of transfer learning.) Consider prediction for a regression dataset of some size $n$, often referred to as the \textit{task size},  which will be taken to equal $N_j$ for the task of predicting for the dataset $D_j$ in the multi-task problem of prediction for datasets in the training sample.   The linear prediction of a response $Y$, based on $n$ observations from a random $G_J \in \cal G$, and when the subset $\mathpzc P$ is used 
is given by $({\bf X}^{(\mathpzc P)})'\, \widehat{\bs \beta}^{(\mathpzc P)}_{J,n}$. In order to select a subset for regression tasks of size $n$ we make the counterfactual assumption that all datasets in the training sample are of size $n$. Then 
the corresponding expected prediction error or risk is given by 
\begin{equation}\label{eq:rtild}
{\bf R}(n,\mathpzc P):=\frac{1}{{\cal J}} \sum_{j=1}^{\cal J} R_j(n,\mathpzc P):=\frac{1}{{\cal J}} \sum_{j=1}^{\cal J} E_{G_{j}}\big(Y-({\bf X}^{(\mathpzc P)})'\, \widehat{\bs \beta}^{(\mathpzc P)}_{j,n}\big)^2,
\end{equation}
where $({\bf X},Y)\sim G_j$ independently of $\widehat{\bs \beta}^{(\mathpzc P)}_{j,n}$, and the expectation on the right-hand side of \eqref{eq:rtild}  applies to both $({\bf X},Y)$ and  $\widehat{\bs \beta}^{(\mathpzc P)}_{j,n}$.
This expression has the alternative interpretation where instead of predicting for a random $D_J$ we predict for all $D_j$, $j=1,\ldots,\cal J$, assuming that there is a common sample size $n$, and now  ${\bf R}(n,\mathpzc P)$ represents the average prediction error. With either interpretation,
our goal is to estimate ${\bf R}(n,{\mathpzc P})$ and related quantities, in order to select (with high probability) an optimal common subsets ${\mathpzc P}^*(n)$  for prediction for any dataset in the training set (taking $n=N_j$) and also for out of sample datasets when we later consider the case that $\cal J < K$, on the basis of $n$ observations, where the subset selection is based on the pooled training set $\cal T$.

It is natural to choose common subsets for prediction if the different regression datasets arise from a common source; besides efficiency in  subset selection due to pooling, common subsets lead to  computational efficiency. However, we assert that in 
a variety of situations (but obviously not always) it is advantageous   to choose  standard common sets of  covariates to be used for prediction even if the regression datasets do not arise from a homogeneous source. In this case we are trying to select   compromise subsets that can be used for the different regressions (and may not be optimal for some or  any of them). For example, a large health organization with $\cal K$ clinics  often recommends a common standard set of tests for the purpose of certain  diagnoses, thereby simplifying the instructions to participating clinics and doctors. In our notation, the set of tests is based on a sample of size $\cal J$, which is in general 
$\le \cal K$. The regression coefficients used for prediction based on this common set of tests may differ between communities or doctors, who may attach different weights to different tests. Concerning economics models, consider the OECD, where ${\cal J=\cal K} =37$ (as of 2021) since all countries are sampled and economic prediction are made in all of them.  The OECD attempts to standardize sets of common economic indicators to be used for economic predictions (e.g., forecast   of GDP growth) for its member countries, which are to be estimated by their bureaus of statistics by the same methodology. In general, it makes sense to assume that in different countries, economic variables may have different weights in economic predictions. For example, oil prices must weigh differently for economic predictions between oil importing and exporting countries.

\subsection{GENO, a measure of usefulness}\label{sec:GENOInt}
In order to compare the  quality of different models, we introduce a new measure, GENO, which is inspired by the measure ENO (equivalent number of observations) of \citet{Erev}. 
To describe ENO in the context of experimental economics, consider an experiment where a game is  played by a sample of subjects in order to study the average behavior of players, and predict future play. ENO is based on a comparison between the empirical statistics of past actions of the players, and a given model for predicting players' actions. 
The more subjects who have already played the game, the better the estimate that past play will give of the mean behavior of the subject population on this game. ENO  measures the usefulness of the prediction of a particular model by asking how many prior observations of subjects playing the game, say $m$, would be needed to make the empirical statistics as  accurate as the prediction by the model. ENO of a model is this number $m$.

While ENO compares a given model to the relevant empirical model, GENO generalizes ENO to comparing any two data-based models. Thus, let now ${\bf R}(n,\mathpzc P)$ denote the prediction error of some model $\mathpzc P$ in a very general setup. For our present purposes, one can have regression models in mind, with ${\bf R}(n,\mathpzc P)$ defined above; however, the definition of GENO below is more general.
Given two models   $\mathpzc P$  and $\mathpzc Q$, define  GENO$(n;\mathpzc P,\mathpzc Q)$ to be the value of $m$ satisfying ${\bf R}(m,\mathpzc Q)={\bf R}(n,\mathpzc P)$. In words, GENO$(n;\mathpzc P,\mathpzc Q)$ is the number of observations  required in order for a model based on the covariates in $\mathpzc Q$ to predict equally well as a model based on the covariates in $\mathpzc P$, when the parameters of the latter model are estimated on the basis of $n$ observations. In Section \ref{sec:GENO} we shall use an approximation to ${\bf R}(n,\mathpzc P)$ to formally define and  estimate GENO.  Such a measure allows us to decide between a set of covariates that may be good for prediction but costly to obtain, and another set of more accessible covariates that we may consider using, even if their predictive value is lower and therefore may require more observations. See the recent paper - Andrade et al.  \cite{Andrade} and the references therein for a formal Bayesian approach to minimizing cost of classification in the presence of costly covariates. A comparison in terms of the sample size required by one model (for prediction,  testing  or estimation) to be as good as another with a given sample size is closely akin to the notion of Pitman efficiency; see, e.g., \citet{Zacks}.
Our approach to quantifying the value of a model is close in spirit, but not in detail, to the work of \citet{Lindsay} who define a ``model credibility index'' as the sample size $N^*$, where data from the model and from the true generating process are indistinguishable in the sense that for a given goodness of fit test of the model with $N^*$ observations, the probability of rejection under the model  is, say, 50\%. 

A different approach to measuring the usefulness of a model is by Akaike weights, which are defined by the likelihood function of each model evaluated at the MLE, standardized by their sum; see \citet{Anderson} (Page 75), where these weights are referred to informally as ``the weight of evidence in favor of model." The AIC weights are sometimes \citep[e.g.,][]{Wagen} interpreted as probabilities of a model to be the best in terms of the AIC criterion. With a uniform prior on the set of models this interpretation could be meaningful if we believe that one of the models is true. Otherwise, the weights are still informative, but their interpretation is less clear. GENO, on the other hand, is measured in units of number of observations, which are easy to grasp. Another advantage of GENO is that it accounts for the number of observations in the data to which the model is applied. This makes sense as the usefulness of a model for a given dataset is also a function of the size of the data.

In Section \ref{sec:predd} we restate the problem and provide some basic results and notation for a single regression dataset, as a preliminary to the main part, Section \ref{sec:Jind}, where we consider the multi-task problem of model selection for several regression datasets. 
In Section \ref{sec:GENO} we discuss the  GENO measure of the relative quality of models. In Section \ref{sec:sims}
we demonstrate the  results by simulations, and in Section \ref{sec:med} we discuss an application to a medical management problem of predicting service times, that is, visit durations of patients in hospital. Section \ref{sec:AppA} is an appendix  containing the proofs. Appendix B summarizes the notation used in the paper.

\section{Prediction error with random covariates: A  single dataset} \label{sec:predd}
We start with   $|{\cal G}|={\cal J}=1$, that is, with selection of a model for prediction given a training set consisting of a single regression dataset. This  case is treated in the standard model selection literature. Although our real interest is in   results for large $\cal J$,  we consider ${\cal J}=1$ as a starting point which  simplifies the notation while allowing us to present some of the ideas used in the general case.  
For now  our training set $\cal T$ consists of a single dataset $D_1=\{({\bf X}_i, Y_i)\}_1^N$  of  $N:=N_1$ iid pairs from some distribution $G:=G_1$.
The distinction between $N$  and $n$ may seem artificial in this case, but we shall make it and consider prediction based on any sample size $n$ for later purposes.  We use $D_1$  for selecting a subset of covariates for linear prediction of a future $Y$ from ${\bf X}$ distributed by $G$, with  parameters that will be estimated using a dataset $D=\{({\bf X}_i, Y_i)\}_1^n$  of $n$  observations from $G$.

We derive some results that will be needed for the general case ${\cal J}>1$, to be discussed in Section \ref{sec:Jind}.
Subsets of the covariates are denoted by letters like $\mathpzc P$, $\mathpzc Q$, etc., and their sizes by $p$ and $q$, etc.   We refer to the associated linear model as model $\mathpzc P$. For now
we fix $\mathpzc P$ and  suppress it in most of our notation and instead of ${\bf X}^{(\mathpzc P)}$ we write ${\bf X}$ and assume it is in $\mathbb{R}^p$. The same holds for other vectors and matrices. Later we shall assume that ${\bf X} \in \mathbb{R}^d$, and consider different subsets of covariates.

\subsection{Preliminaries} \label{subsec:single}
Consider a dataset $D=\{({\bf X}_i, Y_i)\}_1^n$ of iid pairs from some distribution $G$, where ${\bf X}_i$ is a column vector in ${\mathbb{R}}^p$, $i=1, \ldots,n$. Let $({\bf X},Y)$ without indexes  denote one such ``generic" observation, distributed independently of the dataset $D$ as any $({\bf X}_i, Y_i)$ according to $G$.   The first entry of each ${\bf X}_i$ may be 1, so that the models may include an intercept term.

Set $\mathbb{Q}:=E({\bf X} {\bf X}')$ and
let ${\bf{Y}}_n \in {\mathbb{R}}^n$ denote the $n$-column vector of the $Y_i$'s, and set $m({\bf X}):=E(Y|{\bf X})$ for some function $m$. Assuming that both ${\bf X}$ and $Y$ have finite second moments and that $\mathbb{Q}$ is invertible, the best linear approximation of $m(\bf{X})$ is $\bf{X}'{\bs \beta}$, where
\begin{equation}\label{eq:beta} 
{\bs \beta}:=\arg\min_{{\bf b} \in \mathbb{R}^p} E(m({\bf X})-{\bf X}'{\bf b})^2 =\mathbb{Q}^{-1}E({\bf X}Y).
\end{equation}
The same projection coefficient vector ${\bs \beta}$ also satisfies 
${\bs \beta}=\arg\min_{\bf b}E(Y-{\bf X}'{\bf b})^2$; hence $\bf{X}'{\bs \beta}$ is the best linear predictor of $Y$. Our assumptions imply that the minimizer ${\bs \beta}$ is unique.
Set $e_i:=Y_i-{\bf X}'_i{\bs \beta}$, with ${\bs \beta}$  defined in \eqref{eq:beta}. By (2.25) in \citet{Hansen}, where most of our notation and the standard results we use can be found, we have $E({\bf X}e)=\bf 0$, where again ${\bf X}$ and $e$ are ``generic"
${\bf X}_i$ and $e_i$.

Define  $\mathbb{X}_n$
to be the $n \times p$ matrix whose $n$ rows are  the row vectors ${\bf X}'_i$.  In this common notation  the standard linear model will be written as $\mathbb{X}_n{\bs \beta}$, whereas each  of its rows  as ${\bf X}'_i{\bs \beta}$, and  $\mathbb{X}_n'\mathbb{X}_n=\sum_{i=1}^n {\bf X}_i{\bf X}_i'$. 
Under standard assumptions, the least squares estimator is 
\begin{equation}\label{eq:linbeta}
\widehat{\bs \beta}_n:=\arg\min_{{\bf b} \in \mathbb{R}^p}||{\bf Y}_n-\mathbb{X}_n{\bf b}||^2=(\mathbb{X}_n'\mathbb{X}_n)^{-1}\mathbb{X}_n'{\bf{Y}}_n.
\end{equation}
The assumption that $(\mathbb{X}_n'\mathbb{X}_n)^{-1}$ exists (with probability 1) holds if we assume that ${\bf X}$ has a continuous distribution. For the existence of certain moments required later we shall assume that the distribution of ${\bf X}$ is a mixture of normals. See \citet{Hansen}, pp. 102--3, for a discussion of the existence of $(\mathbb{X}_n'\mathbb{X}_n)^{-1}$ and its moments. Without assuming continuity, the assumption that $\mathbb{Q}$ is invertible implies that $(\mathbb{X}_n' \mathbb{X}_n)^{-1}$ exists with probability converging to 1 as $n\to \infty$; however, for discrete distributions this probability is smaller than one, and thus $\widehat{\bs \beta}_n$ may not exist, and has no finite moments, a ``conundrum" in the words of Hansen.
In Section \ref{sec:BG} we extend our discussion to discrete covariates by conditioning  on the existence of a bounded inverse, and showing that under simple conditions this amounts to neglecting a set having an exponentially small probability, thus providing some solution to the  above conundrum. 

We now assume that $(\mathbb{X}_n'\mathbb{X}_n)^{-1}$ exists and has sufficiently many moments so that expressions like \eqref{eq:kprederr} below are finite. If ${\bf X}$ and $Y$ have finite fourth moments,  then by Theorem 7.3 in \citet{Hansen}
\begin{equation}\label{eq:CLT}
\sqrt{n}(\widehat{\bs \beta}_n-{\bs \beta}) \tendd N({\bf 0}, \mathbb{Q}^{-1}\mathbb{W}\mathbb{Q}^{-1}), 
\end{equation}
where  $\mathbb{W}:=E({\bf X}{\bf X}'e^2)$, 
a $p \times p$ matrix assumed to be positive definite.
For a single distribution $G$ and a dataset $D$ as above, the prediction error  incurred by a model ${\mathpzc P}$ based on all
$p$ covariates with linear regression coefficients computed from  a sample of size $n$ is 
\begin{equation}
\label{eq:kprederr}
R(n,{\mathpzc P})=E_G\big( Y-{\bf X}'\, \widehat{{\bs \beta}}_n\big)^2.
\end{equation}
Later we assume that ${\bf X} \in \mathbb{R}^d$ and set 
${\bf X}^{(\mathpzc P)} \in \mathbb{R}^p$ 
to be the vector consisting of the covariates of ${\bf X}$ in the subset of covariates $\mathpzc P$ of size $p$.  
When we consider several models,  we set, for example, 
$\mathbb{X}^{(\mathpzc P)}_n$
to be the $n \times p$ matrix whose $n$ rows are  the row vectors ${\bf X}^{'(\mathpzc P)}_i$, 
$\widehat{\bs \beta}^{(\mathpzc P)}_n:=(\mathbb{X}^{'(\mathpzc P)}_n\mathbb{X}^{(\mathpzc P)}_n)^{-1}\mathbb{X}^{'(\mathpzc P)}_n{\bf{Y}}_n$, 
$\mathbb{W}^{(\mathpzc P)}:=E({\bf X}^{(\mathpzc P)}{\bf X}^{(\mathpzc P)'}e^2)$, and likewise for $\mathbb{Q}$, etc. We then have 
\begin{equation}
\label{eq:kpr2}
R(n,\mathpzc P)=E_G\big( Y-{\bf X}^{'(\mathpzc P)}\, \widehat{{\bs \beta}}^{(\mathpzc P)}_n\big)^2.
\end{equation}

\subsection{Equally good sequences of models} \label{sec:Neighborhoods} 
When selecting the best model for a given $n$, that is, the  subset of covariates that minimizes $R(n,\mathpzc P)$, we should take into account that different samples yield different estimators $\widehat{{\bs \beta}}_n$, leading to different prediction errors; thus, there is no gain in optimizing more precisely than the difference between such errors.  
Consider  the prediction error        conditioned on the estimated regression coefficients
\[
R\big(n,{\mathpzc P};\widehat{{\bs \beta}}_n\big):=E\Big[\big( Y-{\bf X}'\, \widehat{{\bs \beta}}_n\big)^2\Big|\widehat{{\bs \beta}}_n \Big].
\]
Note that $R(n,{\mathpzc P})=E\{R(n,{\mathpzc P};\widehat{{\bs \beta}}_n)\}.$
By using the relation  $\big( Y-{\bf X}'\,\widehat{{\bs \beta}}_n\big)^2=\big( Y-{\bf X}'\,{{\bs \beta}}+{\bf X}'\,{{\bs \beta}} - {\bf X}'\,\widehat{{\bs \beta}}_n\big)^2$, expanding the latter term, 
and taking conditional expectation noting that  $E[(Y-{\bf X}^{'}{{\bs \beta}}){\bf X}^{'} | \widehat{{\bs \beta}}_n]=E[e {\bf X}']=0$, we obtain
\begin{equation}\label{eq:RnkBeta}
R\big(n,{\mathpzc P};\widehat{{\bs \beta}}_n\big)
=E\left( Y-{\bf X}'\,{{\bs \beta}}\right)^2+E\Big\{\Big[{\bf X}'\, (\widehat{{\bs \beta}}_n-{{\bs \beta}})\Big]^2 \Big| \widehat{{\bs \beta}}_n \Big\}.
\end{equation}
The first term in \eqref{eq:RnkBeta} is a constant and the second equals $ (\widehat{{\bs \beta}}_n-{{\bs \beta}})'\mathbb{Q} (\widehat{{\bs \beta}}_n-{{\bs \beta}})$, which is of order $O_p(1/n)$ since $\sqrt{n}(\widehat{{\bs \beta}}_n-{{\bs \beta}})=O_p(1)$; see \eqref{eq:CLT}. This means that $R\big(n,{\mathpzc P};\widehat{{\bs \beta}}_n\big)$ varies between different 
$\widehat{\bs \beta}_n$ by a quantity of order $O_p(1/n)$. Hence, if two sequences of models ${\mathpzc P}(n)$ and ${\mathpzc Q}(n)$ satisfy $$|R(n,{\mathpzc P}(n))-R(n,{\mathpzc Q}(n))| =o(1/n), \text{ i.e,  }
\lim_{n \to \infty}n|R(n,{\mathpzc P}(n))-R(n,{\mathpzc Q}(n))|=0,$$ we consider them to be  \textit{equally good}. If ${\mathpzc P}(n)$ is best in the sense of minimizing $R(n,{\mathpzc P}(n))$ and ${\mathpzc Q}(n)$ is equally good, we say that ${\mathpzc Q}(n)$ is \textit{adequate}, and rather than choose ``best models" we settle for adequate models. See, e.g., \citet{Nevo} for a related approach.

\subsection{Versions of Mallows $C_p$ for random covariates}\label{subsec:mallows}

Given a dataset $D_1=\{({\bf X}_i, Y_i)\}$ of size $N$ (which constitutes the training set when ${\cal J}=1$), we first estimate the prediction error \eqref{eq:kprederr} incurred if prediction is to be based on $n$ observations.  We shall consider two types of asymptotics: one when $n$ is considered to be large, and the other when
$n$ is fixed, and $N$ is large. For now ${\cal J}=1$; asymptotics in ${\cal J}$ will be considered later.

We use the following notation: set
$\widehat{\mathbb{Q}}_N:=\frac{1}{N}\mathbb{X}_N'\mathbb{X}_N$, and let
${\bf Y}_N$ denote the $N$-vector of the $Y_i$'s. Recalling the notation  ${e}_i=Y_i-{\bf X}_i'{{\bs \beta}}$, let ${\bf e}_N$ denote the $N$-vector having components $e_i$. Set  $\widehat{\mathbb{W}}_N :=\frac{1}{N}\sum_{i=1}^N{\bf X}_i{\bf X}_i'\widehat{e}_i^{\,2}$  with\, $\widehat{e}_i=Y_i-{\bf X}_i'\widehat{{\bs \beta}}_{N}$, where $\widehat{\bs \beta}_{N}$ is given by \eqref{eq:linbeta} upon replacing $n$ by $N$. Thus in \eqref{eq:cpp} below, $\frac{1}{N}||{\bf Y}_N-\mathbb{X}_N\widehat{\bs \beta}_N||^2=\frac{1}{N}\sum_{i=1}^N \widehat{e}_i^2$.
%and let $\widehat{e}=Y-{\bf X}'\widehat{{\bs \beta}}_{N}$
%without subindexes  denote one such ``generic" observation for some $i$. 
Let ${\mathbb{V}}:={\mathbb{W}}{\mathbb{Q}}^{-1}$, and  ${\widehat{\mathbb{V}}_N}:=\widehat{\mathbb{W}}_N\widehat{\mathbb{Q}}_N^{-1}$.  
In addition we define ${\bf U}_N := \frac{1}{\sqrt{N}}\mathbb{X}_N' {\bf e}_N=\frac{1}{\sqrt{N}}\sum_{i=1}^N {\bf X}_i e_i$. Note that ${\bf U}_N$ is not a statistic and that $E({\bf U}_N {\bf U}_N')= {\mathbb{W}}$ since E$({\bf X}e)=\bf 0$ implies that the expectations of mixed terms vanish.
In all the  vectors and matrices above and below the index ${\mathpzc P}$  was suppressed unless otherwise indicated.

Akin to \eqref{eq:kprederr}, we define the approximate prediction error to be
\begin{equation}\label{eq:approx_predara} AR(n,{\mathpzc P}):=E\big(Y-{\bf X}^{'}{{\bs \beta}}\big)^2+\frac{1}{n}tr({{\mathbb{V}}}),
\end{equation}
where $tr$ denotes trace, and equation \eqref{eq:approx_pred3} of Theorem \ref{thm:o1n} below shows that it is an approximation to the quantity $R(n,{\mathpzc P})$ of \eqref{eq:kprederr}. 
Clearly $E\big( Y-{\bf X}'\, \widehat{{\bs \beta}}_n\big)^2  \ge E\big(Y-{\bf X}^{'}{{\bs \beta}}\big)^2$ and the trace is an approximation of the difference with precision of order $O(1/n^{3/2})$; see \eqref{eq:approx_pred3}. %as follows from \eqref{eq:R_AR} and the ensuing calculations. 
Next we define the statistic $C^{({\mathpzc P})}(n,N)$ as an estimator of $AR(n,{\mathpzc P})$ by
\begin{equation}
\label{eq:cpp}
C^{({\mathpzc P})}(n,N):=\frac{1}{N}||{\bf Y}_N-\mathbb{X}_N\widehat{\bs \beta}_N||^2
+ {tr}({\widehat{\mathbb{V}}_N})\left( \frac{1}{n} +\frac{1}{N}\right).
\end{equation}
The new term $\frac{1}{N}{tr}({\widehat{\mathbb{V}}_N})$ is an approximately (up to $o_p(1/N)$) unbiased estimator of $\frac{1}{N}||{\bf Y}_N-\mathbb{X}_N\widehat{\bs \beta}_N||^2-E\big(Y-{\bf X}^{'}{{\bs \beta}}\big)^2$, as  shown in \eqref{eq:stam1} and \eqref{eq:tracesg}.
The fact that ${tr}({\widehat{\mathbb{V}}_N})$ is a biased estimator of ${tr}({{\mathbb{V}}})$ entails a bias of order $1/n$ for  the estimator $C^{({\mathpzc P})}(n,N)$ as an estimator of $AR(n,{\mathpzc P})$. We shall study the latter estimator, and when we use it, we shall apply a standard jackknife correction for its bias; see  \citet{Efron}, Equation (2.8). We denote the bias-corrected $C^{({\mathpzc P})}(n,N)$ by $\mathbb{C}^{({\mathpzc P})}(n,N)$. It suffices to bias-correct only ${tr}({\widehat{\mathbb{V}}_N})$ in \eqref{eq:cpp} as explained in the fourth paragraph after Theorem \ref{thm:o1n}.

The superscript ${\mathpzc P}$ in the statistic $C^{({\mathpzc P})}$ refers to the set of covariates in ${\mathpzc P}$ and for now we have
${\bf X}_i^{({\mathpzc P})}={\bf X}_i \in \mathbb{R}^p$ and the subset ${\mathpzc P}$ is fixed and suppressed. The statistic $C^{(\mathpzc P)}(n,N)$ is a counterpart of Mallows $C_p$, but here we consider random covariates. Furthermore, we distinguish between the number  $N$ of observations used for the choice of  the model and the sample size  $n$ of observations used for estimating the model's parameters. The classic Mallows $C_p$ concerns  nonrandom covariates, where $n=N$, and the true model is assumed to be linear.
To see the relation to Mallows $C_p$, assuming a homoskedastic linear model, we have that 
$e_i=Y_i-{\bf X}'_i{\bs \beta}$ is uncorrelated with the covariates, with variance $\sigma^2$, and $\mathbb{W}=\sigma^2 \mathbb{Q}$, and therefore $\widehat{\mathbb{V}}_N=\widehat{\mathbb{W}}_N
\big(\widehat{\mathbb{Q}}_N\big)^{-1}$
will converge to $\sigma^2 I_p$ and ${tr}({\widehat{\mathbb{V}}_N})$ to $\sigma^2 p$. If we use $\sigma^2 p$  as an approximation of ${tr}({\widehat{\mathbb{V}}_N})$ (and therefore we only have to estimate $\sigma^2$ rather than a trace), then $C^{({\mathpzc P})}$ in the case $N=n$ coincides with Mallows $C_p$.

The following theorem provides the rate of approximation of $AR(n,{\mathpzc P})$ to $R(n,{\mathpzc P})$, and then analyzes $C^{({\mathpzc P})}$ as an estimator of $AR(n,{\mathpzc P})$;  
some of its conditions and implications are discussed below. All proofs are in the Appendix.
Our proof shows that  Assumption (i) below can be replaced by the assumption that ${\bf X}$ and $Y$ have 24 finite  moments, and a careful inspection of the proof shows that this number can be somewhat reduced.  
\begin{theorem}\label{thm:o1n}
	Assume that\\ {\rm (i)} The coordinates of ${\bf X}$ and $Y$ have finite  moments of all orders. \\
	{\rm (ii)} The entries of $(\mathbb{X}_n'\mathbb{X}_n/n)^{-1}$ have third moments that are bounded uniformly in $n$.
	Then
	\begin{equation}\label{eq:approx_pred3}
	|R(n,{\mathpzc P}) - AR(n,{\mathpzc P})| = O(1/{n^{3/2}}),
	\end{equation}
	and
	\begin{equation}\label{eq:final}
	AR(n,{\mathpzc P})-C^{({\mathpzc P})}(n,N) = {\cal E}_N+\frac{1}{n}\left\{{tr}({\mathbb{V}})-{tr}({\widehat{\mathbb{V}}_N})\right\} +o_p(1/{N}),
	\end{equation} 
	where %\ycomm{changed order on RHS of (11) to make product of matrices it look like $\mathbb{V}$ ALSO in proof?}
	\begin{equation}\label{eq:CAPEPS}
	{\cal E}_N =E(Y-{\bf X}'{{\bs \beta}})^2-\frac{1}{N}||{\bf Y}_N-\mathbb{X}_N{\bs \beta}||^2 + \frac{1}{N}\{tr({\bf U}_N{\bf U}_N'\mathbb{Q}^{-1})-tr({\mathbb{V}})\}.
	\end{equation} 
	Furthermore, 
	\begin{equation}\label{eq:further}
	(a) \quad	{\cal E}_N=O_p(1/\sqrt{N}),\qquad (b) \quad  {tr}({\mathbb{V}})-{tr}({\widehat{\mathbb{V}}_N})=O_p(1/\sqrt{N}),
	\end{equation} 
	and 
	\begin{equation}\label{eq:further2}	
	\sqrt{N}\big(C^{({\mathpzc P})}(n,N)- AR({\mathpzc P},p)\big ) \tendd N(0,\uptau^2) 
	\end{equation}  
	for some asymptotic variance $\uptau^2$ as $N \rightarrow \infty$, and $n$ is fixed. 
\end{theorem} 
Since there is only a finite number of models, the above terms  $O$, $O_p$, and $o_p$ do not depend on the subset of covariates ${\mathpzc P}$. For example, we could replace \eqref{eq:approx_pred3} by $|R(n,{\mathpzc P}) - AR(n,{\mathpzc P})| \le B/{n^{3/2}}$ for all $n$ and ${\mathpzc P}$, where $B$ is a constant. Moreover, the term $o_p(1/{N})$ in \eqref{eq:final} does not depend on $n$.

Condition (i) of Theorem of \ref{thm:o1n} is standard, and  Lemma \ref{lem:normal} below shows that Condition (ii) is satisfied if ${\bf X}$ is distributed as a mixture of normals; see \citet{Sampson}. Such mixtures form a dense family of distributions with respect to weak convergence in the space of distribution on  $\mathbb{R}^p$. As the distribution of ${\bf X}$ is never known exactly, it makes sense to assume, as an approximation, that the data satisfy such a condition. 
The case where $\bf X$ has discrete  components is discussed in Section \ref{sec:BG}.

We shall later compare models consisting of different subsets of covariates. Equation \eqref{eq:approx_pred3} suggests that choosing a model by minimizing a  good estimate of $AR(n,\mathpzc P)$ with respect to $\mathpzc P$ can lead to a model for which $R(n,\mathpzc P)$ is within $o({1}/{n})$ of the best model, and thus $\mathpzc P$ is an adequate model in the sense of Section \ref{sec:Neighborhoods}.  This is stated formally in Proposition \ref{prop:single}.

In view of \eqref{eq:final} we use $C^{({\mathpzc P})}(n,N)$ of \eqref{eq:cpp}  as an estimator of the approximate prediction error $AR(n,{\mathpzc P})$ and hence of the prediction error $R(n,{\mathpzc P})$. 
This is formalized in Propositions \ref{prop:cons} and \ref{prop:nfixed} below. 
We now briefly discuss Equations \eqref{eq:final} and \eqref{eq:CAPEPS}. First  consider the bias of $C^{({\mathpzc P})}(n,N)$ as an estimator of $AR(n,{\mathpzc P})$. It is easy to see that  $E{\cal E}_N=0$. By \eqref{eq:further} (b), ${tr}({\mathbb{V}})-{tr}({\widehat{\mathbb{V}}_N}) =O_p(1/\sqrt{N})$, and after dividing the latter term by $n$ as in  \eqref{eq:final}, it  is of a smaller order than the term  ${tr}({\widehat{\mathbb{V}}_N})\left( \frac{1}{n} +\frac{1}{N}\right)$ appearing in $C^{({\mathpzc P})}(n,N)$.
This shows that the latter term   contributes to  reducing the bias of $C^{({\mathpzc P})}(n,N)$ as an estimator of $AR(n,{\mathpzc P})$. 

Our main interest is in the case of ${\cal J}>1$ regressions, and in choosing a model that minimizes an average of $\cal J$ values of $AR$. Averaging (nearly) unbiased estimates can result in consistency in $\cal J$, which explains why we care about correcting the bias of $C^{({\mathpzc P})}(n,N)$. In this case, a further bias correction using the jackknife is useful (see Section \ref{sec:simJ}).  The above discussion implies that it suffices to bias-correct the estimator ${tr}({\widehat{\mathbb{V}}_N})$, which is what we do when using the jackknife. 

Choosing a good model can be reduced to choosing between two models, say, $\mathpzc P$ and $\mathpzc Q$ at a time, by approximating the difference $AR(\mathpzc P)-AR(\mathpzc Q)$ using $C^{(\mathpzc P)}(n,N)-C^{(\mathpzc Q)}(n,N)$. The leading terms in the latter expression will be the difference between the relevant values of ${\cal E}_N$ for the two models, and it is easy to see that the leading term of this difference is the difference between the values of $\frac{1}{N}||{\bf Y}_N-\mathbb{X}_N{\bs \beta}||^2$ for the corresponding models, which is of order $O_p(1/\sqrt{N})$ by the central limit theorem.
However, when two models having  very similar prediction values are compared 
by differencing their corresponding values of $C^{(\mathpzc P)}(n,N)$,  their leading terms will approximately cancel, and in this case the second term on the right-hand side of \eqref{eq:cpp} plays a role. This holds also for Mallows $C_p$ and the AIC, \citep{Akaike}, and will be exploited formally in the Propositions \ref{prop:cons} and \ref{prop:nfixed} below.  

The following lemma shows that Condition (ii) of Theorem \ref{thm:o1n} holds when ${\bf X}$ is distributed as a mixture of normals.

\begin{lemma}\label{lem:normal} Let the distribution of the  covariate vectors (excluding the first coordinate in the case that it is a constant 1) be normal, or a finite mixture of normals, or an  infinite mixture of normals with covariance matrices in a set $\varXi$, and $\inf_{\Sigma \in \varXi} \lambda_{min}(\Sigma)>0$, where $\lambda_{min}$ denotes the smallest eigenvalue. Then, for $n>p+5$, Condition {\rm (ii)} of Theorem \ref{thm:o1n}
	is satisfied.
\end{lemma}
More generally, the $r$th moments of the entries of $(\mathbb{X}_n'\mathbb{X}_n/n)^{-1}$ are bounded under the conditions of Lemma \ref{lem:normal} provided that $n>p+2r-1$ (see \citet{von_Rosen}, Theorem 4.1). Note that the condition on  $\lambda_{min}$ guarantees that $\bf X$ is bounded away from exact multicollinearity.

\subsection{Discrete covariates}\label{sec:BG}

When ${\bf X}$ contains discrete covariates, the probability that the matrix $(\mathbb{X}_n'\mathbb{X}_n/n)^{-1}$ does not exist is  positive, and expressions like 
$\widehat{\bs \beta}_n$ of \eqref{eq:linbeta} and hence
$R(n,p)$ of  \eqref{eq:kprederr} may not exist. When the components of $\bf X$ are bounded, we provide the following limiting approach.
Set 
\begin{equation}\label{eq:setGn}
H_n:=\big\{\mathbb{X}_n\,:\,\lambda_{min}(\mathbb{X}_n'\mathbb{X}_n/n)
\ge\lambda_{min}(\mathbb{Q})/2\big\},
\end{equation} 
where $\lambda_{min}$ is the smallest eigenvalue, and  $\widetilde R(n,{\mathpzc P}):=E\big[\big( Y-{\bf X}'\, \widehat{{\bs \beta}}_n\big)^2 \mid H_n\big]$.
%\begin{equation*}
%\label{eq:kprederrcond}
%\end{equation*}
We have 
\begin{theorem}
	\label{thm:BeGood}
	Suppose that $Y$ has all moments,   the components of $\bf X$ are bounded, and $\mathbb Q$ is invertible; then for some $a \in (0, 1)$,
	\begin{equation*}
	\big|\widetilde R(n,{\mathpzc P}) - AR(n,{\mathpzc P})\big|=O(1/n^{3/2}) \text{    and  } P(H_n) > 1-a^{n\lambda_{min}(\mathbb{Q})}.
	\end{equation*}
	Moreover, all quantities appearing in Theorem \ref{thm:o1n} are well defined on $H_N$, and can be defined in an arbitrary way outside of $H_N$, and the results \eqref{eq:final}--\eqref{eq:further2} hold.
\end{theorem}
Thus, apart from the complement $H^c_n$, which has exponentially small probability, the approximation rate of $AR(n,{\mathpzc P})$ to the prediction error is  the same as in 
\eqref{eq:approx_pred3} and the rest of Theorem \ref{thm:o1n} still holds. The result follows from Theorem \ref{thm:o1n} and Lemma 2.3  given in the Appendix.

\subsection{Approximations and consistency}\label{sec:approx}
The focus of this section is on choosing a subset of covariates  for prediction of future responses on the basis of a single dataset of size $N$.  The linear model parameters  are estimated from a sample of size $n$, with the understanding that different $n$'s may (and should) lead to different choices of subsets; more specifically,  a larger $n$ naturally gives rise to a larger set of covariates.   Asymptotic results in $n$ are not of major interest in this context; however, they may contribute some understanding when $n$ is not small. Such results are discussed in this section.

In Proposition \ref{prop:cons} we show  that under the conditions of Theorem \ref{thm:o1n}, choosing a subset of covariates in the set $\arg \min_{\mathpzc P} C^{(\mathpzc P)}(n,N)$ guarantees that for increasing  $n$ and $N$ we choose the best linear model with probability converging to 1, that is, the model minimizing $R(n,\mathpzc P)=E\big(Y-{\bf X}^{(\mathpzc P)'}\widehat{{\bs \beta}}_n^{(\mathpzc P)}\big)^2$, with notation  defined after \eqref{eq:kprederr}.   In Proposition \ref{prop:nfixed} we show that for fixed $n$, using $C^{(\mathpzc P)}(n,N)$, we choose 
an adequate model  in the sense defined in  Section \ref{sec:Neighborhoods},
with probability converging to 1 as $N \to \infty$. 

Below $\arg\min_{\mathpzc P}$ is taken over all subsets of covariates.
For a given $n$,  define the following sets: 
\begin{align*}
&{\mathpzc P}^*(n) :=\arg \min_\mathpzc P R(n,\mathpzc P)=\arg \min_{\mathpzc P} E\big(Y-{\bf X}^{(\mathpzc P)'}\widehat{{\bs \beta}}_n^{(\mathpzc P)}\big)^2,\\
&\pi^*(n)  := \arg \min_{\mathpzc P} AR(n,{\mathpzc P})=\arg \min_{\mathpzc P}\Big\{E\big(Y-{\bf X}^{({\mathpzc P})'}{{\bs \beta}^{({\mathpzc P})}}\big)^2+\frac{1}{n}tr(\mathbb{V}^{({\mathpzc P})})\Big\},\\
&{\mathpzc P}^*:= \arg\min_{{\mathpzc P} \in {\cal M}} |{\mathpzc P}|,\,\,\text{where}\,\, {\cal M} :=\arg\min_{\mathpzc P} E\big(Y-{\bf X}^{({\mathpzc P})'}{{\bs \beta}^{({\mathpzc P})}}\big)^2\,\,\text{and}\,\,|{\mathpzc P}|\,\, \text{denotes}\\& \text{the number of covariates in the model} \,\,{\mathpzc P},\\
&\widehat{\pi^*}(n,N) := \arg \min_{\mathpzc P} C^{({\mathpzc P})}(n,N).
\end{align*}  
The following proposition shows that the first two sets defined above by $\arg\min$ converge to the third, which is  a singleton.  Note that ${\mathpzc P}^*$ is the best linear model in the sense of being the most parsimonious model minimizing the expected square of the projection error $Y-{\bf X}^{({\mathpzc P})'}{{\bs \beta}^{({\mathpzc P})}}$. We deal with the convergence of $\widehat{\pi^*}(n,N)$ in Proposition \ref{prop:cons}.
\begin{proposition}  
	\label{prop:single}
	Suppose that the conditions of Theorem \ref{thm:o1n} hold. Then \\
	\noindent {\rm (i)} 
	Any two sequences in $\pi^*(n)$ and ${\mathpzc P}^*(n)$ are equally good, that is, any sequence of models in $\pi^*(n)$ is adequate in the sense of Section \ref{sec:Neighborhoods}.

	\noindent {\rm (ii)}  The set ${\mathpzc P}^*$ is a singleton, and the sets   $\pi^*(n)$ and ${\rm {\mathpzc P}}^*(n)$ converge to the singleton ${\mathpzc P}^*$ as $n \rightarrow \infty$. 
\end{proposition}
The proof shows that essentially $\cal M$ is a singleton; that is, besides ${\mathpzc P}^*$, $\cal M$ may only contain  models having the same covariates and  regression  coefficients as those of ${\mathpzc P}^*$, and further covariates whose coefficients vanish. 
Note that since the number of models is finite, it follows that ${\mathpzc P}^*(n)=\pi^*(n)={\mathpzc P}^*$ for large enough $n$; that is, the same model ${\mathpzc P}^*$ minimizes both $R(n,{\mathpzc P})$ and $AR(n,{\mathpzc P})$.
The model ${\mathpzc P}^*$ is the minimal best linear predictive model that one would ideally use if the projection coefficients ${\bs \beta}^{({\mathpzc P})}$ were known. 

The next proposition shows that minimizing the statistic $C^{({\mathpzc P})}(n,N)$  leads to correct selection asymptotically, that is, to selecting 
the model that minimizes the prediction error $R(n,{\mathpzc P})$ with probability converging to 1.
\begin{proposition}\label{prop:cons} 
	Under the conditions of Theorem \ref{thm:o1n}, with  both $n,N \rightarrow \infty$,$  $ and  $n/N \to 0$, we have $P\big( \widehat{\pi^*}(n,N)={{\mathpzc P}}^*(n)\big) \to 1$.
\end{proposition}
The proof is given in the Appendix, where we also show that the condition $n/N \to 0$ is necessary. The case $n=N$ (with nonrandom covariates) corresponds to the standard Mallows $C_p$, which is inconsistent; more specifically, it is well known that for $n=N$, the choice $\widehat{\pi^*}(n,N)$ may lead to models $\mathpzc Q$ that strictly contain ${\mathpzc P}^*$; see, e.g., \citet{Nishii}. 
The equality $\widehat{\pi^*}(n,N)= {\mathpzc P}^*(n)$, which holds for large enough $n$ and $N$ with high probability,
implies that $\widehat{\pi^*}(n,N)$ is a singleton (by Proposition \ref{prop:single} (ii)), and that selecting a model according to the statistic $\widehat{\pi^*}(n,N)$ yields a model that minimizes the prediction error. Furthermore,  the choice of a model by $\widehat{\pi^*}(n,N)$ leads asymptotically to the choice
of ${\mathpzc P}^*$, the smallest model in terms of the number of covariates in $\cal M$, that is, the most parsimonious model ${\mathpzc P}$ that minimizes $E\big(Y-{\bf X}^{({\mathpzc P})'}{{\bs \beta}^{({\mathpzc P})}}\big)^2$. This property is often referred to as consistency; see, e.g., \citet{Shao}.

In the case of fixed $n$, Equation \eqref{eq:further2}
readily implies that  
$C^{({\mathpzc P})}(n,N) - AR(n,{\mathpzc P}) 
=O_p(1/\sqrt{N})$.
Therefore, as $N$ goes to infinity, the left-hand side converges to zero (at a rate of $1/\sqrt{N}$), implying 
\begin{proposition}\label{prop:nfixed} 
	Under the condition of Theorem \ref{thm:o1n}, we have for any fixed $n$,
	$P\Big( \widehat{\pi^*}(n,N) \subseteq \pi^*(n)\Big){\stackrel{N \to \infty}{\longrightarrow}} 1$.
\end{proposition}
In words, Proposition \ref{prop:nfixed} says that a model that minimizes $C^{({\mathpzc P})}$ will minimize $AR(n,{\mathpzc P})$ with high probability for fixed $n$ and a suitably large $N$. Proposition \ref{prop:single} (i)  asserts that minimizing $AR(n,{\mathpzc P})$  by $\pi^*(n)$  is close to minimizing $R(n,{\mathpzc P})$ by ${\mathpzc P}^*(n)$, which is our goal.

\section{Several datasets} \label{sec:Jind}
Our main focus is on the case where several regression datasets are observed. We first discuss the case where we observe datasets   from all the regressions of interest, and then, in Section \ref{sec:supersuper}, we consider a hierarchical situation where the data consist of a random sample of  regression datasets from a structured collection of regression models.
\subsection{Model selection observing  all regressions}\label{subsec:mallowssev}
We consider a population of distributions
${\cal G}=\{G_j:j=1,\ldots,{\cal K}\}$ with $\cal J = \cal K < \infty$, that is, the training set comprises of all regression datasets in the population. Thus, we observe data $D_j=\{({\bf X}_{ij},Y_{ij})\sim G_j,\,\,i=1,\ldots,N_j\}$, $j=1,\ldots,\cal J$, and ${\bf X}_{ij} \in \mathbb{R}^d$.  

For a given $n$, the goal is to select a common set of covariates $\mathpzc P$ to be used for prediction of the response $Y$ from ${\bf X}={\bf X}^{(\mathpzc P)}$ (the subvector with coordinates in $\mathpzc P$) for each  individual distribution $G_j$ from the population, or equivalently for a random $G_J$, see below \eqref{eq:rtilda}, where the coefficients $\widehat{{\bs \beta}}_{j,n}^{(\mathpzc P)}$, which are allowed to vary with $j$, are estimated with  a sample of size $n$. The relevant prediction error for this task is \eqref{eq:rtilda} below.
When predicting for individual $j$,  it may be natural to set $n=N_j$. However, other values of $n$ may be of interest in studying the contribution of covariates as a function of the sample size. Later (in Section \ref{sec:supersuper}), we use the $\cal J$ datasets as a training set for choosing a model to predict for any out-of-sample  $G_J$ on the basis of $n$ future observations, where $n$ is not determined in advance  since $J$ is not in the training set. In this case we use the chosen subset of covariates, and estimate its parameters on the basis of a dataset of size $n$ from $G_J$. The value of $n$ may vary, being the size of the dataset $G_J$.

Let ${\bf X}:={\bf X}^{(\mathpzc P)} \in {\mathbb{R}}^p$, where for now ${\mathpzc P}$ and its size   $p$ are suppressed in the notation. For each $j$ and generic observation  $({\bf X},Y)$ from the distribution $G_j$, we define 
\begin{equation*}%\label{eq:beta} 
{\bs \beta}_j:=\arg\min_{\bs \beta} E_{G_j}(Y-{\bf X}'{\bs \beta})^2 ={\mathbb{Q}}_j^{-1}E_{G_j}({\bf X}Y);
\end{equation*}
see \eqref{eq:beta}, where $\mathbb{Q}_j:=E_{G_j}({\bf X}{\bf X}')$. 
Assuming  finite fourth moments, we have for a sample size $n \to \infty$, for each $j$, as in \eqref{eq:CLT},
\begin{equation*}\label{eq:asympbeta}
\sqrt{n}(\widehat{\bs \beta}_{j,n}-{\bs \beta}_j) \tendd N({\bf 0}, \mathbb{Q}_j^{-1}\mathbb{W}_j\mathbb{Q}_j^{-1}) \,\,
{\rm where  \,} \,\, \widehat{\bs \beta}_{j,n}:=(\mathbb{X}_{j,n}'\mathbb{X}_{j,n})^{-1}\mathbb{X}_{j,n}'{\bf Y}_{j,n},
\end{equation*}
${\mathbb{X}}_{j,N}$ and ${\bf{Y}}_{j,N}$ are the $j$th versions of ${\mathbb{X}}_{N}$, and ${\bf Y}_N$,
and  $\mathbb{W}_j:=E_{G_j}({\bf X}{\bf X}'e^2)$, 
a $p \times p$ matrix, assumed to be positive definite.
We further use the notation $\mathbb{V}_j$  for the $j$th version of $\mathbb{V}$, that is, when expectations are taken with respect to $G_{j}$, and similar notation when $N=N_j$ observations are used for the estimators $\widehat{\mathbb{Q}}_{j,N}$, $\widehat{\mathbb{V}}_{j,N}$, and $\widehat{\mathbb{W}}_{j,N}$ instead of $\widehat{\mathbb{Q}}_{N}$, $\widehat{\mathbb{V}}_{N}$, and $\widehat{\mathbb{W}}_{N}$.

We consider prediction for a random individual 
regression dataset of size $n$ from the population ${\cal G}$, based on a model, that is, a subset of covariates $\mathpzc P$. As above we suppress $\mathpzc P$ and write ${\bf X}$ and ${\bs \beta}$ rather than ${\bf X}^{(\mathpzc P)}$ and ${\bs \beta}^{\mathpzc P}$, etc. The relevant prediction error (see \eqref{eq:rtild} and around for a discussion) is
\begin{equation}\label{eq:rtilda}
{\bf R}(n,\mathpzc P):=\frac{1}{{\cal J}} \sum_{j=1}^{\cal J} R_j(n,\mathpzc P):=\frac{1}{{\cal J}} \sum_{j=1}^{\cal J} E_{G_{j}}(Y-{\bf X}^{'} \widehat{\bs \beta}_{j,n})^2,
\end{equation}
where $({\bf X},Y)\sim G_j$ independently of $\widehat{\bs \beta}_{j,n}$, and the expectation on the right-hand side of \eqref{eq:rtilda} is also applied to $\widehat{\bs \beta}_{j,n}$. The risk ${\bf R}(n,\mathpzc P)$ can be interpreted as an expectation over $G_J$ for a uniform choice of a single $J \in \cal G$ or equivalently, as the  risk per task average for the multi-task of predicting for all $G_j \in \cal G$ if all datasets sizes (or task size) were $n$.  In ${\bf R}(n,\mathpzc P)$ above and similar expressions below, we suppress the number of datasets $\cal J$. 
In the case that any of the distributions $G_j$ involves discrete covariates, we replace $E_{G_{j}}(Y-{\bf X}^{'} \widehat{\bs \beta}_{j,n})^2$ by a conditional expectation as in Section \ref{sec:BG}, where the conditioning is on a set whose complement is exponentially small. In the definition given in Equation \eqref{eq:rtild}, \eqref{eq:rtilda}, and others below we use boldface letters when ${\cal J}>1$.
Next define
\begin{equation}\label{eq:AR_J}
{\bf AR}(n,{\mathpzc P}):=\frac{1}{{\cal J}} \sum_{j=1}^{\cal J} AR_j(n,{\mathpzc P}):=\frac{1}{{\cal J}} \sum_{j=1}^{\cal J} \Big\{ E_{G_{j}}(Y-{\bf X}' {\bs \beta}_{j})^2  +\frac{tr({\mathbb{V}}_{j})}{n}\Big\}.
\end{equation}

Using \eqref{eq:approx_pred3} we have 
\begin{equation}\label{eq:approx_predJJ}
{\bf R}(n,{\mathpzc P})={\bf AR}(n,{\mathpzc P})+O{(1/n^{3/2})}.
\end{equation}
Set 
\begin{align}\label{eq:JCPP}
{C}_j^{({\mathpzc P})}(n,N_j):=&\frac{1}{N_j}||{\bf Y}_{j,N_j}-\mathbb{X}_{j,N_j} \widehat{\bs \beta}_{j,N_j}||^2 + {tr}(\widehat{\mathbb{V}}_{j,N_j}) (1/n + 1/N_j),\,\,\,\text{and}\nonumber\\& {\bf C}^{({\mathpzc P})}(n,{\bf N}):=\frac{1}{{\cal J}}\sum_{j=1}^{\cal J} C_j^{({\mathpzc P})}(n,N_j),
\end{align} 
where ${\bf N}=(N_1,\ldots,N_{\cal J})$. 
We define the jackknife bias-corrected ${\bf C}^{({\mathpzc P})}$  by 
\begin{equation}\label{eq:CpJackJ}
\pmb{\mathbb{C}}^{({\mathpzc P})}(n,{\bf N}):=\frac{1}{{\cal J}}\sum_{j=1}^{\cal J} \mathbb{C}_j^{({\mathpzc P})}(n,N_j),
\end{equation}
where $\mathbb{C}_j^{({\mathpzc P})}(n,N_j)$ is the bias-corrected ${C}_j^{({\mathpzc P})}(n,N_j)$; see  \citet{Efron}, Equation (2.8), for a precise definition of the jackknife correction we use.

Theorem \ref{prop:C_pj1} below parallels Theorem \ref{thm:o1n} concerning the error of   ${\bf C}^{({\mathpzc P})}(n,{\bf N})$ as an estimator of ${\bf AR}(n,\mathpzc P)$.

\begin{theorem}\label{prop:C_pj1}
	Suppose that the conditions of Theorem \ref{thm:o1n} are satisfied when $({\bf X},Y)\sim G_j$ for each $j=1,\ldots,{\cal J}$. Then,
	\begin{align*}&{\bf AR}(n,{\mathpzc P})-{\bf C}^{({\mathpzc P})}(n,{\bf N})\\&= \frac{1}{{\cal J}}\sum_{j=1}^{{\cal J}} {\cal E}_{j,N_j}+\frac{1}{n{\cal J}}\sum_{j=1}^{{\cal J}} \left\{tr({\mathbb{V}}_{j})-tr(\widehat{\mathbb{V}}_{j,N_j})\right\} +\frac{1}{{\cal J}}\sum_{j=1}^{{\cal J}}o_p\Big(\frac{1}{{N}_j}\Big),
	\end{align*}
	where ${\cal E}_{j,N_j}$ is the $j$th version of ${\cal E}_N$  defined in \eqref{eq:CAPEPS}, and the $o_p$ terms do not depend on $n$.
	
	Moreover, assume that $\lim_{{\bf N}\to \infty} N_1/N_j:=a_j$ exists for all $j$, where $0<a_j <\infty$; then 
	$$\sqrt{ N_1} \big\{{\bf AR}(n,{\mathpzc P})- {\bf C}^{({\mathpzc P})}(n,{\bf N}) \big\}  \tendd N(0,\uptau_{\cal J}^2),$$
	as ${\bf N} \to \infty$, where $\uptau_{\cal J}^2= \frac{1}{{\cal J}^2}\sum_{j=1}^{\cal J} a_j\uptau_j^2$ and  $\uptau^2_j$  is the asymptotic variance under $G_j$ as in Theorem \ref{thm:o1n}, Equation \eqref{eq:further2}.
\end{theorem}
Notice that if $\uptau^2_j$ and $a_j$ are bounded (in $j$), then the asymptotic  variance of $\sqrt{ N_1} \left\{{\bf AR}(n,{\mathpzc P})- {\bf C}^{({\mathpzc P})}(n,{\bf N}) \right\}$ decreases like $1/{\cal J}$, which means that the error is decreasing in ${\cal J}$. Theorem \ref{prop:C_pj1} and \eqref{eq:approx_predJJ}  imply properties of ${\bf C}^{({\mathpzc P})}(n,{\bf N})$ as an estimator of ${\bf R}(n,{\mathpzc P})$ as discussed next.

\subsection{Consistency}\label{subsec:conssev}

Analogously to the definitions in Section \ref{sec:approx}, where now the optimal sets of the multi-task problem are denoted using boldface, define 
\begin{align}\nonumber& \contour[2]{black}{$\mathpzc{P}$}^*(n) :=\arg \min_\mathpzc P {\bf R}(n,\mathpzc P)=\arg \min_p \sum_{j=1}^{\cal J} E_{G_{j}}(Y-{\bf X}^{(\mathpzc P)'} \widehat{\bs \beta}^{(\mathpzc P)}_{j,n})^2,\\
\label{eq:pi_J} &{\bs \pi}^*(n)  := \arg \min_\mathpzc P {\bf AR}(n,\mathpzc P)=\arg \min_\mathpzc P \sum_{j=1}^{\cal J}\Big\{ E_{G_{j}}(Y-{\bf X}^{({\mathpzc P})'} {\bs \beta}^{(\mathpzc P)}_{j})^2  +\frac{tr({\mathbb{V}}^{(\mathpzc P)}_{j})}{n}\Big\},\\ \nonumber &\contour[2]{black}{$\mathpzc{P}$}^*:= \arg\min_{\mathpzc P \in \bs{\cal M}} |\mathpzc P| \text{ where } \bs{\cal M} :=\arg\min_\mathpzc P\sum_{j=1}^{\cal J} E_{G_{j}}(Y-{\bf X}^{(\mathpzc P)'} {\bs \beta}^{(\mathpzc P)}_{j})^2 ,\\ \nonumber  &\widehat{{\bs \pi}^*}(n,{\bf N}) := \arg \min_\mathpzc P {\bf C}^{(\mathpzc P)}(n,{\bf N}).
\end{align}
The next result is similar to Proposition \ref{prop:single}, with essentially by the same proof. The notions  \textit{equally good} and \textit{adequate} are the same as that of Section \ref{sec:Neighborhoods}.
\begin{proposition}\label{prop:Cpj}  
	\label{prop:several}
	Suppose that the conditions of the first part of Theorem \ref{prop:C_pj1} hold. Then \\
	\noindent {\rm (i)} 
	Any two sequences in ${\bs \pi}^*(n)$ and $\contour[2]{black}{$\mathpzc{P}$}^*(n)$ are equally good, that is, any sequence of models in ${\bs \pi}^*(n)$ is adequate.\\
	{\rm (ii)}
	The set $\contour[2]{black}{$\mathpzc{P}$}^*$ is a singleton and the sets $\contour[2]{black}{$\mathpzc{P}$}^*(n)$ and  ${\bs \pi}^*(n)$ converge to the singleton $\contour[2]{black}{$\mathpzc{P}$}^*$ as $n \rightarrow \infty$.
\end{proposition}

The following proposition generalizes Propositions \ref{prop:nfixed} and \ref{prop:cons} to ${\cal J} >1$.
Here we  consider a uniform bound \eqref{eq:Cbound}. Technically, the constant $C$ provides a measure of the notion of  ``sufficiently homogeneous" of Section \ref{sec:forset} when referring to the set of distributions $\cal G$;  informally we  mean that the regression datasets have enough in common to justify common subsets for prediction.

\begin{proposition}\label{prop:cons_J_fixed}
	\begin{enumerate}
		\item
		Assume that the conditions of the first part of Theorem \ref{prop:C_pj1} hold. Then
		for fixed $n$ we have  
		\[
		\lim_{{\bf N}\to \infty} P\left( \widehat{{\bs \pi}^*}(n,{\bf N}) \subseteq {\bs \pi}^*(n)\right)  = 1.
		\]
		\item Let $n/N_j$ be bounded for all $j=1,\ldots,\cal J$, and let $C$ be a constant satisfying
		for all $j$ and $\mathpzc P$  
		\begin{equation}\label{eq:Cbound}
		n/N_j , \lambda_{max}(\mathbb{W}_j^{(\mathpzc P)}), 1/\lambda_{min}(\mathbb{W}_j^{(\mathpzc P)}),\lambda_{max}(\mathbb{Q}_j^{(\mathpzc P)}) \le C.
		\end{equation}
		Then 
		\[
		\lim \inf_{n/N_j \le C, \,n\to \infty,\,{\bf N} \to \infty} P\left( \widehat{{\bs \pi}^*}(n,{\bf N})=\contour[2]{black}{$\mathpzc{P}$}^*(n)\right)  \ge 1-K_C/{\cal J},\] 
	 where $K_C$  depends only on $C$.
	\end{enumerate}
\end{proposition}
The existence of $C$  follows from the assumption on $n/N_j$ since only a finite number of bounded terms appear in \eqref{eq:Cbound} besides $n/N_j$.

Part 1 of the above proposition extends Propositions \ref{prop:nfixed}.
Part 2 extends  Proposition \ref{prop:cons}; however, a stronger condition was needed before, namely that $n/N \to 0$, to obtain consistency. Here, we obtain approximate consistency for large ${\cal J}$, assuming only that $n/N_j$ is bounded, along with the other terms in \eqref{eq:Cbound}.
This is useful since in our application it is natural to consider the possibility that $n=N_j$.

The rate $K_C   /\cal J$ in the theorem was achieved by using Chebyshev's inequality. Since under our assumptions all moments are bounded, a similar argument using a bound on $2m$ moments leads in the same way to the rate $K_{C,m}/{\cal J}^m$, where $K_{C,m}$ depends also on $m$, and with further effort,  a large deviation rate (in ${\cal J}$) can be achieved.

\subsection{A population of distributions}
\label{sec:supersuper}
We now consider the situation where we have a sample of  $\cal J$ regression datasets from a given, finite or infinite, population of such datasets, and we are interested in predictions for a random (possibly out-of-sample) further regression or several regressions from the same population. In terms of the application considered in this paper, this situation corresponds to the case that we have a training  sample of ${\cal J}$ doctors out of many more, and our goal is to select a subset of covariates to be used to predict service durations for a random doctor from the population who may not be in the training sample.

Formally, let $(\Theta,{\mathscr T},{\mathscr P})$ be a probability space and let $\{G_\theta:\theta \in \Theta\}$ be a family of distributions; see, e.g.,  \citet{Cinlar}, Chapter VI for a formulation of random measures.  Let $\{\theta_1,\ldots,\theta_{\cal J}\}$ be a sample  where $\theta_j \sim {\mathscr P}$, and as in Section \ref{subsec:mallowssev}   we observe a training set consisting of datasets $D_j=\{({\bf X}_{ij},Y_{ij})\sim G_j,\,\,i=1,\ldots,N_j\}$, $j=1,\ldots,\cal J$, where $G_j$ stands for $G_{\theta_j}$.
Given $\theta \in \Theta$, we consider  $D=\{({\bf X}_{i},Y_{i})\sim G_\theta, \, i=1,\ldots,n\}$.
For any function $f$ for which the conditional expectation $E_{G_\theta}f(D)$ of $f(D)$ given $G_\theta$ is well defined, we assume that so is 
$E_{\mathscr P} E_{G_\theta}f(D)=\int E_{G_\theta}f(D) {\mathscr P}(d\theta)$, where the outer expectation is over $\theta \sim {\mathscr P}$.

We now fix a set of covariates ${\mathpzc P}$, which is suppressed in most of the notation as before.
For any $G_\theta$ we define $\widehat{\bs \beta}_{\theta,n}$ to be the least squares estimator for the given dataset  $D$.
If  $G_{\theta}$  is sampled randomly from ${\mathscr P}$
then the population prediction error is defined as
\begin{equation}\label{eq:predrrr}
{\bf R}_{pop}(n,{\mathpzc P}):=\int R_{\theta}(n,{\mathpzc P}){\mathscr P}(d\theta):=\int E_{G_{\theta}}(Y-{\bf X}^{'} \widehat{\bs \beta}_{\theta,n})^2 {\mathscr P}(d\theta),
\end{equation}
where the expectation $E_{G_\theta}$ is over $\widehat{\bs \beta}_{\theta,n}$ and $({\bf X},Y) \sim G_\theta$ that are independent of $\widehat{\bs \beta}_{\theta,n}$.
Let ${\bs \beta}_\theta:=\arg\min_{\bs \beta } E_{G_\theta}(Y-{\bf X}'{\bs \beta} )^2$, 
$\mathbb{Q}_\theta:=E_{G_\theta}({\bf X}{\bf X}')$, $\mathbb{W}_\theta:=E_{G_\theta}({\bf X}{\bf X}'e^2)$, and $\mathbb{V}_\theta:=
\mathbb{W}_\theta \mathbb{Q}^{-1}_\theta$.
As before,
${\bf R}_{pop}(n,{\mathpzc P})$ is approximated  by 
\begin{equation}\label{eq:predrand}
{\bf AR}_{pop}(n,{\mathpzc P}):=\int AR_{\theta}(n,{\mathpzc P}){\mathscr P}(d\theta):=\int \Big\{ E_{G_{\theta}}(Y-{\bf X}' {\bs \beta}_{\theta})^2  +\frac{tr({\mathbb{V}}_{\theta})}{n}\Big\}{\mathscr P}(d\theta),
\end{equation}
where the latter integrand defines $AR_{\theta}(n,{\mathpzc P})$ as in  \eqref{eq:approx_predara}. 
The quantity $AR_{\theta}(n,{\mathpzc P})$, whose estimation was already discussed, is now a random variable, since it depends on $G_{\theta}$ with $\theta \sim {\mathscr P}$; its expectation,  given by \eqref{eq:predrand},  is the basis of our estimation of ${\bf R}_{pop}(n,{\mathpzc P})$ of \eqref{eq:predrrr}. 
Lemma \ref{lem:unif} below generalizes \eqref{eq:approx_pred3}. 

\begin{lemma}\label{lem:unif}
	Suppose that the conditions of Theorem \ref{thm:o1n} hold uniformly in  $\theta \in \Theta$; that is, for each $k$, the $k$th moment with respect to $G_\theta$ of each entry of $\bf X$ and $Y$   is bounded uniformly in $\theta$, and the entries of $(\mathbb{X}_n'\mathbb{X}_n/n)^{-1}$ have third moments with respect to $G_\theta$ that are bounded uniformly in $n$ and $\theta$. 
	Then 
	\[
	{\bf R}_{pop}(n,{\mathpzc P})={\bf AR}_{pop}(n,{\mathpzc P})+O({1}/{n^{3/2}}).
	\]	
\end{lemma}
The lemma clearly holds if $\Theta$ is finite, and in general it follows readily by the uniform boundedness of moments in $\theta$ and the proof of \eqref{eq:approx_pred3} given in the Appendix.  
Recall Lemma \ref{lem:normal}, where we showed that the moment conditions of Theorem \ref{thm:o1n} hold when ${\bf X}$ is a mixture of normals and $\inf_{\Sigma \in \Xi}\lambda_{min}({\Sigma}) >0$. 
For the bound on moments as assumed in  Lemma \ref{lem:unif} to hold uniformly,  it suffices that 
$\inf_\Sigma \lambda_{min}({\Sigma}) >0$, where  now the infimum is over all covariance matrices of all the mixing normal distributions involved in all the distributions $G_\theta$ for all $\theta \in \Theta$. This technical assumption means that the covariates that are taken into account for the model selection are ``bounded away" from multicollinearity. 
For discrete variables we redefine the prediction error by conditioning as in Section \ref{sec:BG}. Theorem \ref{thm:BeGood} extends easily when we assume that all covariates are uniformly bounded in $\theta$, and that  $\lambda_{min}(\mathbb{Q}_\theta)>c$ for some $c>0$, for all $\theta$.

Lemma \ref{lem:unif} suggests that a consistent estimator of ${\bf AR}_{pop}(n,\mathpzc P)$ will lead to selection of an adequate model  in the sense of Section \ref{sec:Neighborhoods}, that is, a model that is as good as the model that minimizes ${\bf R}_{pop}(n,\mathpzc P)$.

Recall the definition of ${\bf AR}(n,\mathpzc P)$ in \eqref{eq:AR_J}; now this quantity is considered random as it is a function of the sampled distributions $G_1,\ldots,G_{\cal J}$. In order to generalize the consistency results of Theorem \ref{prop:C_pj1} to this case, we need to bound ${\bf AR}_{pop}(n,\mathpzc P)-{\bf AR}(n,\mathpzc P)$ as in the lemma below.

\begin{lemma}\label{lem:tild_AR}
	Under the conditions of Lemma \ref{lem:unif},
	\begin{equation}\label{eq:R_j1}
	{\bf AR}_{pop}(n,{\mathpzc P})-{\bf AR}(n,{\mathpzc P}) =O_p(1/\sqrt{{\cal J}})
	\end{equation}
	uniformly in $n$. Moreover, for any fixed $n$, $\sqrt{{\cal J}} \big({\bf AR}_{pop}(n,{\mathpzc P})-{\bf AR}(n,{\mathpzc P})\big)$ is asymptotically normal.
\end{lemma}

We now consider the population versions of the quantities defined in Section \ref{subsec:conssev}.
\begin{align*}&\contour[2]{black}{$\mathpzc{P}$}^*_{pop}(n) :=\arg \min_\mathpzc P {\bf R}_{pop}(n,\mathpzc P)=\arg \min_\mathpzc P \int E_{G_{\theta}}(Y-{\bf X}^{(\mathpzc P)'} \widehat{\bs \beta}^{(\mathpzc P)}_{\theta,n})^2 {\mathscr P}(d\theta),\\ &{\bs \pi}_{pop}^*(n)  := \arg \min_\mathpzc P {\bf AR}_{pop}(n,\mathpzc P)\\& \qquad \qquad \qquad \qquad=\arg \min_\mathpzc P \int \Big\{ E_{G_{\theta}}(Y-{\bf X}^{(\mathpzc P)'} {\bs \beta}^{(\mathpzc P)}_{\theta})^2  +\frac{tr({\mathbb{V}}^{(\mathpzc P)}_{\theta})}{n}\Big\}{\mathscr P}(d\theta),\\& \contour[2]{black}{$\mathpzc{P}$}^*_{pop}:=\arg\min_{\mathpzc P \in \bs{\cal M}_{pop}} |\mathpzc P| \text{ where } \bs{\cal M}_{pop} :=\arg\min_\mathpzc P\int E_{G_{\theta}}(Y-{\bf X}^{(\mathpzc P)'} {\bs \beta}^{(\mathpzc P)}_{\theta})^2 {\mathscr P}(d\theta).
\end{align*}
$\widehat{{\bs \pi}^*}(n,{\bf N}) \text{  and  }  {\bs \pi}^*(n) $ are defined as in \eqref{eq:pi_J}, however the fact that  now the $G_j$'s are random adds randomness to $\widehat{{\bs \pi}^*}(n,{\bf N})$, and makes  ${\bs \pi}^*(n)$ a random variable.
\vspace{.15cm}

Proposition \ref{prop:popolo} parallels  Proposition \ref{prop:cons_J_fixed}; it shows consistency properties  of  $\widehat{{\bs \pi}^*}(n,{\bf N})$, as defined in \eqref{eq:pi_J} using \eqref{eq:JCPP}. 
Below, the probability $P$ is obtained by first conditioning on $\theta_1,\ldots,\theta_{\cal J}$, and then unconditioning by taking expectation over $\theta_1,\ldots,\theta_{\cal J}$ with respect to the product measure ${\mathscr P}^{\cal J}$.

\begin{proposition}
	\label{prop:popolo}
	Assume that the conditions of Lemma \ref{lem:unif} hold, and in addition that $n/N_j , \lambda_{max}(\mathbb{W}_\theta^{(\mathpzc P)}),$  $1/\lambda_{min}(\mathbb{W}_\theta^{(\mathpzc P)}), \lambda_{max}(\mathbb{Q}_\theta^{(\mathpzc P)}) \le C$ for all $\theta$ and $\mathpzc P$.
	\begin{enumerate}
		\item
		When $n$ is fixed, 
		\[
		\lim \inf_{{\bf N}\to \infty} P\Big( \widehat{{\bs \pi}^*}(n,{\bf N}) \subseteq {\bs \pi}_{pop}^*(n)\Big)  \ge 1 - \frac{K_C}{\cal J}.
		\] 
		
		\item Letting $n, {\bf N}\to \infty$,
		\[
		\lim \inf_{n/N_j \le C, \,n\to \infty, \,{\bf N} \to \infty} P\left( \widehat{{\bs \pi}^*}(n,{\bf N})=\contour[2]{black}{$\mathpzc{P}$}_{pop}^*(n)\right)  \ge 1-\frac{K_C}{\cal J},
		\]
		where $K_C$ depends only on $C$.
	\end{enumerate}
\end{proposition}

The proof of Proposition \ref{prop:popolo} also shows that $\contour[2]{black}{$\mathpzc{P}$}^*_{pop}$ is a singleton and $\contour[2]{black}{$\mathpzc{P}$}^*_{pop}(n)$ and ${\bs \pi}^*_{pop}(n)$ converge to it when $n\to \infty$.

\section{GENO}\label{sec:GENO}
\subsection{Definition of GENO}

Given a model (i.e., a set of covariates) $\mathpzc P$ with  coefficients estimated by a sample of $n$ observations, we can say that it is equivalent to  another model $\mathpzc Q$ with $m$ observations if their expected prediction errors satisfy $R(m,\mathpzc Q) = R(n, \mathpzc P)$. 
Using the  approximation  $AR(n,\mathpzc P)$ to $R(n, \mathpzc P)$ given in \eqref{eq:approx_predara}, \eqref{eq:AR_J}, and \eqref{eq:predrand} for each of the  scenarios we consider, 
we define GENO by 
\begin{equation} \label{eq:nostarsi} 
\text{GENO}(n;\mathpzc P,\mathpzc Q) :=\big\{ m : AR(m,\mathpzc Q)=AR(n,\mathpzc P) \big\}.
\end{equation} 
If $AR(m,\mathpzc Q)> (<) AR(n,\mathpzc P)$ for all $m$,  we set GENO$(n,\mathpzc P,\mathpzc Q)=\infty (0)$, indicating that model $\mathpzc P$ with $n$ observations is better than model $\mathpzc Q$ with any number of observations (model $\mathpzc Q$ with any number of observations is better than $\mathpzc P$ with $n$). A direct calculation shows that for  ${\cal J}=1$ we have 
$$\text{GENO}(n;\mathpzc P,\mathpzc Q)=tr(\mathbb{V}^{(\mathpzc Q)})\big\{AR(n,\mathpzc P)-AR(n,\mathpzc Q)+\frac{1}{n}tr\mathbb{V}^{(\mathpzc Q)}\big\}^{-1}.$$ 
For ${\cal J}>1$ with ${\bf AR}$  defined in \eqref{eq:AR_J} we have 
\begin{equation*} 
\text{GENO}(n;\mathpzc P,\mathpzc Q) =\Big[\frac{1}{{\cal J}}\sum_j tr(\mathbb{V}_j^{(\mathpzc Q)})\Big]\Big\{{\bf AR}(n,\mathpzc P)-{\bf AR}(n,\mathpzc Q)+\frac{1}{{\cal J}n}\sum_j tr(\mathbb{V}_j^{(\mathpzc Q)})\Big\}^{-1}.
\end{equation*} 
For the case of \eqref{eq:predrand}, $j$ is replaced by $\theta$, and the averages by integrals ${\mathscr P}(d\theta)$.

{GENO}$(n;\mathpzc P,\mathpzc Q)=m$ means that model $\mathpzc P$ with $n$ observations is equivalent in terms of  expected prediction error to model $\mathpzc Q$ with $m$ observations.
Note that the larger {GENO}$(n;\mathpzc P,\mathpzc Q)$ is, the better model $\mathpzc P$ (with $n$ observations) is relative to model $\mathpzc Q$. 
For each model $\mathpzc P$ and sample size $n$, we define  
\begin{equation}
\label{eq:ridell1} {\rm GENO}(n,\mathpzc P)=\min_{\mathpzc R}{\rm GENO}(n;\mathpzc P,\mathpzc R),
\end{equation}
where the minimum is over all subsets of covariates $\mathpzc R$.
It follows that the inequality ${\rm GENO}(n,\mathpzc P) \le n$ holds always, where equality means that $\mathpzc P$ is the best model for $n$ observations, as no other model can achieve the same prediction error with fewer observations. On the  other hand,  ${\rm GENO}(n,\mathpzc P) = m <n$ means that there is a model that achieves, with $m <n$ observations, the same prediction error as $\mathpzc P$ with $n$ observations.  Thus,  small values of GENO$(n,\mathpzc P)$ suggest considering another model. 
By the monotonicity of $AR(n,\mathpzc P)$ in $n$, if the inequality ${\rm GENO}(n;\mathpzc P,\mathpzc R) \ge {\rm GENO}(n;\mathpzc Q,\mathpzc R)$ holds for some model $\mathpzc R$, then it holds all $\mathpzc R$. This readily implies 
\begin{multline}\label{eq:implicationsGENO}
AR(n,\mathpzc P) \le AR(n,\mathpzc Q) \Leftrightarrow {\rm GENO}(n;\mathpzc P,\mathpzc R) \ge {\rm GENO}(n;\mathpzc Q,\mathpzc R)\,\, \text{ for all}\,\, {\mathpzc R} \\
 \Leftrightarrow {\rm GENO}(n,\mathpzc P) \ge {\rm GENO}(n,\mathpzc Q).
\end{multline}

\subsection{Estimation of GENO}\label{subsec:estGENO}

In the case ${\cal J}=1$, \eqref{eq:further2}  shows the consistency of $C^{(\mathpzc P)}(n,N)$ as an estimator of $AR(n,\mathpzc P)$
for fixed $n$ as $N \to \infty$. In view of  \eqref{eq:nostarsi} we define a consistent estimator of ${\rm GENO}(n;p,q)$ by
$$\widehat{{\rm GENO}}(n;p,q):=\big\{m: C^{(\mathpzc Q)}(m,N)=C^{(\mathpzc P)}(n,N)\big\}.$$ 
To avoid cumbersome notation we suppress $N$ in $\widehat{{\rm GENO}}$.
Using \eqref{eq:cpp} we obtain,  as before,
\begin{equation}\label{GENO_AIC}
\widehat{{\rm GENO}}(n;\mathpzc P,\mathpzc Q) = {tr}(\widehat{\mathbb{V}}^{(\mathpzc Q)}_N)
\Big\{C^{(\mathpzc P)}(n,N)-C^{(\mathpzc Q)}(n,N) +\frac{{tr}(\widehat{\mathbb{V}}^{(\mathpzc Q)}_N)}{n}\Big\}^{-1},
\end{equation}
setting it to be $\infty$ if the expression in curly brackets is negative or zero.  

In the case of ${\cal J}>1$ datasets of Section \ref{subsec:mallowssev}, or in the population case of Section \ref{sec:supersuper}, the above expression \eqref{GENO_AIC} remains unchanged except that now $C^{(\mathpzc P)}(n,N)$ is replaced by ${\bf C}^{(\mathpzc P)}(n,{\bf N})$ defined in \eqref{eq:JCPP}, and $\widehat{\mathbb{V}}^{(\mathpzc Q)}_N$ is replaced by $\frac{1}{{\cal J}}\sum_j tr(\widehat{\mathbb{V}}_{j,N_j}^{(\mathpzc Q)})$.  We can also define the estimator of \eqref{GENO_AIC} in terms of the jackknife bias-corrected $\pmb{\mathbb{C}}^{(\mathpzc P)}(n,{\bf N})$ of 
\eqref{eq:CpJackJ}. This is done in estimating GENO in Section \ref{sec:ms}.
The results below hold in the same way for all these cases.
Similarly to \eqref{eq:ridell1}, we define the statistic
$$\widehat{\rm GENO}(n,\mathpzc P):=\min_{\mathpzc R}\widehat{\rm GENO}(n;\mathpzc P, \mathpzc R),$$
which is an
estimate the minimal number of observations  required by the best competing model to achieve the same prediction error as model $\mathpzc P$ with sample size $n$.

As in \eqref{eq:implicationsGENO}, we have 
\begin{multline*}%\label{eq:implicationsGENO}
C^{(\mathpzc P)}(n,N) < C^{(\mathpzc Q)}(n,N) \Leftrightarrow \widehat{{\rm GENO}}(n;\mathpzc P,\mathpzc R) \ge \widehat{{\rm GENO}}(n;\mathpzc Q,\mathpzc R)\,\, \forall \mathpzc R \\
 \Leftrightarrow \widehat{{\rm GENO}}(n,\mathpzc P) \ge\widehat{{\rm GENO}}(n,\mathpzc Q).
\end{multline*}

The next proposition follows from \eqref{eq:further2} by applying the $\delta$-method to the inverse function in \eqref{GENO_AIC}. In particular, it shows the consistency of $\widehat{{\rm GENO}}(n;\mathpzc P,\mathpzc Q)$ for fixed $n$ as $N \rightarrow \infty$.
\begin{proposition}\label{prop:cons_GENO} 
	Under the conditions of Theorem \ref{prop:C_pj1} (which include the case ${\cal J}=1$), we have for any fixed $n$ 
	\begin{equation*}%\label{eq:cons_GENO}
	\sqrt{N_1}\big(\widehat{{\rm GENO}}(n;\mathpzc P,\mathpzc Q)-{{\rm GENO}}(n;\mathpzc P,\mathpzc Q)\big) \tendd N(0,\upeta^2)\,\, \rm{ as}\,\,  {\bf N} \rightarrow \infty,
	\end{equation*}
	for some $\upeta^2 >0$.
\end{proposition}
The variance $\upeta^2$ is not given explicitly since it is too complicated to be useful, and it can be be computed by the bootstrap.
See Theorem \ref{prop:C_pj1} and the ensuing comment, which show that (under certain conditions) the variance decreases at a rate of $1/\cal J$.

A similar problem is to estimate for a given model $\mathpzc P$ and a certain prescribed prediction error $E$ the sample size $n$ that satisfies $AR(n,\mathpzc P)=E$. When ${\cal J}=1$, using \eqref{eq:approx_predara} this quantity is given by 
$\displaystyle
\frac{tr({\mathbb{V}}^{(\mathpzc P)})}{E- E\big(Y-{\bf X}^{(\mathpzc P)'}{\bs \beta}^{(\mathpzc P)}\big)^2}$ and can be estimated by 
\begin{equation}\label{eq:achieve} \frac{{tr}( {\widehat{\mathbb{V}}_N})}{E-\frac{1}{N}\{||{\bf Y}_N-\mathbb{X}_N^{(\mathpzc P)}\widehat{\bs \beta}_N^{(\mathpzc P)}||^2
	+ {tr}({\widehat{\mathbb{V}}^{(\mathpzc P)}_N})\} }
\end{equation}  (Since $\frac{1}{N}\{||{\bf Y}_N-\mathbb{X}_N^{(\mathpzc P)}\widehat{\bs \beta}_N^{(\mathpzc P)}||^2
+ {tr}({\widehat{\mathbb{V}}^{(\mathpzc P)}_N}) \}$ is an unbiased estimator of $E\big(Y-{\bf X}^{(\mathpzc P)'}{\bs \beta}^{(\mathpzc P)}\big)^2$); the extensions to the cases ${\cal J}>1$ and to the population setup are straightforward.

\section{Simulations}\label{sec:sims}
In this section we evaluate by simulations the prediction error $R(n,\mathpzc P)$, its approximation $AR(n,\mathpzc P)$, and its estimation using $C^{(\mathpzc P)}$. We start with a single dataset  (${\cal J}=1$) and then we consider the case of several datasets. This simple example demonstrates the well-known difficulty involved in model selection for a single given dataset with  methods such as 
Mallows $C_p$, AIC, BIC, as well as our version  $C^{({\mathpzc P})}$. 
In Section \ref{sec:simJ} we compare the case of model selection for one dataset to that of choosing a common model for successful prediction  on the average when we have data from several datasets, that is, a multi-task. Section \ref{sec:compare} compares the prediction error when model selection is done according to $\pmb{\mathbb{C}}^{({\mathpzc P})}$ to the prediction error under other methods. 

\subsection{A single dataset}
\begin{figure}[ht!]
	\begin{center}
		\includegraphics[width=.9\textwidth]{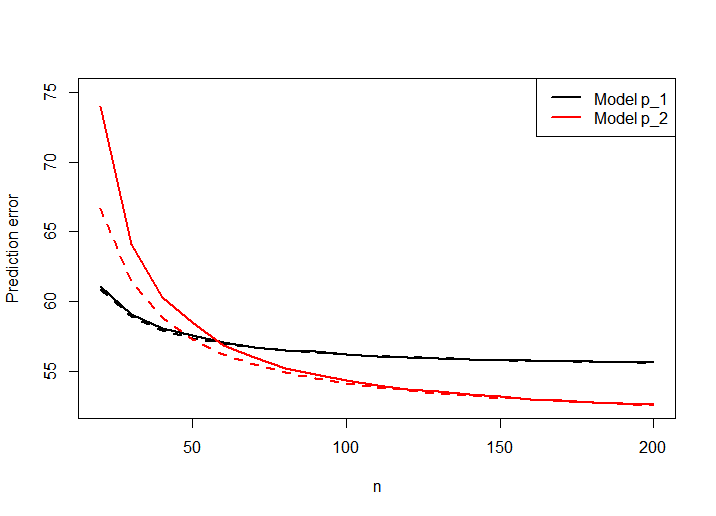}
		\caption{\footnotesize Simulation estimates of $R(n,{\mathpzc P}_1)$ and $R(n,{\mathpzc P}_2)$ (solid line) as well as the approximations AR (dashed line) given in  \eqref{eq:appr}.  } \label{fig:pred}
	\end{center}
\end{figure}

Suppose that the distribution of $({\bf X},Y)$ for ${\bf X} \in \mathbb{R}^{5}$ is given by
\begin{equation}\label{eq:xy}
Y= b_0 + b_1 X_1 + \ldots+ b_{5} X_{5}+ a(X_1^2-1)+\sigma\ve,
\end{equation}
with $X_1,\ldots,X_{5},\,\ve \sim^{iid} N(0,1)$. Setting all models to include the intercept,  there are $2^{5}$ possible submodels; for simplicity, we focus for now on two models consisting of the subsets of covariates ${\mathpzc P}_1 =\{1, X_1\}$, ${\mathpzc P}_2 =\{1, X_1,\ldots, X_5\}$; more explicitly, we have  
model ${\mathpzc P}_1$: $Y=\beta_0+\beta_1X_1+e$ and  model ${\mathpzc P}_2$: $Y=\beta_0+\beta_1X_1+\ldots\beta_5 X_5 +e$. These two models are wrong (as linear conditional expectation function models, see \citet{Hansen} Section 2.15) since the residual $e$ includes the nonlinear term $X_1^2-1$. By the orthogonality of the variables in \eqref{eq:xy}, the projection parameters $\beta_k$ are equal to $b_k$ for these models; see \eqref{eq:beta}.
This is used in computing the first part of $AR(n,{\mathpzc P}_\ell)$ for $\ell=1,2$, and since in this case $\mathbb{Q}=I$, it is also easy to compute $tr(\mathbb{V})$ for each model.
We obtain
\begin{align}\label{eq:appr}
AR(n,{\mathpzc P}_1)&= \sum_{k=2}^{5} b_k^2+ 2a^2 +\sigma^2 +\frac{2(\sum_{k=2}^{5} b_k^2+\sigma^2)+12a^2}{n},\\
AR(n,{\mathpzc P}_2) &=2a^2 + \sigma^2 + \frac{ 6\sigma^2 +20 a^2 }{n};
\end{align}
notice that the above functions do not depend on $b_0,b_1$. 
For a concrete example, we set in \eqref{eq:xy} 
\begin{equation}\label{eq:bs}
(b_0,b_1)=(1,3), ~(b_2,\ldots,b_5)=(1,\ldots,1),~a=1,~\sigma=7.
\end{equation}
Figure \ref{fig:pred} plots $R(n,{\mathpzc P}_\ell)$ (see \eqref{eq:kprederr}-\eqref{eq:kpr2}) (solid lines) and $AR(n,p_\ell)$ (see \eqref{eq:appr}) (dashed line), $\ell=1,2$, as functions of $n$ for the above parameters. We evaluated  $R(n,{\mathpzc P}_\ell)$, where $\ell=1,2$, by a simulation based on $10^3$ repetitions and using the decomposition (see \eqref{eq:predX} and recall that $\mathbb{Q}=I$)
\[
R(n,{\mathpzc P}_\ell)=E\big(Y-({\bf X}^{({\mathpzc P}_\ell)})'{\bs \beta}^{({\mathpzc P}_\ell)}\big)^2 + E\| \hat{\bs \beta}^{({\mathpzc P}_\ell)} - {\bs \beta}^{({\mathpzc P}_\ell)}\|^2; 
\] 
the first expectation can be computed explicitly and the second is evaluated using simulations.
For small $n$, $R(n,{\mathpzc P}_2)$  differs from $AR(n,{\mathpzc P}_2)$, and the approximation improves as $n$ increases. For $n$   smaller than about 50, model ${\mathpzc P}_1$ has a smaller prediction error; for large $n$ model ${\mathpzc P}_2$ is better. This holds approximately for both $R$ and $AR$. This makes sense as models with fewer parameters have a smaller prediction error for small $n$.  The rest of the models are not optimal for any $n$ (this observation is not shown in Figure \ref{fig:pred}). 

Consider GENO as defined in \eqref{eq:nostarsi}. Careful inspection of Figure \ref{fig:pred} shows, for example,  that  GENO$(49;{\mathpzc P}_1,{\mathpzc P}_2)= 49$, which means that in order to achieve the same prediction error as model ${\mathpzc P}_1$ with $n=49$ observations (the value of $n$ where the dashed black line and red the line intersect), model ${\mathpzc P}_2$ requires the same number of observations. Also,
GENO$(60;p_2,p_1)=95$, and therefore, to achieve the same prediction error as model ${\mathpzc P}_2$ with $n=60$, model ${\mathpzc P}_1$ would require 95 observations (the value of $n$ where the dashed black line has the same level as the dashed red line at 60). Since the decrease of $AR(n,{\mathpzc P}_1)$ (the black line) in $n$ is   slow, a small increase in $n$, will result in a much larger value of GENO$(n;{\mathpzc P}_2,{\mathpzc P}_1)$; for example, GENO$(65;{\mathpzc P}_2,{\mathpzc P}_1)=142$. As mentioned before, GENO allows the statistician to compare the cost of additional observations to the cost of measuring additional variables, which may be expensive, or harmful, such as in the case of an invasive medical procedure or imaging that involves radiation. 

By \eqref{eq:achieve}, the numbers of observations for models ${\mathpzc P}_1$ and ${\mathpzc P}_2$ to obtain a prediction error of 59 are about 29 and 39, respectively; i.e., model ${\mathpzc P}_1$ can achieve this prediction error with a sample size that is smaller by 10 observations. On the other hand, for a prediction error of 56, model ${\mathpzc P}_1$ requires 118 observations, while model ${\mathpzc P}_2$ needs only 63 observations.

\begin{figure}[ht!]
	\subfigure[Boxplots of $C^{({\mathpzc P}_1)}(n,N) - C^{({\mathpzc P}_2)}(n,N)$]{%
		\includegraphics[width=0.485\textwidth]{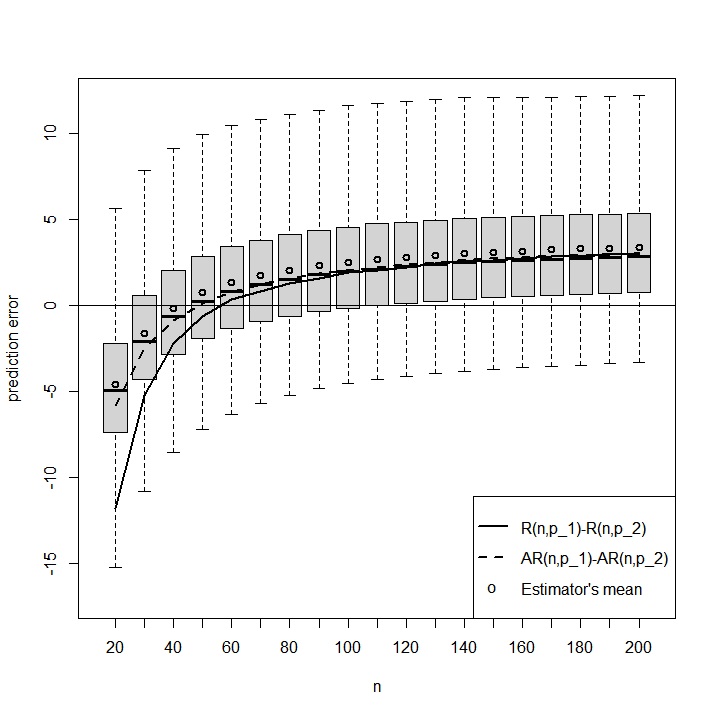}
	}
	\subfigure[Boxplots of $\mathbb{C}^{({\mathpzc P}_1)}(n,N) - \mathbb{C}^{({\mathpzc P}_2)}(n,N)$]{%
		\includegraphics[width=0.485\textwidth]{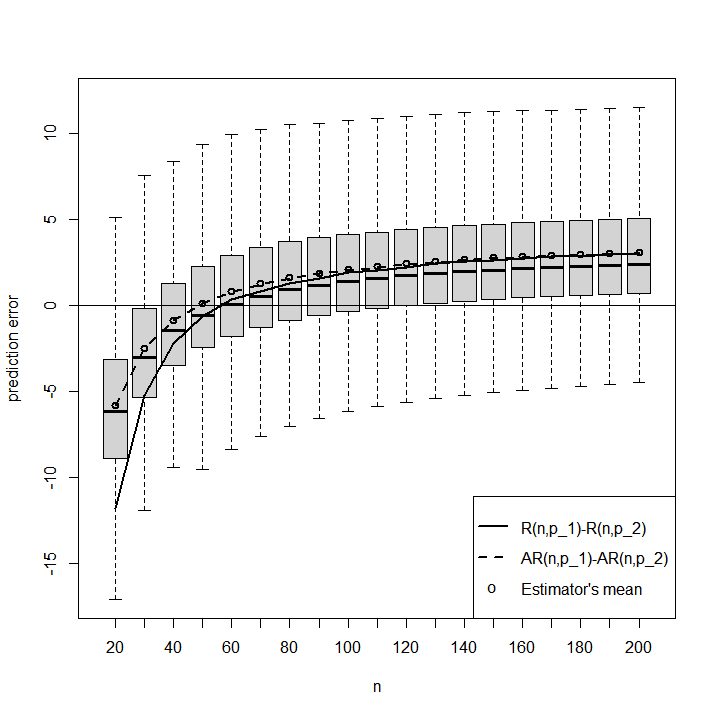}
	}\\
	\caption{\footnotesize (a) Boxplots of the simulation results of $C^{({\mathpzc P}_1)}(n,N) - C^{({\mathpzc P}_2)}(n,N)$, where {\tiny $\circ$} (circle) denotes the mean, and  $R(n,{\mathpzc P}_1)-R(n,{\mathpzc P}_2)$ (respectively, $AR(n,{\mathpzc P}_1)-AR(n,{\mathpzc P}_2)$) is a solid (respectively, dashed) line. (b) Same as (a) for jackknifed version $\mathbb{C}^{({\mathpzc P}_1)}(n,N)-\mathbb{C}^{({\mathpzc P}_2)}(n,N)$. } \label{fig:Cp}
\end{figure}

We now discuss estimation of the prediction error using $C^{({\mathpzc P})}(n,N)$ based on a single dataset of size $N=100$. Figure \ref{fig:Cp} plots $R(n,{\mathpzc P}_1)-R(n,{\mathpzc P}_2)$ (solid line), $AR(n,{\mathpzc P}_1)-AR(n,{\mathpzc P}_2)$ (dashed line), and boxplots of the estimators $C^{({\mathpzc P}_1)}(n,N) - C^{({\mathpzc P}_2)}(n,N)$ on the left-hand side, and  $\mathbb{C}^{({\mathpzc P}_1)}(n,N) - \mathbb{C}^{(p_2)}(n,N)$, the jackknife bias-corrected version, on the right-hand side, based on $10^3$ simulations for each $n=20, 40,\ldots,200$. Their means are given by circles.  We see that the jackknife corrects the bias of $C^{({\mathpzc P}_1)}(n,N) - C^{({\mathpzc P}_2)}(n,N)$ as an estimator of $AR(n,{\mathpzc P}_1)-AR(n,{\mathpzc P}_2)$; see the discussion following \eqref{eq:cpp}.   Recall that the bias itself and the correction  decrease in $n$. The mean of the difference ${C}^{({\mathpzc P}_1)}(n,N) - {C}^{({\mathpzc P}_2)}(n,N)$ and $\mathbb{C}^{({\mathpzc P}_1)}(n,N) - \mathbb{C}^{({\mathpzc P}_2)}(n,N)$ equals 0 at about $n=40$ and  $n=50$, respectively; thus the jackknife leads to correct selection on average since it is optimal to select model ${\mathpzc P}_1$ for about $n \le 50$.

\begin{figure}[h!]
	\subfigure[Probability of selection by $C^{({\mathpzc P})}(n,N)$]{%
		\includegraphics[width=0.48\textwidth]{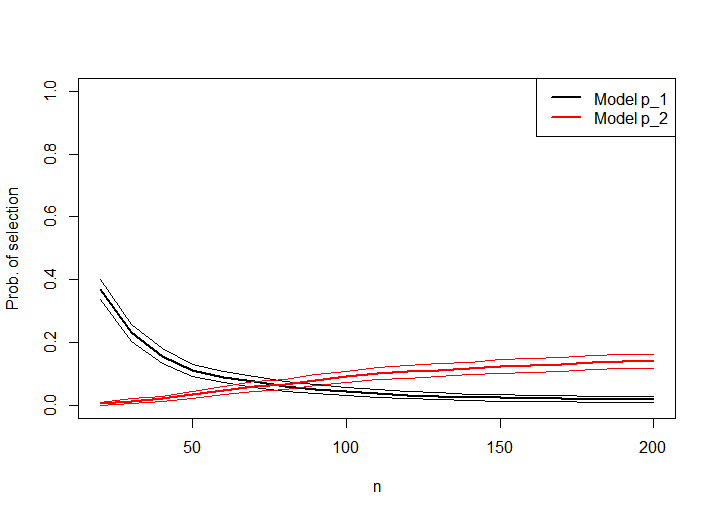}
	}
	\subfigure[Probability of selection by $\mathbb{C}^{({\mathpzc P})}(n,N)$]{%
		\includegraphics[width=0.48\textwidth]{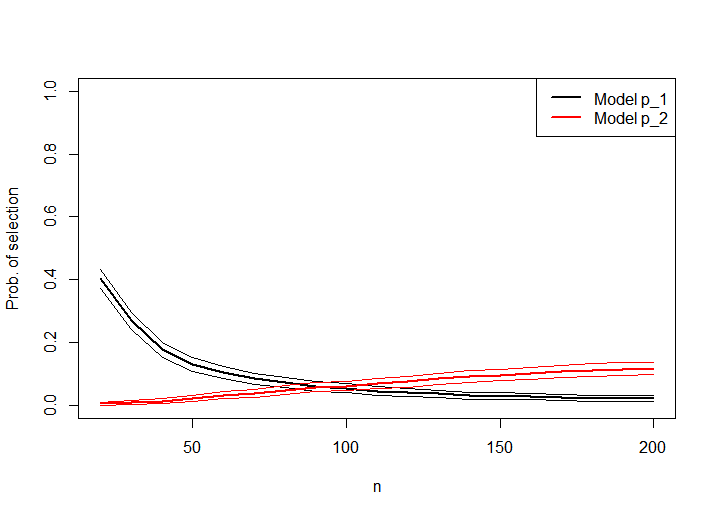}
	}\\
	\caption{\footnotesize Probability of selecting model $p_1$; the thick line is the simulation mean and the thin lines are plus and minus two standard errors. } \label{fig:Prob_Cp}
\end{figure}

Figure \ref{fig:Prob_Cp} depicts simulation estimates of the probability of selecting models ${\mathpzc P}_1$ and ${\mathpzc P}_2$ as a function of $n$, using ${C}^{({\mathpzc P})}(n,{\bf N})$ and the jackknifed  $\mathbb{C}^{({\mathpzc P})}(n,{\bf N})$, where all $2^5$ possible sub-models ${\mathpzc P}$ are considered; for clarity we present the curves of ${\mathpzc P}_1$ and ${\mathpzc P}_2$ only.  For each $n$ and for each simulated dataset, ${C}^{({\mathpzc P})}(n,{\bf N})$ and  $\mathbb{C}^{({\mathpzc P})}(n,{\bf N})$ are calculated for all ${\mathpzc P}$.  
The empirical averages over 100 simulations of selecting  models ${\mathpzc P}_1$ and ${\mathpzc P}_2$ out of the $2^5$ sub-models for each $n$ are plotted in Figure \ref{fig:Prob_Cp}.
The bias correction increases the probability of selecting model ${\mathpzc P}_1$ for small $n$. This improves the selection for small or moderate values of $n$.  
For the problem of selecting a common model for $\cal J$ datasets, the bias correction becomes more significant, as demonstrated next.

\subsection{Multiple datasets}\label{sec:simJ}

We now consider the case of ${\cal J}>1$ datasets. Suppose that $G_\theta$ is given by model \eqref{eq:xy} with $b_{\theta,0}=1$, $b_{\theta,k}= W_k(b_k+Z_k)$ for $k=1,\ldots,5$, where $b_k$ is given in \eqref{eq:bs}, $(Z_1,\ldots,Z_5) \sim N(0,0.2^2)$, $W_k$ is $\pm1$ with equal probability, and all the above random variables are independent (thus determining the distribution ${\mathscr P}$ of Section \ref{sec:supersuper}), and then fixed throughout this section.
The expected $b_{\theta,k}^2$ is approximately equal to $b_k^2$ in \eqref{eq:bs}, but about half of the $b_{\theta,k}$'s are positive and half are negative.  
The number of regression datasets is ${\cal J}=100$, and
$N_j=20, 100$, and 200  for $1 \le j \le 33$, $34 \le j \le 66$, $67 \le j \le 100$, respectively. 

In the case of observing   all regressions  (see Section \ref{subsec:mallowssev}, Equation \eqref{eq:rtild}), we wish to estimate ${\bf R}(n,{\mathpzc P})$, whereas in the case of observing a sample of regressions from the distribution ${\mathscr P}$ (see Section \ref{sec:supersuper}, Equation \eqref{eq:predrrr}), the relevant quantity is ${\bf R}_{pop}(n,{\mathpzc P})$. Computing the latter quantity is difficult, and instead we use the approximation  ${\bf R}(n,{\mathpzc P})$, which is justified by the law of large numbers and the central limit theorem (see Lemma \ref{lem:tild_AR}). Thus we now focus on estimating ${\bf R}(n,{\mathpzc P})$ and selecting according to its estimate.
The plot of ${\bf R}(n,{\mathpzc P})$ for $\cal J$=100 and ${\mathpzc P}={\mathpzc P}_1, {\mathpzc P}_2$ is similar to Figure \ref{fig:pred} and therefore is not presented here. 

\begin{figure}[ht!]
	\subfigure[Boxplots of ${\bf C}^{({\mathpzc P}_1)}(n,\bf{N}) - {\bf C}^{({\mathpzc P}_2)}(n,\bf{N})$]{%
		\includegraphics[width=0.48\textwidth,height=0.3\textheight]{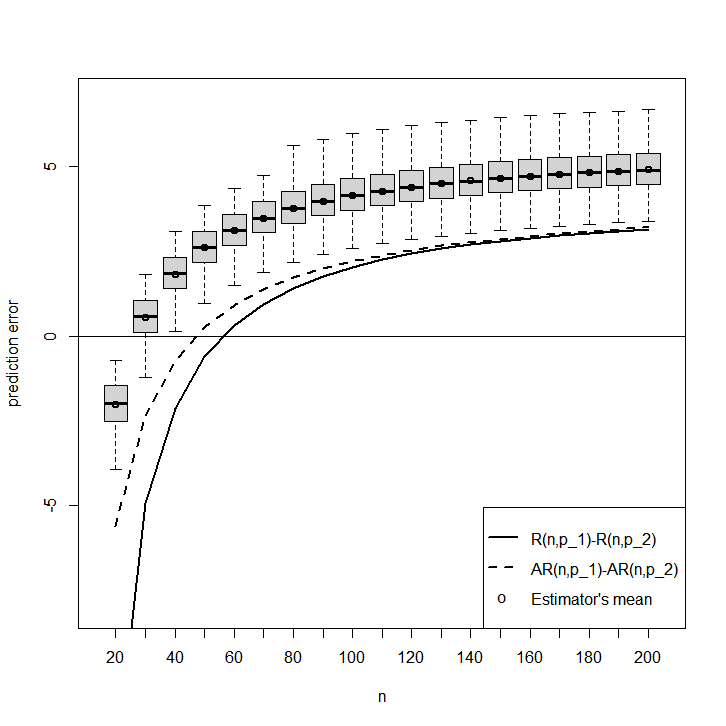}
	}
	\subfigure[Boxplots of  $\pmb{\mathbb{C}}^{({\mathpzc P}_1)}(n,\bf{N}) - \pmb{\mathbb{C}}^{({\mathpzc P}_2)}(n,\bf{N})$]{%
		\includegraphics[width=0.48\textwidth, height=0.3\textheight]{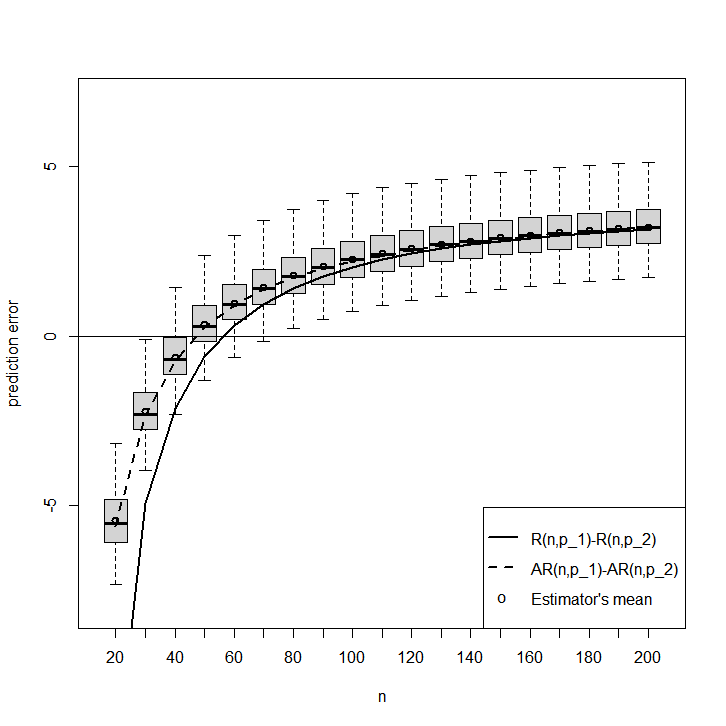}
	}\\
	
	\caption{\footnotesize Same plots as in Figure \ref{fig:Cp} when there are ${\cal J}=100$ samples.} \label{fig:Cp_multi}
\end{figure}

\begin{figure}[ht!]
	\subfigure[Probability of selection by ${\bf C}^{({\mathpzc P})}(n,\bf{N})$]{%
		\includegraphics[width=0.48\textwidth, height=0.28\textheight]{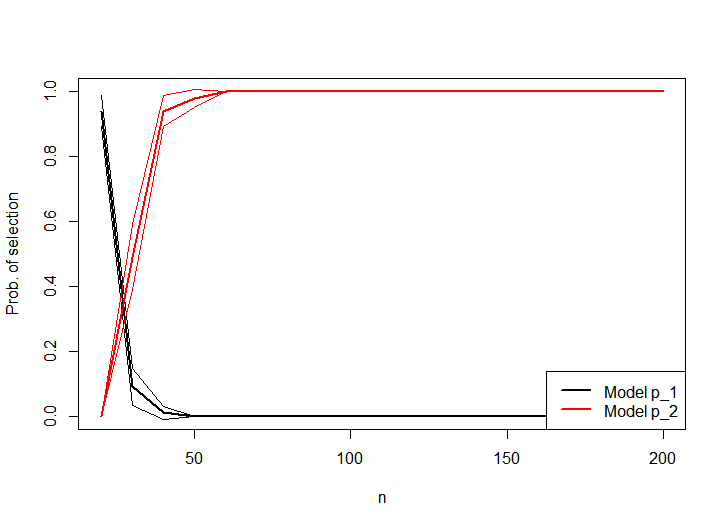}
	}
	\subfigure[Probability of selection by $\pmb{\mathbb{C}}^{({\mathpzc P})}(n,\bf{N})$]{%
		\includegraphics[width=0.48\textwidth,height=0.28\textheight]{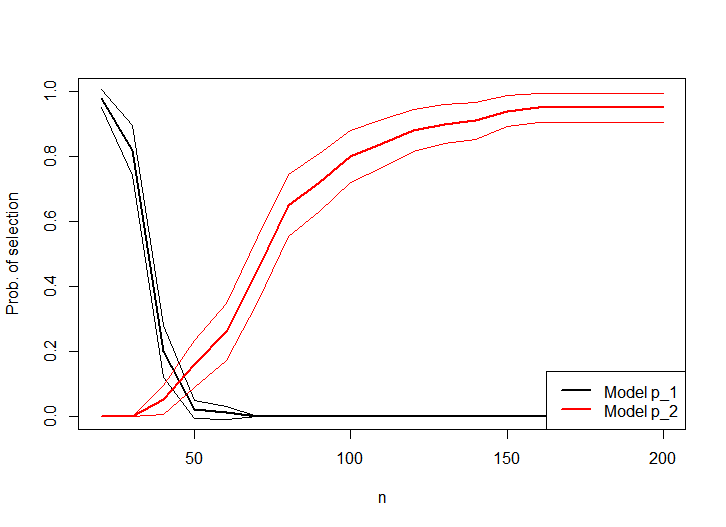}
	}\\
	\caption{\footnotesize Same plots as in Figure  \ref{fig:Prob_Cp} when there are ${\cal J}=100$ samples.} \label{fig:Prob_CpJ}
\end{figure}

Figure \ref{fig:Cp_multi} parallels Figure \ref{fig:Cp}, where now in the case of $\cal J$ datasets,  ${\bf C}^{({\mathpzc P})}(n,{\bf N})$ and $\pmb{\mathbb{C}}^{({\mathpzc P})}(n,{\bf N})$ replace $C^{({\mathpzc P})}(n,{N})$, and   $\mathbb{C}^{({\mathpzc P})}(n,{N})$, respectively; see \eqref{eq:JCPP} and \eqref{eq:CpJackJ}; the number of simulations to evaluate ${\bf R}_{pop}(n,{\mathpzc P})$ and to produce Figures \ref{fig:Cp_multi} and \ref{fig:Prob_CpJ} is 100. We see that the jackknife bias correction works well.
Here the variances of the estimates are much smaller, indicating that several datasets can lead to  better estimates and model selection procedures, as predicted by theory. The $y$-axis scale varies between Figures \ref{fig:Cp_multi} and \ref{fig:Cp}, in a way that undermines their difference.

Figure \ref{fig:Prob_CpJ} plots the selection probabilities as a function of $n$ (out of all $2^5$ sub-models). Unlike the case ${\cal J}=1$ (see Figure \ref{fig:Prob_Cp}), model ${\mathpzc P}_1$ (respectively, model ${\mathpzc P}_2$) is selected with high probability for small $n$ (respectively, large $n$). Recall that it is optimal to select model ${\mathpzc P}_1$ (respectively, ${\mathpzc P}_2$) when $n \le 50$, (respectively, $n\ge 50$). Selecting according to ${\bf C}^{({\mathpzc P})}(n,{\bf N})$ leads to favoring ${\mathpzc P}_2$ (or other models) when $n$ is greater than approximately 25 (instead of 50) and  $\pmb{\mathbb{C}}^{({\mathpzc P})}(n,{\bf N})$ corrects this bias. 
Thus, the probability of correct  model selection is much higher when using the ${\cal J}=100$ datasets (see Figure \ref{fig:pred}). Clearly, the probability of making a correct selection depends on the number of datasets $\cal J$, the similarity among  the $\cal J$ models, the noise level in the models, and the sample size $n$. \\

\vspace{-0.7cm}

\subsection{Comparisons to other approaches}\label{sec:compare}

We considered the possibility of   concatenating the whole training sample and performing a single regression with an intercept for each $j$. In this simulation, since about half of the $b_{\theta,k}$'s are positive and half are negative, the resulting regression model leads to a higher prediction error than the one of $\pmb{\mathbb{C}}^{({\mathpzc P})}(n,{\bf N})$. The latter has estimated prediction error of 56.1 (SE=0.03) (see Table \ref{tab:pred} below), while for the ordinary least squares applied to the concatenated dataset the corresponding number is 63.6 (SE=0.1), computed by averaging the prediction error over 1000 independent datasets with the same distribution.   For ridge and lasso estimators applied to the concatenated dataset (calculated using the glmnet package, where the tuning parameter was computed using cross-validation), the prediction error was slightly higher: 64.1 (SE=0.1) and 63.9 (SE=0.1) for ridge and lasso respectively.

Another approach is to consider a separate model selection algorithm for each of the ${\cal J}$ datasets. We considered three selection criteria: $\mathbb{C}^{({\mathpzc P})}(n,N_j)$ with $n=N_j$ as in \eqref{eq:cpp} (applied to each dataset separately), Mallows' $C_p$ and BIC. 
The means of the resulting prediction errors are given in the Table \ref{tab:pred} below as well as that of $\pmb{\mathbb{C}}^{({\mathpzc P})}(n,N)$ (where the same model is selected for all $j$'s with the same sample size $N_j$). The datasets are divided into three categories according to their sample sizes and the mean is reported for each category separately. Recall that $N_j=20, 100$, and 200  for $1 \le j \le 33$, $34 \le j \le 66$, and $67 \le j \le 100$, respectively. The prediction error $R_j(N_j,{\mathpzc P}_\ell)$  was evaluated as in Figure \ref{fig:pred}.  
Table \ref{tab:pred} shows that $\pmb{\mathbb{C}}^{({\mathpzc P})}(n,N)$ leads to smaller prediction errors and the improvements is higher for smaller sample-sizes, where  borrowing power from other datasets is more important. 

\begin{table} [ht!]
	\caption{The means of the prediction errors $R_j(N_j,{\mathpzc P}_\ell)$, where ${\mathpzc P}_\ell$ is selected by different methods. The standard errors are about 0.03.}
	\begin{center}
		{\small
			\begin{tabular}{|c||c c c c|}
				\hline
				{}&\multicolumn{4}{|c|}{Model selection method}\\
				\hline
				$N_j$ & $\pmb{\mathbb{C}}^{({\mathpzc P})}(n,N)$ &  $\mathbb{C}^{({\mathpzc P})}(n,N_j)$ & Mallows' $C_p$ & BIC\\
				\hline 
				\hline
				$20$ &  61.4 &67.4& 66.4 &66.1\\
				$100$ &  54.4 & 55.2 & 55.2 & 56.1\\
				$200$ & 52.7 & 53.4 & 53.4 & 54.5\\
				\hline
				Average & 56.1 & 58.7 & 58.3 & 58.9 \\
				\hline
		\end{tabular}}\label{tab:pred}
	\end{center}
\end{table}

\section{Prediction of durations of medical examinations}\label{sec:med}
In this section we analyze a dataset of outpatients' hospital visits. Different models are considered in order to predict the actual appointments' durations as opposed to the planned durations. 

\subsection{Description of the data}
The dataset analyzed is taken from the SEE Lab at the Technion. It consists of information on 140,924 hospital visits that took place in a certain US hospital for about two years between 2013 and 2015. For each visit, both the planned time and the actual time are reported. The goal was to provide a more accurate prediction of the actual duration than the planned one. In this dataset there is information on 44,516  patients and 258 doctors, out of whom  34 doctors had fewer than 50 visits. We shall focus on the rest, which corresponds to 99.5\% of all visits. The regression coefficients will differ between doctors, and the goal is to select one common subset of covariates (for each $n$) for all doctors for prediction of visit durations.

\begin{figure}[ht!]
	\begin{center}
		\includegraphics[width=0.65\textwidth,height=0.2\textheight]{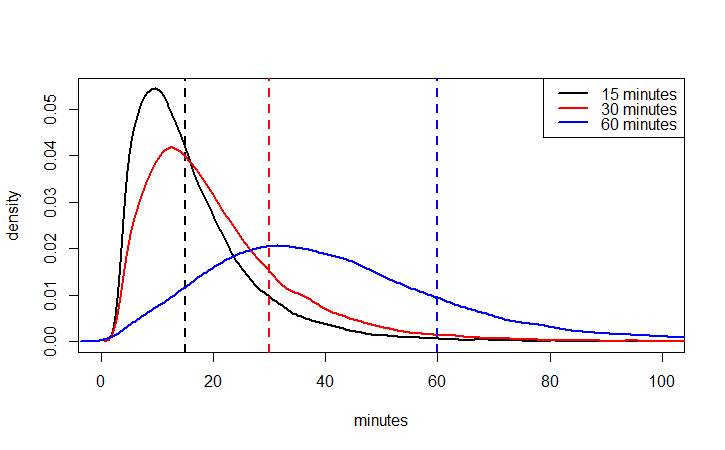}
	\end{center} 
	\caption{\footnotesize Estimated density of the actual duration for the time slots of 15, 30, and 60 minutes. The vertical dashed lines are at 15, 30, and 60 minutes. \label{fig:dens}}
\end{figure}

The distribution of the planned duration is given in Table \ref{tab:1} and Figure \ref{fig:dens} plots the estimated density (a normal kernel estimate using the R command ``density") of the actual durations for the time slots of 15, 30, and 60 minutes. Actual durations are obtained by a real-time location system (RTLS). The means are 16.7, 21.3, and 41.2, respectively.

\begin{table} [ht!]
	\caption{The distribution of the planned duration.}
	\begin{center}
		{\small
			\begin{tabular}{c||c c c c c c}
				%\hline
				minutes & 15 & 30 & 45 & 60 & other\\
				\hline 
				\hline
				percentage & 29.8\% & 52.6\% & 1.8\% & 15.5\% & 0.3\%
		\end{tabular}}\label{tab:1}
	\end{center}
\end{table}

\subsection{A  regression model}\label{sec:global}

The original dataset contains a large number of covariates, of which many did not seem to have any predictive power relative to visit durations. For simplicity of presentation, we focus on a small number of covariates  that seem most relevant.
We aim to predict actual duration, using the following covariates:
\begin{itemize}
	\item duration\_planned = the planned duration of the visit in minutes.
	\item duration\_planned\_2 = the planned duration  in minutes of the visit, squared.
	\item last = the planned  minus the actual  duration of the previous visit of the same patient (taken to be  0 for the first visit of the patient).
	%\item half = whether the exam is planned to start on the hour. It turns out that these kind of exams tend to be longer.
	\item hour\_end = whether the exam is planned to end on the hour. It turns out that these kinds of visits tend to be somewhat longer. 
	\item type = there are two types of examinations: consultation/examination only, or the above plus treatment. In either case, only the first part counts as duration. 
\end{itemize}
%The output of the ``lm'' command in R applied to the whole dataset (ignoring the doctors' index)  is given in Figure \ref{fig:R}. 
Standard statistical inference of the linear regression model of the whole dataset (ignoring the doctors' index) reveals that all of the above covariates besides ``type'' are statistically significant; however, the standard error of the residuals is 15.33, and $R^2=0.227$, suggesting that the prediction error is quite large.

\subsection{$\pmb{\mathbb{C}}^{({\mathpzc P})}$ and model selection}\label{sec:ms}

In our notation, each doctor is indexed by $j$, and $N_j$ is the number of visits to doctor $j$ in the dataset; $N_j$ varies between 50 and 2135.
We demonstrate our approach by focusing on four candidate models that have the smallest (or nearly smallest) $\pmb{\mathbb{C}}^{({\mathpzc P})}$ from all submodels of the five covariates (all models included the intercept term) for relevant sample sizes $n$. These models are ${\mathpzc P}_1$ -- the model with the covariates:  duration\_planned,  duration\_planned\_2;  ${\mathpzc P}_2$ -- the model with the same covariates as in ${\mathpzc P}_1$ and additionally, the variable ``last''; ${\mathpzc P}_3$ -- the model with the same covariates of ${\mathpzc P}_2$ and additionally, the variable ``type''; and ${\mathpzc P}_4$ -- the full model.  For certain submodels estimation is possible only for a subset of the doctors since ${\mathbb{X}}_{j,N_j}^{({\mathpzc P})'}{\mathbb{X}}_{j,N_j}^{({\mathpzc P})}$ is not always invertible. Therefore $\cal J$ varies between the models. For the models ${\mathpzc P}_1$ and ${\mathpzc P}_2$, invertibility held for 96 
doctors and for the models ${\mathpzc P}_3$ and ${\mathpzc P}_4$, the corresponding number is 95, and so for these models ${\cal J} =96$ or ${\cal J} =95$. 
In this case, $\pmb{\mathbb{C}}^{({\mathpzc P})}(n,{\bf N})$ is based only on this subset.

\begin{figure}[ht!]
	\begin{center}
		\includegraphics[width=.5\textwidth]{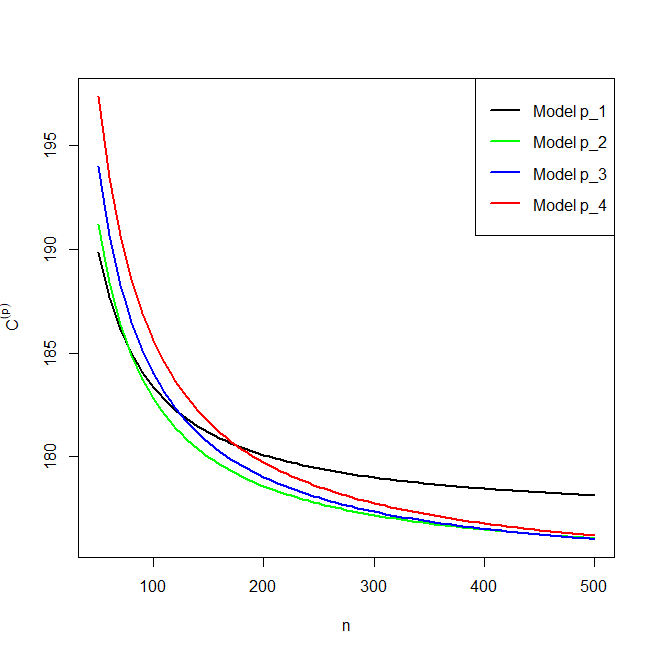}
	\end{center} 
	\caption{\footnotesize A plot of $\pmb{\mathbb{C}}^{(p)}(n,{\bf N})$ for $p=p_1,p_2,p_3,p_4$ and $n=50,55,\ldots,500$. \label{fig:C_p_1_p_2}}
\end{figure}

Figure \ref{fig:C_p_1_p_2} plots $\pmb{\mathbb{C}}^{({\mathpzc P})}(n,{\bf N})$ for ${\mathpzc P}={\mathpzc P}_\ell$ where $\ell=1,\ldots,4$ and $n$ is between 50 and 500. For $n$  smaller than approximately 80, model ${\mathpzc P}_1$ is the best among the candidate models; for $n$ between 80 and 450, ${\mathpzc P}_2$ has a smaller $\pmb{\mathbb{C}}^{({\mathpzc P})}$, and for larger $n$, ${\mathpzc P}_3$ is the best, but ${\mathpzc P}_2$ is very close. In terms of GENO, we have, for example, that for $n=50$, $\widehat{\rm GENO}(n,{\mathpzc P}_1,{\mathpzc Q})$ for ${\mathpzc Q}={\mathpzc P}_2, {\mathpzc P}_3,$ and  ${\mathpzc P}_4$ equals 54, 63, 73, respectively. The latter number means that model ${\mathpzc P}_4$ (the full model) would require 73 observations to achieve the same prediction error as model ${\mathpzc P}_1$ with $n=50$ observations. Also,
$\widehat{\rm GENO}(200,{\mathpzc P}_2,{\mathpzc P}_1)=370$; if one considers using only the planned duration (${\mathpzc P}_1$) or using model ${\mathpzc P}_2$, that is, adding the variable ``last" with the information on the last visit, which may not be available for some patients, then the estimated prediction error by the model ${\mathpzc P}_2$ with $n=200$ observations can be achieved without knowing  ``last" by ${\mathpzc P}_1$, with $n=370$. It is then left to the user to decide whether to invest in measuring ``last" or in using a larger sample, if such a sample is available.

Table \ref{tab:2} reports $\pmb{\mathbb{C}}^{({\mathpzc P}_\ell)}(n,{\bf N})$ for different sample sizes $n$. Standard deviations estimated by the bootstrap, and cross-validation estimates of $\textbf{R}(n,{\mathpzc P})$, are also provided. The latter estimates are computed  only for $j$'s where $N_j > n$. For each such $j$, the data were split at random into a training set with $n$ observations, and a testing set of size $N_j-n$. The estimates $\hat{\bs \beta}^{({\mathpzc P})}_{j,n}$ are based on the training set and the prediction error ${R}_j(n,{\mathpzc P})$ is evaluated using the testing set. This procedure was repeated 1,000 times and the average prediction error is reported. The cross-validation estimates are mostly within one standard error of the  $\pmb{\mathbb{C}}^{({\mathpzc P})}$ values, and the two approaches lead to selection of the same models.

Table \ref{tab:3} reports the values of $\pmb{\mathbb{C}}^{({\mathpzc P})}(n,{\bf N})-\pmb{\mathbb{C}}^{({\mathpzc Q})}(n,{\bf N})$ together with a bootstrap estimate of the standard deviation for different values of $n$ and various pairs of models. Also the differences of the corresponding cross-validation estimates are given. The standard deviations of Table \ref{tab:3} are much smaller than those  of Table \ref{tab:2}. This is consistent with our theoretical results that comparison of two similar models leads to a small  estimation error (see the discussion after Theorem \ref{thm:o1n}). The table shows which pairs ${\mathpzc P}$, ${\mathpzc Q}$ differ significantly, and for which values of $n$.

 \begin{table} [ht!]
	\caption{ \footnotesize{$\pmb{\mathbb{C}}^{({\mathpzc P})}(n,{\bf N})$ for ${\mathpzc P}={\mathpzc P}_1,{\mathpzc P}_2,{\mathpzc P}_3,{\mathpzc P}_4$ and for $n=50,150,500$. Bootstrap standard deviations (SD) and cross-validation (CV) estimates are also provided.}}
	\begin{center}
		{
			\scriptsize\resizebox{\columnwidth}{!}{\begin{tabular}{|c||c|c||c|c||c|c||c|c|}
				\hline
				{}&\multicolumn{2}{|c||}{Model ${\mathpzc P}_1$} & \multicolumn{2}{|c||}{Model ${\mathpzc P}_2$}& \multicolumn{2}{|c||}{Model ${\mathpzc P}_3$}&\multicolumn{2}{|c|}{Model ${\mathpzc P}_4$}\\
				\hline
				n & $\pmb{\mathbb{C}}^{({\mathpzc P}_1)}$ (SD) &  CV & $\pmb{\mathbb{C}}^{({\mathpzc P}_2)}$ (SD) &  CV & $\pmb{\mathbb{C}}^{({\mathpzc P}_3)}$ (SD) &  CV & $\pmb{\mathbb{C}}^{({\mathpzc P}_4)}$ (SD) &  CV\\
				\hline 
				\hline
				50 & 189.9 (3.4) & 183.9 & 191.2 (3.5) & 185.9 & 194.0 (3.5) & 190.1 & 197.4 (3.8) & 194.5\\
				150 & 181.1 (3.3) & 181.1 & 180.0 (3.3) & 180.4 & 180.7 (3.4) & 181.3 & 181.7 (3.4) & 183.0\\
				500 & 178.1 (3.3) & 183.0 & 176.1 (3.1) & 181.3  & 176.0 (3.4) & 181.2 & 176.2 (3.3) & 181.3 \\
				\hline
		\end{tabular}}}\label{tab:2}
	\end{center}
\end{table}

\begin{table} [ht!]
	\caption{ \footnotesize{The values of $\pmb{\mathbb{C}}^{({\mathpzc P})}(n,{\bf N})-\pmb{\mathbb{C}}^{({\mathpzc Q})}(n,{\bf N})$ for $n=50,150,500$. Bootstrap standard deviations (SD) and cross-validation (CV) estimates are also provided. Boldface numbers indicate differences that are significantly (more than two SD's) non-zero.}}
	\begin{center}
		{
		\scriptsize{	\begin{tabular}{|c|c||c|c||c|c||c|c|}
				\hline
				\multicolumn{2}{|c||}{}&\multicolumn{2}{|c||}{$n=50$} & \multicolumn{2}{|c||}{$n=150$}& \multicolumn{2}{|c|}{$n=500$}\\
				\hline
				${\mathpzc P}$ & ${\mathpzc Q}$ & $\pmb{\mathbb{C}}^{({\mathpzc P})}-\pmb{\mathbb{C}}^{({\mathpzc Q})}$ (SD) &  CV & $\pmb{\mathbb{C}}^{({\mathpzc P})}-\pmb{\mathbb{C}}^{({\mathpzc Q})}$ (SD) &  CV & $\pmb{\mathbb{C}}^{({\mathpzc P})}-\pmb{\mathbb{C}}^{({\mathpzc Q})}$ (SD)  &  CV \\
				\hline 
				\hline
				${\mathpzc P}_1$ & ${\mathpzc P}_2$ &\quad {\bf -1.3}\quad \quad (0.3) & -2.0 & \quad  {\bf 1.2} \quad \quad (0.3)  & 0.7 & \quad {\bf 2.1} \quad \quad (0.3)  &1.7\\
				${\mathpzc P}_1$ & ${\mathpzc P}_3$ &\quad  {\bf -4.1} \quad \quad (0.4)  & -6.2 & \quad  0.5 \quad \quad (0.4)  & -0.2 & \quad{\bf 2.1} \quad \quad (0.4)  & 1.8\\
				${\mathpzc P}_1$ & ${\mathpzc P}_4$ &\quad  {\bf -7.5} \quad \quad (0.5)  & -10.6 & \quad  -0.5 \quad \quad (0.5)  & -1.9 & \quad {\bf 1.9} \quad \quad (0.5) & 1.7\\
				${\mathpzc P}_2$ & ${\mathpzc P}_3$ &\quad  {\bf -2.7} \quad \quad (0.2) & -4.2 & \quad  {\bf -0.7 } \quad \quad (0.2) & -0.9 & \quad 0.0 \quad \quad (0.2) & 0.1 \\
				${\mathpzc P}_2$ & ${\mathpzc P}_4$ &\quad  {\bf -6.2} \quad \quad (0.4) & -8.6 & \quad {\bf -1.7}  \quad \quad (0.3)  & -2.6 & \quad -0.1  \quad \quad (0.3) & 0.0 \\ 
				${\mathpzc P}_3$ & ${\mathpzc P}_4$ &\quad  {\bf -3.4}  \quad \quad (0.3)& -4.4 & \quad  {\bf -1.0}   \quad \quad (0.2) & -1.7 & \quad {-0.2}  \quad \quad (0.2) & -0.1\\
				\hline		
		\end{tabular}}}\label{tab:3}
	\end{center}
\end{table}

\subsection{Comparisons to other approaches}\label{sec:comparison}

As in Section \ref{sec:compare} we compare our method to other approaches. One possibility is to concatenate the whole training sample and add a categorical variable for the doctors.
The (10-fold) cross-validation estimate of the prediction error of the OLS is 204.5. The corresponding numbers for the ridge and lasso estimates (applied to the concatenated data) are similar: 205.0 and 204.8. The estimates of prediction errors of our method are smaller: they vary between 190 and 176 for $50 \le n \le 500$ (See Figure \ref{fig:C_p_1_p_2} and Table \ref{tab:2}).

\begin{figure}[ht!]
	\begin{center}
		\includegraphics[width=.8\textwidth]{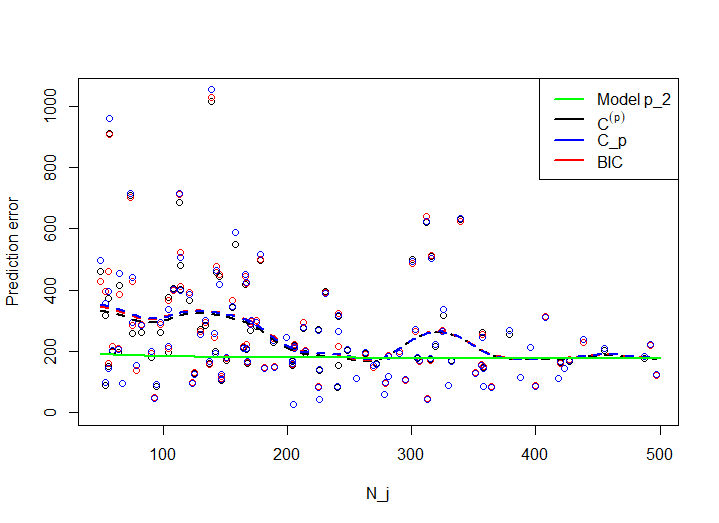}
	\end{center} 
	\caption{\footnotesize Prediction errors (estimated by cross-validation) of the model selection methods $\mathbb{C}^{(p)}(N_j,N_j)$, Mallows' $C_p$ and BIC applied to each dataset separately compared to the common model ${\mathpzc P}_2$. The dashed lines are Gaussian-kernel smoothing. The green line is $\pmb{\mathbb{C}}^{({\mathpzc P}_2)}(n,\bf{N})$ as a function of $n \in [50,500]$. \label{fig:sep}}
\end{figure}

A different approach is to preform a separate model selection for each of the ${\cal J}$ datasets. As in Section \ref{sec:compare}, three selection criteria are considered, $\mathbb{C}^{({\mathpzc P})}(N_j,N_j)$ (the bias-corrected $C^{({\mathpzc P})}(n,N)$ with $n=N_j$), Mallows' $C_p$ and BIC. Figure \ref{fig:sep} plots the cross-validation estimates of the prediction errors of the selected models by the three criteria as a function of the sample size $n=N_j$. A normal-kernel smoothing is drawn to illustrate the average prediction error as a function of the sample size. Also plotted is the prediction error of the common model selection ${\mathpzc P}_2$ (which is close to optimal for sample sizes between 50 and 500) as estimated by $\pmb{\mathbb{C}}^{({\mathpzc P}_2)}$. Table \ref{tab:2} shows that the latter estimate is rather close to its cross validation estimate. Figure \ref{fig:sep} shows that: a. the differences between the three  selection criteria are small; b. a common model selection by $\pmb{\mathbb{C}}^{({\mathpzc P})}$ is better on average than a separate model selection; c. the latter statement is especially true for small sample sizes where borrowing strength is more important.

\noindent\textbf{{Acknowledgment}:}
We are grateful to Avishai Mandelbaum for providing access to the SEE Lab dataset we analyzed, and to Ella Nadjharov for creating the files we needed. We also wish to thank the associate editor and the two referees for their very useful comments.

	\section{Appendix A: Proofs}\label{sec:AppA}
	
	Recall that we use ${\mathpzc P}$ to denote a subset of the covariates, to which we sometimes refer  as a model, and denote its size with the corresponding letter $p$. We suppress it in most of our notation below and instead of ${\bf X}^{(\mathpzc P)}$ we write ${\bf X}$ and assume it is in $\mathbb{R}^p$.

	\noindent{\bf Proof of Theorem \ref{thm:o1n}}. We first prove \eqref{eq:approx_pred3}.  
	For $\widehat{{\bs \beta}}_n$ computed from a sample $D=\{({\bf X}_i, Y_i):i=1,\ldots,n\}$, and a pair of new observations from the same distribution $({\bf X}, Y)$, independent of $D$, we have
	\begin{multline}\label{eq:predX}
	E(Y-{\bf X}^{'}\widehat{{\bs \beta}}_n)^2=E(Y-{\bf X}^{'}{{\bs \beta}})^2+E[{\bf X}^{'}(\widehat{{\bs \beta}}_n-{\bs \beta})]^2-2E[(Y-{\bf X}^{'}{{\bs \beta}}){\bf X}^{'}(\widehat{{\bs \beta}}_n-{\bs \beta})]\\
	=E(Y-{\bf X}^{'}{{\bs \beta}})^2+E[{\bf X}^{'}(\widehat{{\bs \beta}}_n-{\bs \beta})]^2,
	\end{multline}
	where the last term in the first line of \eqref{eq:predX} vanishes since $E[(Y-{\bf X}^{'}{{\bs \beta}}){\bf X}]=E[e{\bf X}]={\bf 0}$ and $Y$ and ${\bf X}$ are independent of $\widehat{{\bs \beta}}_n$. This argument holds also if we condition on $H_n$ (see \eqref{eq:setGn}, and Theorem \ref{thm:BeGood}).
	By \eqref{eq:predX} we have that
	\begin{equation}\label{eq:R_AR}
	n[R(n,{\mathpzc P})-AR(n,{\mathpzc P})]=E[{\bf X}^{'}\sqrt{n}(\widehat{{\bs \beta}}_n-{\bs \beta})]^2-tr({\mathbb{V}}).
	\end{equation}
	Using independence of ${\bf X}$ and $\widehat{{\bs \beta}}_n$ again we have
	\[
	E[{\bf X}^{'}\sqrt{n}(\widehat{{\bs \beta}}_n-{\bs \beta})]^2=tr \{E[{\bf X} {\bf X}']E [n (\widehat{{\bs \beta}}_n-{\bs \beta}) (\widehat{{\bs \beta}}_n-{\bs \beta})']\}=tr\{\mathbb{Q}  E [n (\widehat{{\bs \beta}}_n-{\bs \beta}) (\widehat{{\bs \beta}}_n-{\bs \beta})']\}. 
	\]
	For the proof of Theorem \ref{thm:BeGood} the expectations should be conditioned on the set $H_n$, whose probability is large, and the conditioning does not affect the rates we obtain.
	
	By Equation (7.3) of \citet{Hansen},
	\begin{equation*}%\label{eq:beta-beta}
	\sqrt{n}(\widehat{{\bs \beta}}_{n}-{\bs \beta})=\widehat{\mathbb{Q}}_n^{-1}\frac{1}{\sqrt{n}}\sum_{i=1}^n {\bf X}_i e_i=\widehat{\mathbb{Q}}_n^{-1}{\bf U}_n.
	\end{equation*}
	Since $E\{tr({\bf U}_n{\bf U}_n'\mathbb{Q}^{-1})\}=tr({\mathbb{V}})$ we rewrite the right-hand side of \eqref{eq:R_AR} as 
	\begin{multline}\label{eq:tr_Q}
	E\{tr(\mathbb{Q} \widehat{\mathbb{Q}}_n^{-1}{\bf U}_n {\bf U}_n' \widehat{\mathbb{Q}}_n^{-1} ) - tr (\mathbb{Q} \mathbb{Q}^{-1} {\bf U}_n {\bf U}_n' \mathbb{Q}^{-1} )\}\\=tr(\mathbb{Q} E\{ \widehat{\mathbb{Q}}_n^{-1}{\bf U}_n {\bf U}_n' \widehat{\mathbb{Q}}_n^{-1}- \mathbb{Q}^{-1} {\bf U}_n {\bf U}_n' \mathbb{Q}^{-1}\}).
	\end{multline}
	In order to prove \eqref{eq:approx_pred3} we now show that the latter expectation is of order $O(1/n^{1/2})$. To this end, notice that
	\begin{align*}
	&\widehat{\mathbb{Q}}_n^{-1}{\bf U}_n {\bf U}_n' \widehat{\mathbb{Q}}_n^{-1}- \mathbb{Q}^{-1} {\bf U}_n {\bf U}_n' \mathbb{Q}^{-1}\\&=(\widehat{\mathbb{Q}}_n^{-1} - \mathbb{Q}^{-1}){\bf U}_n {\bf U}_n' \widehat{\mathbb{Q}}_n^{-1}+
	\mathbb{Q}^{-1} {\bf U}_n {\bf U}_n' (\widehat{\mathbb{Q}}_n^{-1}-\mathbb{Q}^{-1}) .
	\end{align*}
	
	We now deal with the first term on the right-hand side above, the other term being similar, and simpler. We have
	\begin{equation}\label{eq:Q_hat_minus1}
	\widehat{\mathbb{Q}}_n^{-1}-\mathbb{Q}^{-1}=\mathbb{Q}^{-1}(\mathbb{Q}-\widehat{\mathbb{Q}}_n)\widehat{\mathbb{Q}}_n^{-1}
	\end{equation}
	{and therefore (recall \eqref{eq:tr_Q}) we consider the expectation of
		\[
		tr\Big(\mathbb{Q}(\widehat{\mathbb{Q}}_n^{-1} - \mathbb{Q}^{-1}){\bf U}_n {\bf U}_n' \widehat{\mathbb{Q}}_n^{-1}\Big)=tr\Big( (\mathbb{Q}-\widehat{\mathbb{Q}}_n)\widehat{\mathbb{Q}}_n^{-1}{\bf U}_n {\bf U}_n' \widehat{\mathbb{Q}}_n^{-1}\Big).
		\] 
		This matrix is a product of random matrices of the form $ABCD$ where $A=\mathbb{Q}-\widehat{\mathbb{Q}}_n$, $B=D=\widehat{\mathbb{Q}}_n^{-1}$, and $C={\bf U}_n {\bf U}_n'$. 
		The trace is a sum of products of entries from all the matrices appearing in the 
		product.  Different choices of powers can be made, but 
		we use H\"{o}lder's inequality in the form $E|abcd| \le (Ea^{12})^{1/12} (E|b|^3)^{1/3} (Ec^{4})^{1/4} (E|d|^3)^{1/3}$ for simplicity. Here $a$ is an entry from $A$, $b$ an entry from from $B$, etc., and the triangle inequality can then be used to bound the sum comprising the trace.}
	
	 For each element $j,k$ of $\mathbb{Q}-\widehat{\mathbb{Q}}_n$ we have 
		\[
		E(\mathbb{Q}-\widehat{\mathbb{Q}}_n)_{j,k}^{12} =  E \Big( \frac{1}{n}\sum_{i=1}^n [E(X_j X_k)-X_{ij} X_{ik}]\Big)^{12}.
		\]
		The summands $E(X_j X_k)-X_{ij} X_{ik}$ have zero expectation;  expanding 	$(\mathbb{Q}-\widehat{\mathbb{Q}}_n)_{j,k}^{12}$  we see that the number of nonvanishing terms when the expectation is taken is of order $n^6$, and all these terms are bounded by our assumptions.  Therefore, $[E(\mathbb{Q}-\widehat{\mathbb{Q}}_n)_{j,k}^{12}]^{1/12}$ is of order $1/\sqrt{n}$. (Actually, it is easy to see that 24  bounded moments suffice for this argument, and also for bounding the remaining terms, and 24 can be somewhat reduced by a better but more cumbersome choice of the powers in H\"{o}lder's inequality.) A similar computation for the matrix $C$ shows that the required moments of its entries are bounded. The rest of the terms are bounded by our assumptions.
		Now \eqref{eq:approx_pred3} follows.
	
	The proof required the bounded third power of $B=\widehat{\mathbb{Q}}_n^{-1}$, which means that with a mixture of normals we need $n>p+5$. See the Proof of Lemma \ref{lem:normal}.

	We now show \eqref{eq:final}. The definitions of $AR(n,{\mathpzc P})$ and $C^{({\mathpzc P})}(n,N)$ imply that
	\begin{align}\label{eq:Ar-Cp}
	&AR(n,{\mathpzc P})-C^{({\mathpzc P})}(n,N)\nonumber\\&=E(Y-{\bf X}^{'}{{\bs \beta}})^2-\frac{1}{N} ||{\bf Y}_N-\mathbb{X}_N\widehat{\bs \beta}_N||^2+\frac{1}{n}\left\{{tr}({\mathbb{V}})-{tr}({\widehat{\mathbb{V}}_N})\right\}-\frac{{tr}({\widehat{\mathbb{V}}_N})}{N}.
	\end{align}
	Starting with the second term on the right-hand side of \eqref{eq:Ar-Cp}, we have
	\begin{equation*}%\label{eq:eq11}
	||{\bf Y}_N-\mathbb{X}_N{\bs \beta}||^2=||{\bf Y}_N-\mathbb{X}_N\widehat{\bs \beta}_N||^2+||\mathbb{X}_N(\widehat{\bs \beta}_N-{\bs \beta})||^2
	-2({\bf Y}_N-\mathbb{X}_N\widehat{\bs \beta}_N)'(\mathbb{X}_N(\widehat{\bs \beta}_N-{\bs \beta})),
	\end{equation*} 
	where the last term vanishes since $\mathbb{X}_N'({\bf Y}_N-\mathbb{X}_N\widehat{\bs \beta}_N)=0$. Hence,
	\begin{align}\label{eq:stam1}
	&E(Y-{\bf X}^{'}{{\bs \beta}})^2-\frac{1}{N} ||{\bf Y}_N-\mathbb{X}_N\widehat{\bs \beta}_N||^2\nonumber\\&=E(Y-{\bf X}^{'}{{\bs \beta}})^2-\frac{1}{N}||{\bf Y}_N-\mathbb{X}_N{\bs \beta}||^2+\frac{1}{N}||\mathbb{X}_N(\widehat{\bs \beta}_N-{\bs \beta})||^2.
	\end{align}
	
	Recall the notation  ${\bf U}_N=\mathbb{X}_N' {\bf e}_N/\sqrt{N}$. Since by Equation (7.3) of \citet{Hansen} $N^{-1/2}(\widehat{\bs \beta}_N-{\bs \beta})=N^{-1/2}(\mathbb{X}_N'\mathbb{X}_N)^{-1} \mathbb{X}_N' {\bf e}_N = (\mathbb{X}_N'\mathbb{X}_N)^{-1}{\bf U}_N$,  we have
	\begin{multline*}%\label{eq:X_beta}
	\frac 1N||\mathbb{X}_N(\widehat{\bs \beta}_N-{\bs \beta})||^2= \frac 1N(\widehat{\bs \beta}_N-{\bs \beta})'\mathbb{X}_N'\mathbb{X}_N(\widehat{\bs \beta}_N-{\bs \beta})= {\bf U}_N' (\mathbb{X}_N'\mathbb{X}_N)^{-1}  {\bf U}_N\\
	=\frac 1N tr\{  {\bf U}_N{\bf U}_N' (\mathbb{X}_N'\mathbb{X}_N/N)^{-1}\}
	=\frac 1N tr\{  {\bf U}_N{\bf U}_N' \mathbb{Q}^{-1} \}+o_p(1/N),
	\end{multline*}
	where the last equality holds true since $(\mathbb{X}_N'\mathbb{X}_N/N)^{-1}-\mathbb{Q}^{-1}=o_p(1)$. Also,
	\begin{multline}\label{eq:tracesg}
	\frac 1N ||\mathbb{X}_N(\widehat{\bs \beta}_N-{\bs \beta})||^2=
	\frac 1N\Big[tr({\bf U}_N{\bf U}_N'\mathbb{Q}^{-1})-tr({\mathbb{W}}\mathbb{Q}^{-1})+tr({\mathbb{W}}\mathbb{Q}^{-1})\Big]+o_p(1/N)\\
	=\frac 1N \Big[tr({\bf U}_N{\bf U}_N'\mathbb{Q}^{-1})-tr({\mathbb{V}})+tr(\widehat{\mathbb{V}}_N)\Big]+o_p(1/N),
	\end{multline}
	where for the last equality it suffices that
	$\widehat{\mathbb{Q}}_N$ and $\widehat{\mathbb{W}}_N$ are consistent estimates of $\mathbb{Q}$ and $\mathbb{W}$, and therefore $tr(\widehat{\mathbb{V}}_N)$ is consistent for $tr({\mathbb{V}})=tr({\mathbb{W}}\mathbb{Q}^{-1})$.
	Equations \eqref{eq:Ar-Cp}, \eqref{eq:stam1}, and \eqref{eq:tracesg} imply \eqref{eq:final}.
	
	Next we show that $\sqrt{N}\{tr(\widehat{\mathbb{V}}_N)-tr(\mathbb{V})\}$ is asymptotically normal starting with  the asymptotic normality of  $\sqrt{N}(\widehat{\mathbb{W}}_N - \mathbb{W})$.  We have, 
	\begin{align} \label{eq:Wnorm}
	\widehat{\mathbb{W}}_N &=\frac{1}{N}\sum_{i=1}^N{\bf X}_i{\bf X}_i'\widehat{e}_i^{\,2}=\frac{1}{N}\sum_{i=1}^N{\bf X}_i{\bf X}_i'(Y_i-{\bf X}_i'\widehat{{\bs \beta}}_{N})^2\nonumber\\&=\frac{1}{N}\sum_{i=1}^N{\bf X}_i{\bf X}_i'\big(Y_i-{\bf X}_i'{{\bs \beta}} -\{{\bf X}_i'\widehat{{\bs \beta}}_{N}-{\bf X}_i'{{\bs \beta}}\}\big)^2={\frac{1}{N}\sum_{i=1}^N{\bf X}_i{\bf X}_i'\big(Y_i-{\bf X}_i'{{\bs \beta}}\big)^2}\nonumber\\
	&-{\frac{2}{N}\sum_{i=1}^N{\bf X}_i{\bf X}_i'\big(Y_i-{\bf X}_i'{{\bs \beta}}\big) \left( {\bf X}_i'(\widehat{{\bs \beta}}_{N}-{{\bs \beta}}) \right)}+{\frac{1}{N}\sum_{i=1}^N{\bf X}_i{\bf X}_i'({\bf X}_i'\{\widehat{{\bs \beta}}_{N}-{\bs \beta}\})^2}\nonumber\\&=:A-B+C,
	\end{align}
	respectively.
	Starting with the first term we have
	\[
	\sqrt{N}(A-\mathbb{W})=\frac{1}{\sqrt{N}}\sum_{i=1}^n \left({\bf X}_i{\bf X}_i'e_i^2-E[{\bf X}{\bf X}'e^2]\right),
	\]
	which is asymptotically normal. Now  $B$ is obtained by multiplying the sum $\frac{2}{N}\sum_{i=1}^N{\bf X}_i{\bf X}_i'\big(Y_i-{\bf X}_i'{{\bs \beta}}\big)\otimes {\bf X}_i'$  (which converges to a matrix of constants by the law of large numbers) by  ${I}_p \otimes (\widehat{{\bs \beta}}_{N}-{\bs \beta})$, where ${I_p}$ is the identity matrix of order $p$, and $\otimes$ is Kronecker's product. By Equation (7.3) of \citet{Hansen},
	\begin{equation*}%\label{eq:beta-beta}
	\sqrt{N}(\widehat{{\bs \beta}}_{N}-{\bs \beta})=\widehat{\mathbb{Q}}_N^{-1}\frac{1}{\sqrt{N}}\sum_{i=1}^N {\bf X}_i e_i.
	\end{equation*}
	By the law of large numbers and the fact that ${N}(\widehat{{\bs \beta}}_{N}-{{\bs \beta}})_j(\widehat{{\bs \beta}}_{N}-{{\bs \beta}})_k=O_p(1)$, we have that the term $C$ in \eqref{eq:Wnorm} is $O_p(\frac{1}{N})$. Summing up,
	\begin{multline*}
	\sqrt{N}(\widehat{\mathbb{W}}_N-\mathbb{W})\\=\frac{1}{\sqrt{N}}\sum_{i=1}^n \left({\bf X}_i{\bf X}_i'e_i^2-E[{\bf X}{\bf X}'e^2]\right)-\mathbb{B}_N [{I_p} \otimes (\widehat{\mathbb{Q}}_N^{-1}\frac{1}{\sqrt{N}}\sum_{i=1}^n {\bf X}_i e_i)]+O_p({1}/{{N}}),
\end{multline*}
	where $\mathbb{B}_N$ is the matrix $\frac{2}{N}\sum_{i=1}^N{\bf X}_i{\bf X}_i'\big(Y_i-{\bf X}_i'{{\bs \beta}}\big)\otimes {\bf X}_i'$. Condition (i) implies that the second moment of ${\bf X}_i{\bf X}_i'e_i^2$ is finite. Hence,  by the central limit theorem, 
	\begin{equation}\label{eq:joint_normal}
	\frac{1}{\sqrt{N}} \Big( \sum_{i=1}^n \left\{{\bf X}_i{\bf X}_i'e_i^2 - E[{\bf X} {\bf X}' e^2]\right\} ,\sum_{i=1}^n {\bf X}_i e_i\Big)
	\end{equation}
	is jointly asymptotically normal, and since $\mathbb{B}_N$ and $\widehat{\mathbb{Q}}_N$ converge to a matrix of constants, a version of Slutsky's theorem implies that $\sqrt{N}(\widehat{\mathbb{W}}_N-\mathbb{W})$ is asymptotically normal.

	Another application of Slutsky's theorem implies that $\sqrt{N}(\widehat{\mathbb{W}}_N\widehat{\mathbb{Q}}_N^{-1}- {\mathbb{W}} {\mathbb{Q}}^{-1})$
	is asymptotically normal (since $\widehat{\mathbb{Q}}_N \rightarrow \mathbb{Q}$ in probability)  and therefore so is
	$$\sqrt{N}\left\{tr\left(\widehat{\mathbb{W}}_N\widehat{\mathbb{Q}}_N^{-1}\right)
	-tr(\mathbb{W}\mathbb{Q}^{-1})\right\}=\sqrt{N}\left\{tr(\widehat{\mathbb{V}}_N)
	-tr(\mathbb{V})\right\}.$$ It follows that $tr(\widehat{\mathbb{V}}_N)-tr(\mathbb{V})=O_p(1/\sqrt{N})$ and it is asymptotically normal.
	
	Similar to previous arguments, the random variables in \eqref{eq:joint_normal} and the first part of $\cal{E}_N$ are jointly asymptotically normal and a version of Slutsky's theorem that allows us to ignore the term $O_p(1/N)$  together with \eqref{eq:final} and \eqref{eq:CAPEPS}, implies the last statement of Theorem \ref*{thm:o1n} about the asymptotic normality of $\sqrt{N}\big(C^{({\mathpzc P})}(n,N)- AR(n,{\mathpzc P})\big )$. \qed

	\noindent{\bf Proof of Lemma \ref{lem:normal}}. First assume that the first coordinate of the covariate vectors is 1 (corresponding to an intercept coefficient). 
	Let $\widetilde{\mathbb{X}}$ denote the $n\times (p-1)$ matrix defined as
	$\mathbb{X}_n$ but without the first column of 1's. We suppress $n$ here and in the following notation.  Let $\overline {\bf X}:=\widetilde{\mathbb{X}'}{\large \bs 1}/n$, where ${\large \bs 1}$ is an $n$-column vector of 1's, that is, $\overline {\bf X}$ is  the $(p-1)$-column of covariates means, and let  $\mathbb{S}:=\widetilde{\mathbb{X}}'\widetilde{\mathbb{X}}/n-
	\overline {\bf X}\,\overline {\bf X}'$. Then by \citet{HJ}, page 25, Equation (0.8.5.6)
	\begin{equation}\label{eq:HornJ}
	(\mathbb{X}'\mathbb{X}/n)^{-1}=\begin{bmatrix} 
	1+\overline {\bf X}'\mathbb{S}^{-1}\overline{\bf X} & -\overline {\bf X}'\mathbb{S}^{-1}\\
	-\mathbb{S}^{-1}\overline{\bf X} & \mathbb{S}^{-1}\\
	\end{bmatrix}.
	\end{equation}
	Let $\widetilde {\bf X}$ denote a $(p-1)$-vector of the covariates without the first 1, and assume first that $\widetilde {\bf X} \sim N({\bs \mu}, \Sigma)$. Then $n\mathbb{S} 
	\sim Wishart_{p-1}(\Sigma,n-1)$, and $\overline {\bf X}$ and $\mathbb{S}$ are independent. The third moments of $\mathbb{S}^{-1}$ are uniformly (in $n$) bounded  by $C\max\{1, 1/\lambda_{min}(\Sigma)^3\}$ for some $C>0$, provided $n-p-5>0$ by 
	Theorem 4.1  of \citet{von_Rosen}. By \eqref{eq:HornJ} the third moments of $(\mathbb{X}'\mathbb{X}/n)^{-1}$ are also bounded. If $\widetilde {\bf X}$ is distributed according to a mixture of normals, the assumption in the lemma that all these normal distributions  have $\lambda_{min}(\Sigma)>c>0$ implies the uniform boundedness for the mixture. 
	
	If the first coordinate is not 1, we append 1 to the covariates and now the matrix of interest in the above notation (with $p$ replacing $p-1$) is $\widetilde{\mathbb{X}}'\widetilde{\mathbb{X}}/n$, which is now $p \times p$, and we wish to show that its inverse has bounded third moments. The eigenvalues of the latter matrix are larger (not strictly) than those of $\mathbb{S}=\widetilde{\mathbb{X}}'\widetilde{\mathbb{X}}/n-
	\overline {\bf X}\,\overline {\bf X}'$. The latter relation is revered for the inverses. Now use the inequality that for a positive definite matrix $A$ we have $|a_{ij}| \le tr(A)/2$ to conclude that the entries of $(\widetilde{\mathbb{X}}'\widetilde{\mathbb{X}}/n)^{-1}$ are bounded by $tr(\mathbb{S}^{-1})$. By the first part of the theorem (with $p$ replacing $p-1$) we know that 
	$\mathbb{S}^{-1}$ has finite third moments and therefore also its trace (by Minkowski inequality), and the result follows. 
	\qed

	\begin{thma}%label by hand
		For $\mathbb{X}_n \in H_n$ (see \eqref{eq:setGn}) the matrix $(\mathbb{X}_n'\mathbb{X}_n/n)^{-1}$ exists and all its entries are bounded uniformly in $n$.
		Moreover, if the components of $\bf X$ are bounded, then for some $a<1$ we have, $P(H_n) >1- a^{n\lambda_{min}(\mathbb{Q})}$, which converge to 1 at an exponential rate in $n$.
	\end{thma}

	\noindent{\bf Proof}.
	When $\mathbb{X}_n \in H_n$, then $\lambda_{min}(\mathbb{X}_n'\mathbb{X}_n/n)
	\ge \frac{1}{2}\lambda_{min}(\mathbb{Q})>0$, and therefore $(\mathbb{X}_n'\mathbb{X}_n/n)^{-1}$ exists. Since the entries of a positive semi-definite matrix are bounded by the maximal eigenvalue, all entries of $(\mathbb{X}_n'\mathbb{X}_n/n)^{-1}$ are bounded in this case by $2/\lambda_{min}(\mathbb{Q})$. For the moreover part, notice that when   
	the elements of ${\bf X}$ are bounded then so is $\lambda_{max}(\mathbb{X}_n'\mathbb{X}_n/n)$.
	By \citet{Tropp}, Theorem 1.1,  $P(H^c_n)\le a^{n\lambda_{min}(\mathbb{Q})}$ for some $a<1$. \qed
	
	\noindent{\bf Proof of Theorem \ref{thm:BeGood}}.
	Conditionally on $H_n$ the arguments in the proof of Theorem \ref{thm:o1n} continue to apply with obvious modifications. Lemma 2.3 provides an exponentially small bound on $P(H^c_n)$, adding a term $O_p(e^{-\gamma n})$ in Equation \eqref{eq:approx_pred3} and $O_p(e^{-\gamma N})$ to \eqref{eq:final} and \eqref{eq:further} (b)  for some $\gamma>0$. Clearly, these terms do not affect the results.  \qed

	\noindent{\bf Proof  of Proposition \ref{prop:single}}.  For Part (i),  set ${\mathpzc P}(n) \in \pi^*(n)$ and $\widetilde {\mathpzc P}(n) \in {\mathpzc P}^*(n)$. 
	We have $R(n,{\mathpzc P}(n))-R(n,\widetilde {\mathpzc P}(n)) \ge 0$, and also 
	$R(n,{\mathpzc P}(n))-R(n,\widetilde {\mathpzc P}(n))= [R(n,{\mathpzc P}(n))-AR(n,{\mathpzc P}(n))]+ 
	[AR(n,{\mathpzc P}(n))-AR(n,\widetilde {\mathpzc P}(n))]+[AR(n,\widetilde {\mathpzc P}(n))-R(n,\widetilde{\mathpzc P}(n))]$. The  middle term is negative and by \eqref{eq:approx_pred3} the two other terms are $o(1/n)$ and Part (i) follows.
	
	For Part (ii), we first show that ${\mathpzc P}^*$ is a singleton. Suppose that there are two models ${\mathpzc P},{\mathpzc Q}$ in ${\mathpzc P}^*$.  By the definition of ${\mathpzc P}^*$, the components of the projection coefficient vectors ${\bs \beta}^{({\mathpzc P})}$ and ${\bs \beta}^{({\mathpzc Q})}$ must all be non-zero. Since the function $(Y-a)^2$ is strictly convex in $a$, we have  
	\begin{equation}\label{eq:convex}
	\left(Y-\frac{{\bf X}^{({\mathpzc P})'}{{\bs \beta}^{({\mathpzc P})}}+{\bf X}^{({\mathpzc Q})'}{{\bs \beta}^{({\mathpzc Q})}}}{2}\right)^2 \le \frac{\big(Y-{\bf X}^{({\mathpzc P})'}{{\bs \beta}^{({\mathpzc P})}}\big)^2+\big(Y-{\bf X}^{({\mathpzc Q})'}{{\bs \beta}^{({\mathpzc Q})}}\big)^2}{2},
	\end{equation}
	with equality if and only if ${\bf X}^{({\mathpzc P})'}{{\bs \beta}^{({\mathpzc P})}}={\bf X}^{({\mathpzc Q})'}{{\bs \beta}^{({\mathpzc Q})}}$.  Unless ${\bf X}^{({\mathpzc P})'}{{\bs \beta}^{({\mathpzc P})}}={\bf X}^{({\mathpzc Q})'}{{\bs \beta}^{({\mathpzc P})}}$ a.s., the model ${\mathpzc P} \cup {\mathpzc Q}$ would contradict the assumption that  ${\mathpzc P}$ and ${\mathpzc Q}$ are in $\cal{M}$ by taking expectations in \eqref{eq:convex}.
	We assumed that $E({\bf X} {\bf X}')$ is invertible and hence  ${\bf X}'{\bs \beta}$ vanishes a.s. only for ${\bs \beta}=\bf 0$. Adding zeros to ${\bs \beta}^{({\mathpzc P})}$ and ${\bs \beta}^{({\mathpzc Q})}$, thus completing them to vectors in $\mathbb{R}^d$, we see that the completed vectors are identical.  Since the two models are in $\cal M$ and their projection coefficients are all non-zero, it follows that  ${\mathpzc P} = {\mathpzc Q}$, and therefore ${\mathpzc P}^*$ is a singleton. The above discussion also shows that any model ${\mathpzc Q} \in \cal M$ must satisfy, as sets of covariates, ${\mathpzc Q} \supseteq {\mathpzc P}^*$.
	
	In order to show that $\pi^*(n)={\mathpzc P}^*$ for large $n$, note first that ${\mathpzc Q} \in \pi^*(n)$ for  large enough $n$ implies ${\mathpzc Q} \in \cal M$; if not
	then there is some $\widetilde {\mathpzc Q} \in \cal M$ such that
	$E\big(Y-{\bf X}^{(\widetilde {\mathpzc Q})'}{{\bs \beta}^{(\widetilde {\mathpzc Q})}}\big)^2 <  E\big(Y-{\bf X}^{( {\mathpzc Q})'}{{\bs \beta}^{({\mathpzc Q})}}\big)^2$. For large $n$ this $\widetilde q$ contradicts ${\mathpzc Q} \in \pi^*(n)$. It follows that ${\mathpzc Q} \supseteq {\mathpzc P}^*$ as sets of covariates, and it suffices to show that $tr(\mathbb{V}^{({\mathpzc P})})$ is minimized over $\cal M$ by ${\mathpzc P}^*$. Indeed we show that 
	if ${\mathpzc P},{\mathpzc Q}\in \cal M$ and ${\mathpzc Q} \supsetneqq {\mathpzc P}$ as sets of covariates, then $tr(\mathbb{V}^{({\mathpzc P})}) < tr(\mathbb{V}^{({\mathpzc Q})})$.  Consider a Gram--Schmidt process on the space of square integrable random variables, with the inner product of two random variables being the expectation of their product. Starting with the indexes in ${\mathpzc P}$, there exist  linear transformations $\widetilde{\bf X}^{({\mathpzc P})}:=A{\bf X}^{({\mathpzc P})}$ and $\widetilde{\bf X}^{({\mathpzc Q})}:=B{\bf X}^{({\mathpzc Q})}$, where $A$ and $B$ are invertible $p \times p$ and $ q\times q$ matrices, such that $E\big(\widetilde{\bf X}^{({\mathpzc P})} \widetilde{\bf X}^{({\mathpzc P})'} \big)$ and $E\big(\widetilde{\bf X}^{({\mathpzc Q})} \widetilde{\bf X}^{({\mathpzc Q})'} \big)$ are both identity matrices (with different dimensions). Also, we can assume that the first $p$ rows of $B$ can be obtained from those of $A$ by adding $q-p$ zeros to each of these rows.  
	Therefore, for $k \in {\mathpzc P}$ we have $\widetilde{X}^{({\mathpzc P})}_k =\widetilde{X}^{(\mathpzc Q)}_k$, where $\widetilde{X}^{({\mathpzc P})}_k$ is the $k$th coordinate of  $\widetilde{\bf X}^{({\mathpzc P})}$. The relation $\widetilde{\bf X}^{({\mathpzc P})}=A{\bf X}^{({\mathpzc P})}$ and  straightforward algebra, using properties of the trace function,   imply that $tr(\mathbb{V}^{({\mathpzc P})})=tr(\widetilde{\mathbb{V}}^{({\mathpzc P})})$, where  
	\[
	\widetilde{\mathbb{V}}^{({\mathpzc P})}:=  E\left(\widetilde{\bf X}^{({\mathpzc P})} \widetilde{\bf X}^{({\mathpzc P})'} \{e^{({\mathpzc P})}\}^2 \right)\left\{E\left(\widetilde{\bf X}^{({\mathpzc P})} \widetilde{\bf X}^{({\mathpzc P})'} \right)\right\}^{-1}= E\left(\widetilde{\bf X}^{({\mathpzc P})} \widetilde{\bf X}^{({\mathpzc P})'} \{e^{({\mathpzc P})}\}^2 \right),
	\]    
	and similarly for $tr(\mathbb{V}^{({\mathpzc Q})})$. We have 
	$e^{({\mathpzc P})}=e^{({\mathpzc Q})}$ for any ${\mathpzc P},{\mathpzc Q}\in \cal M$ (with probability 1) since otherwise,  by the argument in \eqref{eq:convex},  ${\mathpzc P} \cup {\mathpzc Q}$ would contradict the assumption that both ${\mathpzc P}$ and ${\mathpzc Q}$ are in $\cal M$ as above. We conclude that,
	\begin{multline}\label{eq:trace-gram}
	tr(\mathbb{V}^{({\mathpzc P})}) = tr\left\{E\left(\widetilde{\bf X}^{({\mathpzc P})} \widetilde{\bf X}^{({\mathpzc P})'} \{e^{({\mathpzc P})}\}^2 \right)\right\}=\sum_{k \in {\mathpzc P}} E\left(\widetilde{X}_k^{({\mathpzc P})} e^{({\mathpzc P})}\right)^2\\ < \sum_{k \in {\mathpzc Q}} E\left(\widetilde{X}_k^{({\mathpzc Q})} e^{({\mathpzc Q})}\right)^2 = tr(\mathbb{V}^{({\mathpzc Q})}).
	\end{multline}
	The strict inequality follows from the fact that  $\mathbb{W}$ is positive definite, and thus $\widetilde{\bf X}^{({\mathpzc P})} \widetilde{\bf X}^{({\mathpzc P})'}\{e^{({\mathpzc P})}\}^2=A{\bf X}^{({\mathpzc P})} {\bf X}^{({\mathpzc P})'}\{e^{({\mathpzc P})}\}^2A'$ are matrices with positive definite expectations, and therefore positive diagonal elements.
	Summing up, the above discussion shows that ${\mathpzc P}^*$ is a singleton, and that $tr(\mathbb{V}^{({\mathpzc P}^*)})$ is minimal among the models in ${\cal M}$. This implies that $\pi^*(n) \to {\mathpzc P}^*$ as $n \to \infty$. By \eqref{eq:approx_pred3}, for every model $p$, $R(n,{\mathpzc P})-AR(n,{\mathpzc P})=o(1/n)$, and therefore ${\mathpzc P}^*(n)=\pi^*(n)$ for large enough $n$. Hence, also ${\mathpzc P}^*(n) \to {\mathpzc P}^*$ as $n \to \infty$. 		\qed

	\noindent{\bf Proof of Proposition \ref{prop:cons}}.
	It suffices to prove that   $P( \widehat{\pi^*}(n,N)={\mathpzc P}^*) \to 1$ when $n, N\to \infty$  and $n/N \to 0$,
	since by Proposition \ref{prop:single},  ${\mathpzc P}^*(n)={\mathpzc P}^*$ for large enough $n$.
	Equivalently, we claim that for every ${\mathpzc P} \ne {\mathpzc P}^*$ we have $P( \widehat{\pi^*}(n,N)={\mathpzc P}) \to 0$, and since there is a finite number of models, the result follows. 
	The latter claim is proved separately for  ${\mathpzc P} \notin {\cal M}$ and then for ${\mathpzc P} \in {\cal M}$,  (conditions (a) and (b) below):\\
	(a) For ${\mathpzc P} \notin {\cal M}$ we shall show that   
	\begin{equation}\label{eq:prop.a}
	C^{({\mathpzc P})}(n,N)-C^{({\mathpzc P}^*)}(n,N)=A-\frac{tr(\mathbb{V}^{({\mathpzc P})})-tr(\mathbb{V}^{({\mathpzc P}^*)})}{n}+O_p(1/\sqrt{N}),
	\end{equation}
	for a positive constant $A$. Since $\widehat{\pi^*}(n,N)$ is the minimizer of $C^{({\mathpzc P})}(n,N)$, it follows that $P( \widehat{\pi^*}(n,N)={\mathpzc P}) \to 0$ as both $n,N$ go to infinity.
	To prove \eqref{eq:prop.a}
	note that by the definition of $AR(n,{\mathpzc P})$ and Equations \eqref{eq:approx_predara},  \eqref{eq:final}, and \eqref{eq:further}, we have
	\begin{align*}\label{eq:AR_diff}
	C^{({\mathpzc P})}(n,N)-C^{({\mathpzc P}^*)}(n,N)&= AR(n,{\mathpzc P})-AR(n,{\mathpzc P}^*) + O_p(1/\sqrt{N}) \\&=E\big(Y-{\bf X}^{({\mathpzc P})'}{{\bs \beta}^{(p)}}\big)^2 - E\big(Y-{\bf X}^{({\mathpzc P}^*)'}{{\bs \beta}^{({\mathpzc P}^*)}}\big)^2\\&+\frac{tr(\mathbb{V}^{({\mathpzc P})})-tr(\mathbb{V}^{({\mathpzc P}^*)})}{n} +O_p(1/\sqrt{N}).
	\end{align*} 
	Since ${\mathpzc P} \notin {\cal M}$ and ${\mathpzc P}^* \in {\cal M}$, the difference of the expectations,  which we denote by $A$,  is positive.

	(b) For ${\mathpzc P} \in {\cal M}$ and ${\mathpzc P} \ne {\mathpzc P}^*$ we shall show that
	\begin{equation}\label{eq:prop.b}
	C^{({\mathpzc P})}(n,N)-C^{({\mathpzc P}^*)}(n,N)=B/n + O_p(1/{N})+O_p\left(\frac{1}{n\sqrt{N}}\right),
	\end{equation}
	where $B$ is a positive constant implying that $P( \widehat{\pi^*}(n,N)={\mathpzc P}) \to 0$ when both $n,N$ go to infinity and $n/N \to 0$.
	
	We now prove \eqref{eq:prop.b}.
	Consider ${\mathpzc P} \in {\cal M}$ and ${\mathpzc P} \ne {\mathpzc P}^*$.  Since both models are in ${\cal M}$, we have
	\[
	E\big(Y-{\bf X}^{({\mathpzc P})'}{{\bs \beta}^{({\mathpzc P})}}\big)^2 - E\big(Y-{\bf X}^{({\mathpzc P}^*)'}{{\bs \beta}^{({\mathpzc P}^*)}}\big)^2=0
	\]
	and therefore
	\[
	AR(n,{\mathpzc P})-AR(n,{\mathpzc P}^*) =\frac{tr(\mathbb{V}^{({\mathpzc P})})-tr(\mathbb{V}^{({\mathpzc P}^*)})}{n}
	.\]
	In Proposition \ref{prop:single} we showed that $tr(\mathbb{V}^{({\mathpzc P}^*)}) < tr(\mathbb{V}^{({\mathpzc P})})$, and therefore $AR(n,{\mathpzc P})-AR(n,{\mathpzc P}^*) =B/n$, where $B$ is a positive constant.
	Since both ${\mathpzc P}$ and ${\mathpzc P}^*$ are in ${\cal M}$, it follows that ${\bf X}^{({\mathpzc P})'}{{\bs \beta}^{({\mathpzc P})}}={\bf X}^{({\mathpzc P}^*)'}{{\bs \beta}^{({\mathpzc P}^*)}} $ a.s. (see the proof of Proposition \ref{prop:single}). Therefore, the first part of ${\cal E}_N^{({\mathpzc P})}$ and ${\cal E}_N^{({\mathpzc P}^*)}$ is equal, and hence,
	${\cal E}_N^{({\mathpzc P})} - {\cal E}_N^{({\mathpzc P}^*)}=O_p(1/N)$.
	Recalling that $AR(n,{\mathpzc P})-AR(n,{\mathpzc P}^*) =B/n$, \eqref{eq:final} implies that
	\begin{multline*}
	C^{({\mathpzc P})}(n,N)-C^{({\mathpzc P}^*)}(n,N)
	\\=B/n +O_p(1/N)-\frac{tr(\mathbb{V}^{({\mathpzc P})})-tr(\widehat{\mathbb{V}}_N^{({\mathpzc P})})-tr(\mathbb{V}^{({\mathpzc P}^*)})+tr(\widehat{\mathbb{V}}_N^{({\mathpzc P}^*)})}{n},
	\end{multline*}
	which implies \eqref{eq:prop.b} by \eqref{eq:further} (b). \qed

	\noindent{\bf Proof of Theorem \ref{prop:C_pj1}}.
	The first part follows from \eqref{eq:final} of Theorem \ref{thm:o1n}.
	That the $o_p$ terms do not depend on $n$ can be seen by inspecting the proof of \eqref{eq:final} of Theorem \ref{thm:o1n}. The moreover part follows from the asymptotic normality of each $j$; see \eqref{eq:further2}. \qed

	\noindent{\bf Proof of Proposition \ref{prop:cons_J_fixed}}. Part 1 follows from the first part of Theorem \ref{prop:C_pj1}. 
	
	The proof of Part 2 differs from that of Proposition \ref{prop:cons}  only in taking averages over ${\cal J}$ in similar expressions. The only real difference is in case (b) of the proof of Proposition \ref{prop:cons}, with ${\mathpzc P}$ and $\contour[2]{black}{$\mathpzc{P}$}^*$ both in $\bs{\cal M}$. The proof is achieved by showing that there exists ${\cal B}>0$ such that 
	\begin{equation}\label{eq:BCDE}
	\lim \sup_{n/N_j \le C, n,{\bf N} \to \infty} P\left(n\{ {\bf C}^{({\mathpzc P})}(n,{\bf N}) - {\bf C}^{(\contour[2]{black}{$\mathpzc{P}$}^*)}(n,{\bf N})\}<{\cal B}/2\right)<K/{\cal J}.
	\end{equation} 
	We have
	\begin{equation}\label{eq:long}
	n\{ {\bf C}^{({\mathpzc P})}(n,{\bf N}) - {\bf C}^{(\contour[2]{black}{$\mathpzc{P}$}^*)}(n,{\bf N}) \}=\frac{1}{\cal J} \sum_{j=1}^{\cal J} B_j + \frac{1}{\cal J} \sum_{j=1}^{\cal J} C_j + \frac{1}{\cal J} \sum_{j=1}^{\cal J} D_j+ \frac{1}{\cal J} \sum_{j=1}^{\cal J} E_j
	\end{equation} 
	where,
	\begin{align*}
	B_j&:=tr(\mathbb{V}_j^{({\mathpzc P})})-tr(\mathbb{V}_j^{(\contour[2]{black}{$\mathpzc{P}$}^*)}),
	\qquad C_j:=n({\cal E}_{j,N_j}^{({\mathpzc P})}-{\cal E}_{j,N_j}^{(\contour[2]{black}{$\mathpzc{P}$}^*)}), \\
	&\,\, D_j:={tr}(\mathbb{V}_j^{({\mathpzc P})})-{tr}(\widehat{\mathbb{V}}^{({\mathpzc P})}_{j,N_j})-{tr}(\mathbb{V}_j^{(\contour[2]{black}{$\mathpzc{P}$}^*)})+{tr}(\widehat{\mathbb{V}}^{(\contour[2]{black}{$\mathpzc{P}$}^*)}_{j,N_j}),
	\end{align*}
	and $E_j=n o_p(1/N_j)$, arising from the last term in \eqref{eq:final}. 
	The proof of \eqref{eq:BCDE} is accomplished by showing that 
	$\frac{1}{\cal J} \sum_{j=1}^{\cal J} B_j \ge \cal B$, to be defined below, and that the other three sums are small.

	We start with the first term in \eqref{eq:long},
	 $\frac{1}{\cal J} \sum_{j=1}^{\cal J} B_j=\frac{1}{{\cal J}}\sum_j[{tr(\mathbb{V}_j^{({\mathpzc P})})-tr(\mathbb{V}_j^{(\contour[2]{black}{$\mathpzc{P}$}^*)}]}$.
	Since  ${\mathpzc P}$ is in $\bs{\cal M}$ we have ${\mathpzc P} \supseteq \contour[2]{black}{$\mathpzc{P}$}^*$ as sets of covariates.  By \eqref{eq:trace-gram}, 
	$tr(\mathbb{V}_j^{({\mathpzc P})})-tr(\mathbb{V}_j^{(\contour[2]{black}{$\mathpzc{P}$}^*)})$ is bounded below by $E_{G_j}\left(\widetilde{X}_k^{({\mathpzc P})} e^{({\mathpzc P})}\right)^2$ for $k \in {\mathpzc P} \smallsetminus \contour[2]{black}{$\mathpzc{P}$}^*$ (as sets). We have $\widetilde{X}_k= {\bf b}'_k {\bf X}$ where ${\bf b}'_k$ is the $k$th row of the matrix $B$ defined in the proof  of Proposition \ref{prop:single}. We have  $1=E({\bf b}'_k {\bf X})^2={\bf b}_k'\mathbb{Q}_j^{({\mathpzc P})}{\bf b}_k$ and therefore $\| {\bf b}_k'\{\mathbb{Q}_j^{({\mathpzc P})}\}^{1/2}\|=1$. 
	It follows that $\|{\bf b}_k\|^2={\bf b}_k'\{\mathbb{Q}_j^{({\mathpzc P})}\}^{1/2}\{\mathbb{Q}_j^{({\mathpzc P})}\}^{-1}\{\mathbb{Q}_j^{({\mathpzc P})}\}^{1/2}{\bf b}_k \ge
	\lambda_{min}(\{\mathbb{Q}_j^{({\mathpzc P})}\}^{-1}) 
	= 1/\lambda_{max}(\mathbb{Q}_j^{({\mathpzc P})})$ and therefore
	$E_{G_j}\big(\widetilde{X}_k^{({\mathpzc P})} e^{({\mathpzc P})}\big)^2=E_{G_j}({\bf b}_k'{\bf X}^{({\mathpzc P})}e)^2 
	={\bf b}'\mathbb{W}_j^{({\mathpzc P})}\widetilde {\bf b}_k \ge \lambda_{min}(\mathbb{W}_j^{({\mathpzc P})})/\lambda_{max}
	(\mathbb{Q}_j^{({\mathpzc P})})>1/{C}^2>0$. 
	We obtained that $\frac{1}{\cal J} \sum_{j=1}^{\cal J} B_j \ge 1/C^2=:\cal B$.

	We now deal with $\frac{1}{\cal J} \sum_{j=1}^{\cal J} C_j$.
	By \eqref{eq:final} and \eqref{eq:CAPEPS} and the fact that ${\mathpzc P}, {\mathpzc P}^* \in \bs{\cal M}$, this term equals, 
	\begin{multline}\label{eq:difff}
	\frac{n}{{\cal J}}\Big[\sum_j\frac{1}{N_j}\Big\{tr\left[{\bf U}^{(\contour[2]{black}{$\mathpzc{P}$}^*)}_{j,N_j}{\bf U}^{(\contour[2]{black}{$\mathpzc{P}$}^*)'}_{j,N_j}(\mathbb{Q}_j^{(\contour[2]{black}{$\mathpzc{P}$}^*)})^{-1}\right]-tr({\mathbb{V}}^{(\contour[2]{black}{$\mathpzc{P}$}^*)}_j)\\-tr\left[{\bf U}^{(p)}_{j,N_j}{\bf U}^{({\mathpzc P})'}_{j,N_j}(\mathbb{Q}_j^{({\mathpzc P})})^{-1}\right]+tr({\mathbb{V}}^{({\mathpzc P})}_j)\Big\}\Big].
	\end{multline}
	Since ${\bf U}^{({\mathpzc P})}_{j,N_j}$ converges in distribution to 
	${\bf Z}^{({\mathpzc P})}_{j} \sim  N(0,\mathbb{W}_j^{({\mathpzc P})})$ and $C$ is an upper bound on ${n}/{N_j}$ we have that the limit of the probability that the  expression in \eqref{eq:difff} exceeds $\varepsilon$ is bounded by  $P(T_{\cal J}>\varepsilon)$
	where
	\begin{multline*}%\label{eq:difff1}
	T_{\cal J}:=\frac{C}{{\cal J}}\Big|\sum_j\Big\{tr\Big[{\bf Z}^{(\contour[2]{black}{$\mathpzc{P}$}^*)}_{j}{\bf Z}^{(\contour[2]{black}{$\mathpzc{P}$}^*)'}_{j}(\mathbb{Q}_j^{(\contour[2]{black}{$\mathpzc{P}$}^*)})^{-1}\Big]-tr({\mathbb{V}}^{(\contour[2]{black}{$\mathpzc{P}$}^*)}_j)\\-tr\left[{\bf Z}^{({\mathpzc P})}_{j}{\bf Z}^{({\mathpzc P})'}_{j}(\mathbb{Q}_j^{({\mathpzc P})})^{-1}\right]+tr({\mathbb{V}}^{({\mathpzc P})}_j)\Big\} \Big|.
	\end{multline*}	
	Note that the expression within the absolute value sign has mean zero.  Writing  $T_{\cal J} = \frac{C}{{\cal J}} |\sum_{j=1}^{\cal J} A_j|$, Markov's inequality implies that in order to obtain $P(T_{\cal J}>\varepsilon) \le K/{\cal J}$ it is enough to bound $Var{(A_j)}$ uniformly in $j$, which holds  when $(\mathbb{Q}_j^{({\mathpzc P})})^{-1}, \mathbb{W}_j^{({\mathpzc P})}$ are bounded (element-wise) for all models $\mathpzc P$ (of which there is a finite number) and uniformly for all $j$. This follows from our eigenvalue assumptions (see \eqref{eq:Cbound}) and the fact that the entries of a positive-definite matrix are bounded by its maximal eigenvalue.
	Finally, it suffices to show that  $D_j \to 0$ and $E_j \to 0$ as $n, N_j \to \infty$ with $n/N_j$ bounded. The  first follows from  \eqref{eq:further}, and the second is obvious. 
	\qed  
	
	\noindent{\bf Proof of Lemma \ref{lem:tild_AR}}.
	First notice that when the moments appearing in (i) of Theorem \ref{thm:o1n} are bounded uniformly in  $\theta \in \Theta$, then $E_{G_\theta}(Y-{\bf X}'{\bs \beta}_\theta)^2$ is bounded in $\theta$. Also, the matrix $\mathbb{W}_\theta$ is bounded (element-wise) uniformly in $\theta$. Finally, because $ E\{{(\mathbb{X}_n}'\mathbb{X}_n/n)^{-1}\} - \mathbb{Q}_\theta^{-1}$ is positive semi-definite (see \citet{Groves}), then uniform boundedness of the moment condition (ii) of Theorem \ref{thm:o1n} implies that $\mathbb{Q}_\theta^{-1}$ is uniformly bounded and therefore so is $tr(\mathbb{V}_\theta)=tr(\mathbb{W}_\theta \mathbb{Q}_\theta^{-1})$.

	We have 
	\begin{multline}\label{eq:useCLT}
	\textbf{AR}_{pop}(n,{\mathpzc P})-{\bf AR}(n,{\mathpzc P})\\=\Big\{\int E_{G_{\theta}}(Y-{\bf X}' {\bs \beta}_{\theta})^2
	{\mathscr P}(d\theta) - \frac{1}{{\cal J}} \sum_{j=1}^{\cal J} E_{G_{j}}(Y-{\bf X}' {\bs \beta}_{j})^2\Big\}\\
	+\frac{1}{n} \Big\{ \int tr({\mathbb{V}}_{\theta}) {\mathscr P}(d\theta) - \frac{1}{{\cal J}} \sum_{j=1}^{\cal J} tr({\mathbb{V}}_{j}) \Big\}.
	\end{multline}
	The above two sums contain random variables that are bounded, and hence so are their 
	variances. The central limit theorem applied twice, implies \eqref{eq:R_j1} and the claimed asymptotic normality. It is easy to see directly from \eqref{eq:useCLT}  that the $O_p$ term in \eqref{eq:R_j1} is uniform in $n$.\qed

	For the proof of Proposition \ref{prop:popolo} we need the following lemma:
	
	\begin{lemma}\label{lem:pop}
		Suppose that the conditions of Lemma \ref{lem:unif} hold and also that $\lambda_{min}(\mathbb{W}_\theta)$ is bounded away from zero uniformly in $\theta$; then
		\begin{enumerate}
			\item The set  $\contour[2]{black}{$\mathpzc{P}$}^*_{pop}$ is a singleton and as ${n \to \infty}$ both ${\bs \pi}^*_{pop}(n) \to \contour[2]{black}{$\mathpzc{P}$}^*_{pop}$
			and  ${\contour[2]{black}{$\mathpzc{P}$}}^*_{pop}(n) \to \contour[2]{black}{$\mathpzc{P}$}^*_{pop}$, 
			and therefore also ${\bs \pi}^*_{pop}(n) = \contour[2]{black}{$\mathpzc{P}$}^*_{pop}(n)$ for large $n$.
			\item There exists a constant $K_C$ depending only on $C$, such that for ${\bs \pi}^*(n)$ defined in \eqref{eq:pi_J}, 
			\[
			{P}\Big({\bs \pi}^*(n) \subseteq {\bs \pi}_{pop}^*(n)\Big)\ge 1-\frac{K_C}{\cal J}~~\forall n.
			\] 
			\end{enumerate}
	\end{lemma}
	\noindent{\bf Proof of Lemma \ref{lem:pop}}.
	\underline{Part 1.} The proof is similar to that of Proposition \ref{prop:single}. We sketch the proof. Let ${\mathpzc P}$ and ${\mathpzc Q}$ be in $\contour[2]{black}{$\mathpzc{P}$}^*_{pop}$. By convexity as in \eqref{eq:convex},
	\begin{equation}\label{eq:convex1}
	\frac{\big(Y-{\bf X}^{({\mathpzc P})'}{{\bs \beta}^{({\mathpzc P})}}\big)^2+\big(Y-{\bf X}^{({\mathpzc Q})'}{{\bs \beta}^{({\mathpzc Q})}}\big)^2}{2} -\left(Y-\frac{{\bf X}^{({\mathpzc P})'}{{\bs \beta}^{({\mathpzc P})}}+{\bf X}^{({\mathpzc Q})'}{{\bs \beta}^{({\mathpzc Q})}}}{2}\right)^2 \ge 0, 
	\end{equation}
	with equality iff ${\bf X}^{({\mathpzc P})'}{{\bs \beta}^{({\mathpzc P})}}={\bf X}^{({\mathpzc Q})'}{{\bs \beta}^{({\mathpzc P})}}$. This implies that $\contour[2]{black}{$\mathpzc{P}$}^*_{pop}$ is a singleton as in the proof of Proposition \ref{prop:single}, Part (ii). Since ${\mathpzc P}$ and ${\mathpzc Q}$ are in $\bs{\cal M}_{pop}$, the  expectation of the left-hand side of \eqref{eq:convex1} is zero. It follows that 
	$
	\int P_{G_\theta}\Big({\bf X}^{({\mathpzc P})'}{{\bs \beta}^{({\mathpzc P})}}={\bf X}^{({\mathpzc Q})'}{{\bs \beta}^{({\mathpzc Q})}}\Big){\mathscr P}(d\theta)=1,
	$ 
	and therefore for every model ${\mathpzc P}$ in $\bs{\cal M}_{pop}$ we have that $\contour[2]{black}{$\mathpzc{P}$}^*_{pop} \subseteq {\mathpzc P}$.  By  the assumptions on moments being uniformly bounded, it follows  that $\lambda_{max}(\mathbb{Q}_\theta)$ is bounded above and $\lambda_{min}(\mathbb{W}_\theta)$ is positive and bounded  away from zero, both uniformly in $\theta$. Now \eqref{eq:trace-gram} and the discussion in the paragraph above \eqref{eq:difff} imply that if $\contour[2]{black}{$\mathpzc{P}$}^*_{pop} \subseteq {\mathpzc P}$ as sets of covariates,  ${\mathpzc P} \in \bs{\cal M}_{pop}$, and $\contour[2]{black}{$\mathpzc{P}$}^*_{pop} \neq {\mathpzc P} $ then
	$\int tr\Big(\mathbb{V}_\theta^{(\contour[2]{black}{$\mathpzc{P}$}^*_{pop})}\Big) {\mathscr P}(d\theta) < \int tr\left(\mathbb{V}_\theta^{({\mathpzc P})}\right) {\mathscr P}(d\theta)$. Therefore, $\contour[2]{black}{$\mathpzc{P}$}^*_{pop}$ has a minimal trace among $\bs{\cal M}_{pop}$. It follows that ${\bs \pi}^*_{pop}(n) \to \contour[2]{black}{$\mathpzc{P}$}^*_{pop}$ as ${n \to \infty}$.\\
	Furthermore, Lemma \ref{lem:unif} implies that ${\bs \pi}^*_{pop}(n)$ and $\contour[2]{black}{$\mathpzc{P}$}^*_{pop}(n)$ coincide for large $n$. The result now follows from the convergence of ${\bs \pi}^*_{pop}(n)$ to $\contour[2]{black}{$\mathpzc{P}$}^*_{pop}$.
	\vspace{0.15cm}
	
\noindent	\underline{Part 2.} By Part 1, there exists $n_1$ such that for every $n \ge n_1$
	${\bs \pi}^*_{pop}(n)=\contour[2]{black}{$\mathpzc{P}$}^*_{pop}$, and both are singletons. 
	
	 For $n>n_1$ the set $\contour[2]{black}{$\mathpzc{P}$}^*_{pop}$ is a singleton, and we now show that for $n$ sufficiently large
	\begin{equation}\label{eq:pi_pop}
		{P}\left(   \contour[2]{black}{$\mathpzc{P}$}^*_{pop} \notin {\bs \pi}^*(n)\right)\le \frac{K_C}{\cal J}.
	\end{equation}
	We have that
	\[
	{P}\left(   \contour[2]{black}{$\mathpzc{P}$}^*_{pop} \notin {\bs \pi}^*(n)\right) \le 
	\sum_{ {\mathpzc P} \ne \contour[2]{black}{$\mathpzc{P}$}^*_{pop}} {P}\left( \{ {\mathpzc P} \in {\bs \pi}^*(n)\} \cap  \{\contour[2]{black}{$\mathpzc{P}$}^*_{pop} \notin {\bs \pi}^*(n) \} \right).
	\]
	The event $\{ {\mathpzc P} \in {\bs \pi}^*(n) \}$ implies that ${\textbf{AR}}(n,{\mathpzc P}) < {\textbf{AR}}(n,{\mathpzc Q})$ for every ${\mathpzc Q} \notin {\bs \pi}^*(n)$. In particular,
	\begin{equation}\label{eq:pi_star=}
	  {P}\left( \{ {\mathpzc P} \in {\bs \pi}^*(n)\} \cap  \{\contour[2]{black}{$\mathpzc{P}$}^*_{pop} \notin {\bs \pi}^*(n) \}\right) \le P \left({\textbf{AR}}(n,{\mathpzc P}) < {\textbf{AR}}(n,{\contour[2]{black}{$\mathpzc{P}$}^*_{pop}})\right).
	\end{equation}
	We consider now two cases for ${\mathpzc P}$: ${\mathpzc P}\in \bs{\cal M}_{pop}$ and 
	${\mathpzc P}\notin \bs{\cal M}_{pop}$. Starting with the former case, since both ${\mathpzc P}$ and $\contour[2]{black}{$\mathpzc{P}$}^*_{pop}$ are in $\bs{\cal M}_{pop}$, by the argument ensuing \eqref{eq:convex1}, $
	\int P_{G_\theta}\big({\bf X}^{({\mathpzc P})'}{{\bs \beta}^{({\mathpzc P})}}={\bf X}^{(\contour[2]{black}{$\mathpzc{P}$}^*_{pop})'}{{\bs \beta}^{(\contour[2]{black}{$\mathpzc{P}$}^*_{pop})}}\big){\mathscr P}(d\theta)=1,
	$ and therefore for almost every $\theta$,
	$
	P_{G_\theta}({\bf X}^{({\mathpzc P})'}{{\bs \beta}^{({\mathpzc P})}}={\bf X}^{(\contour[2]{black}{$\mathpzc{P}$}^*_{pop})'}{{\bs \beta}^{(\contour[2]{black}{$\mathpzc{P}$}^*_{pop})}})=1;
	$
	hence,
	\begin{equation*}%\label{eq:G_j}
	\sum_{j=1}^{\cal J} E_{G_{j}}(Y-{\bf X}^{(\mathpzc P)'} {\bs \beta}^{(\mathpzc P)}_{j})^2
	=\sum_{j=1}^{\cal J} E_{G_{j}}(Y-{\bf X}^{(\contour[2]{black}{$\mathpzc{P}$}^*_{pop})'} {\bs \beta}^{(\contour[2]{black}{$\mathpzc{P}$}^*_{pop})}_{j})^2,
	\end{equation*}
	with probability 1. By the definition of ${\textbf{AR}}(n,{\mathpzc P})$,
	\[
	{\textbf{AR}}(n,{\contour[2]{black}{$\mathpzc{P}$}^*_{pop}})-{\textbf{AR}}(n,{\mathpzc P})
	=\frac{\frac{1}{\cal J}\sum_{j=1}^{\cal J} tr({\mathbb{V}}^{(\contour[2]{black}{$\mathpzc{P}$}^*_{pop})}_j)- \frac{1}{\cal J}\sum_{j=1}^{\cal J} tr({\mathbb{V}}^{({\mathpzc P})}_j) }{n}
	\] 
	Therefore, going back to \eqref{eq:pi_star=}, we have
	\[
	P \Big({\textbf{AR}}(n,{\mathpzc P}) < {\textbf{AR}}(n,{\contour[2]{black}{$\mathpzc{P}$}^*_{pop}})\Big)=P \Big(\frac{1}{\cal J}\sum_{j=1}^{\cal J} tr({\mathbb{V}}^{(\contour[2]{black}{$\mathpzc{P}$}^*_{pop})}_j)- \frac{1}{\cal J}\sum_{j=1}^{\cal J} tr({\mathbb{V}}^{({\mathpzc P})}_j) >0 \Big).
	\]
	 By Part 1,
	\begin{equation*}%\label{eq:P_tr}
	\int tr\left(\mathbb{V}_\theta^{(\contour[2]{black}{$\mathpzc{P}$}^*_{pop})}\right) {\mathscr P}(d\theta) - \int tr\left(\mathbb{V}_\theta^{({\mathpzc P})}\right) {\mathscr P}(d\theta) \le -\ve,
	\end{equation*}
	where $\ve$ is the difference between $\int tr\big(\mathbb{V}_\theta^{(\contour[2]{black}{$\mathpzc{P}$}^*_{pop})}\big) {\mathscr P}(d\theta)$ and the second best.
	Therefore,  $E\big(\frac{1}{\cal J}\sum_{j=1}^{\cal J} tr({\mathbb{V}}^{(\contour[2]{black}{$\mathpzc{P}$}^*_{pop})}_j)- \frac{1}{\cal J}\sum_{j=1}^{\cal J} tr({\mathbb{V}}^{({\mathpzc P})}_j)\big)\le -\ve$; also, $Var(tr(V_\theta))$ is bounded (by a constant that depends on $C$). Chebyshev's inequality implies that
	\[
	P \Big(\frac{1}{\cal J}\sum_{j=1}^{\cal J} tr({\mathbb{V}}^{(\contour[2]{black}{$\mathpzc{P}$}^*_{pop})}_j)- \frac{1}{\cal J}\sum_{j=1}^{\cal J} tr({\mathbb{V}}^{({\mathpzc P})}_j) >0 \Big) \le K_C/{\cal J},
	\]   
	and therefore, $ {P}\Big( \{ {\mathpzc P} \in {\bs \pi}^*(n)\} \cap  \{\contour[2]{black}{$\mathpzc{P}$}^*_{pop} \notin {\bs \pi}^*(n) \}   \Big) \le K_C/{\cal J}$.

	Next  consider the case ${\mathpzc P}\notin \bs{\cal M}_{pop}$.
	By definition, there exists $\varepsilon>0$ such that  for any ${{\mathpzc P}\notin \bs{\cal M}_{pop}}$
	\begin{equation*}\label{eq:n2}
	 \int  E_{G_{\theta}}(Y-{\bf X}^{(\mathpzc P)'} {\bs \beta}^{(\mathpzc P)}_{\theta})^2 {\mathscr P}(d\theta)-\int  E_{G_{\theta}}(Y-{\bf X}^{(\contour[2]{black}{$\mathpzc{P}$}^*_{pop})'} {\bs \beta}^{(\contour[2]{black}{$\mathpzc{P}$}^*_{pop})}_{\theta})^2 {\mathscr P}(d\theta) > \varepsilon.
	\end{equation*} 
It is easy to see that for $n_2$ large enough this implies	
	\[
	E \big({\textbf{AR}}(n,{\contour[2]{black}{$\mathpzc{P}$}^*_{pop}}) - {\textbf{AR}}(n,{\mathpzc P})\big) < -\ve/2 ~~~~~ \forall n \ge n_2.
	\]
	By an argument as above ${P}\Big( \{ {\mathpzc P} \in {\bs \pi}^*(n)\} \cap  \{\contour[2]{black}{$\mathpzc{P}$}^*_{pop} \notin {\bs \pi}^*(n) \}   \Big) \le K_C/{\cal J}$. Since the number of models is finite, \eqref{eq:pi_pop} follows.
	
	Now, for fixed $n$ that satisfies $n <n_0:=\max\{n_1,n_2\}$ again a similar argument shows that for any ${\mathpzc P} \in {\bs \pi}^*(n)$ ,
	\[
	{P}\Big( {\mathpzc P} \notin {\bs \pi}_{pop}^*(n)\Big) \le \frac{K_C(n)}{\cal J},
	\]
	where $K_C(n)$ may depend on $n$ (and on $C$). Since there are only finite such $n$'s the result of Part 2 follows. \qed

	\noindent{\bf Proof of Proposition \ref{prop:popolo}}.
	The first part of Proposition \ref{prop:popolo} follows from Part 1 of Proposition \ref{prop:cons_J_fixed}, which shows that $\widehat{{\bs \pi}^*}(n,{\bf N}) \subseteq {\bs \pi}^*(n)$ with probability converging to 1, and Part 2 of Lemma \ref{lem:pop}, which shows that ${\bs \pi}^*(n) \subseteq {\bs \pi}^*_{pop}(n)$ with high probability.  
	
	The second part of Proposition \ref{prop:popolo} follows from a combination of several statements: $\widehat{\pi^*}(n,{\bf N})=\contour[2]{black}{$\mathpzc{P}$}^*(n)$ with high probability (Proposition \ref{prop:cons_J_fixed}, Part 2); $\contour[2]{black}{$\mathpzc{P}$}^*(n)={\bs \pi}^*(n)$ for large $n$ (Proposition \ref{prop:Cpj}); ${\bs \pi}^*(n) \subseteq {\bs \pi}^*_{pop}(n)$ with high probability (Lemma \ref{lem:pop} Part 2); and for large $n$, ${\bs \pi}^*_{pop}(n)$ is a singleton, and ${\bs \pi}^*_{pop}(n)=\contour[2]{black}{$\mathpzc{P}$}^*_{pop}(n)$ (Lemma \ref{lem:pop}, Part 1).  \qed

\vspace{1cm}

\newpage

		\section{Appendix B: A table of notation}\label{sec:AppB}
	
	\captionsetup{labelformat=empty} 
	
	\begin{table}[h!]
		\centering
		\begin{tabular}{cl}
			Expression & Description \\% & First place \\ 
			\hline
			$\cal J$ & Number of observed regression datasets\\ %& Page \pageref{P:Dj}\\ 
			$N_{j}$ & {Number of observations in the  the $j$th regression dataset}  \\
			$Y_{ij}$ &{The response of the $i$th observation from the $j$th regression }\\ 
			${\bf X}_{ij} \in \mathbb{R}^d$ & {The covariate vector of the $i$th observation from the $j$th regression } \\
			$({\bf X},Y)$ & A generic  observation (whose distribution is $G_j$)\\
			$D_j=\{({\bf X}_{ij},Y_{ij})\}$ & The $j$th regression dataset \\
			$G_j$ & The distribution of the $j$th regression, i.e., $\{({\bf X}_{ij},Y_{ij})\} \sim^{iid} G_j$ \\%&  Page \pageref{P:Dj}\\
			$\cal G$ & A set of distributions to which $G_j$ belongs (the cases $|{\cal G}|=1$,   \\
			 &$|{\cal G}|={\cal J}$ and ${\cal J} <|{\cal G}| \le \infty$ appear in Sections \ref{sec:predd}, \ref{subsec:mallowssev}, and \ref{sec:supersuper},\\
			& respectively) \\
			${\cal K}$ & The size of ${\cal G}$\\
			${\mathpzc P}$ & A subset of $\{1,\ldots, d\}$, used to denote a subset of covariates\\ 
			& Its size is denoted by $p$.\\  
			$R(n,{\mathpzc P})$ & The prediction error  of the linear model with covariates in ${\mathpzc P}$ with $n$\\& observations for the case $|{\cal G}|=1$;
			${\bf R}(n,{\mathpzc P})$ and ${\bf R}_{pop}(n,{\mathpzc P})$\\
			&  denote the cases of $|{\cal G}|={\cal J}$ and ${\cal J}<|{\cal G}|$, respectively\\
			$AR(n,{\mathpzc P})$ & Approximate prediction error; ${\bf AR}(n,{\mathpzc P})$ and ${\bf AR}_{pop}(n,{\mathpzc P})$ are\\ & approximations of ${\bf R}(n,{\mathpzc P})$ and ${\bf R}_{pop}(n,{\mathpzc P})$, respectively  \\
			\hline
			\multicolumn{2}{c}{In the notation below $j$ and ${}^{({\mathpzc P})}$ are sometimes suppressed}\\
			\hline 
			$\mathbb{X}_{j,N_j}^{({\mathpzc P})}$ & {The $N_j\times p$ design matrix of the $j$th regression}\\
			${\bf Y}_{j,N_j}$ &{The vector of responses for the $j$th regression}\\
			${\bs \beta}_j^{({\mathpzc P})}$ &  {Projection coefficients under $G_j$ for model ${\mathpzc P}$} \\
			${e}_j^{({\mathpzc P})}$ & The residual; ${e}_j^{({\mathpzc P})}=Y-{\bf X}_{j}^{({\mathpzc P})'}{\bs \beta}_j^{({\mathpzc P})}$; ${\bf e}_{j,N_j}$ denotes the vector\\
			& of the residuals of dimension $N_j$\\
			$\widehat{\bs \beta}_{j,n}^{({\mathpzc P})}$ & {The least squares estimate of ${\bs \beta}_j^{({\mathpzc P})}$ based on $n$ observations. }\\
			$\mathbb{Q}_j^{({\mathpzc P})}$ & $E_{G_j}({\bf X}^{({\mathpzc P})}{\bf X}^{({\mathpzc P})'})$  \\
			$\mathbb{W}_j^{({\mathpzc P})}$ & $E_{G_j}({\bf X}^{({\mathpzc P})}{\bf X}^{({\mathpzc P})'}e^2)$  \\ 
			$\mathbb{V}_j^{({\mathpzc P})}$ & $\mathbb{W}_j^{({\mathpzc P})}\{\mathbb{Q}_j^{({\mathpzc P})}\}^{-1}$ \\
			$\widehat{\mathbb{Q}}_{j,N_j}^{({\mathpzc P})}$ & The empirical estimate of $\mathbb{Q}_j^{({\mathpzc P})}$\\
			$\widehat{\mathbb{W}}_{j,N_j}^{({\mathpzc P})}$ & The empirical estimate of $\mathbb{W}_j^{({\mathpzc P})}$\\
			$\widehat{\mathbb{V}}_{j,N_j}^{({\mathpzc P})}$ & The empirical estimate of $\mathbb{V}_j^{({\mathpzc P})}$\\	                               
			${\bf U}^{({\mathpzc P})}_{j,N_j}$ & $\frac{1}{\sqrt{N_j}}\mathbb{X}^{({\mathpzc P})'}_{j,N_j} {\bf e}_{j,N_j}$ (it is not a statistic)\\
			\hline
			$C^{({\mathpzc P})}(n,{N})$ & An estimate of $AR(n,{\mathpzc P})$; ${\bf C}^{({\mathpzc P})}(n,{\bf N})$  corresponds to the case ${\cal J}>1$;\\
			& $\mathbb{C}^{({\mathpzc P})}(n,{N})$ and $\pmb{\mathbb{C}}^{({\mathpzc P})}(n,{N})$ denote a jackknife bias correction\\
			${{\mathpzc P}}^*(n)$ & $\arg \min_{\mathpzc P} R(n,{\mathpzc P})$  (the best model for $n$ observations);\\
			&  $\contour[2]{black}{$\mathpzc{P}$}^*(n)$ corresponds to the case ${\cal J}>1$\\
			$\pi^*(n)$ &  $\arg \min_{\mathpzc P} AR(n,{\mathpzc P})$ ; ${\bs \pi}^*(n)$ corresponds to the case ${\cal J}>1$\\
			${\mathpzc P}^*$ & The limit of both ${{\mathpzc P}}^*(n)$ and  $\pi^*(n)$ as $n \to \infty$;\\
			& $\contour[2]{black}{$\mathpzc{P}$}^*$ corresponds to the case ${\cal J}>1$\\
			$\widehat{\pi^*}(n,{\bf N})$ & $\arg \min_{\mathpzc P} {C}^{({\mathpzc P})}(n,{\bf N})$ ;$\widehat{{\bs \pi}^*}(n,{\bf N})$ corresponds to the case ${\cal J}>1$ \\
		\end{tabular}
	\end{table}

\end{document}